\documentclass{amsproc}
\usepackage{amsfonts}

\theoremstyle{plain}

\newtheorem{definition}{Definition}

\newtheorem{lemma}{Lemma}

\newtheorem{theorem}{Theorem}
\numberwithin{equation}{section}
\input{tcilatex}

\begin{document}
\title[Bounndedness]{On the local everywhere bounndedness of the minima of a
class of integral functionals of the Calculus of the Variations}
\author{Tiziano Granucci}
\address{MIUR, Firenze, 50100, Italy}
\email{tizianogranucci@libero.it}
\urladdr{https://www.tizianogranucci.com/}
\date{30/1/2023}
\subjclass[2000]{ 49N60, 35J50}
\keywords{Everywhere boinndedness, vectorial, minimizer, variational,
integral.}
\dedicatory{For my family: Elsa Cirri, Caterina Granucci and Delia Granucci.}
\thanks{This paper is in final form and no version of it will be submitted
for publication elsewhere. Data sharing not applicable to this article as no
datasets were generated or analysed during the current study. The author has
no conicts of interest to declare that are relevant to the content of this
article. I would like to thank all my family and friends for the support
given to me over the years: Elisa Cirri, Caterina Granucci, Delia Granucci,
Irene Granucci, Laura and Fiorenza Granucci, Massimo Masi and Monia
Randolfi. }

\begin{abstract}
In this paper we study the regularity and the boundedness of the minima of
two classes of functionals of the calculus of variations
\end{abstract}

\maketitle

\section{Introduction}

In this paper we study the local everywhere regularity of the following
integral functional of the Calculus of Variations\bigskip 
\begin{equation}
\mathcal{F}_{0}\left( u,\Omega \right) =\int\limits_{\Omega }G\left(
\left\vert u^{1}\right\vert ,...,\left\vert u^{m}\right\vert ,\left\vert
\nabla u^{1}\right\vert ,...,\left\vert \nabla u^{m}\right\vert \right) \,dx
\end{equation}%
where $\Omega $ is a subset of $%
\mathbb{R}
^{n}$ and $u\in W^{1,p}\left( \Omega ,%
\mathbb{R}
^{m}\right) $\ with $n\geq 2$ and $m\geq 1$. We will make the following
assumptions about the density function $G$.

\begin{description}
\item[H.1] $G:%
\mathbb{R}
_{+,0}^{m}\times 
\mathbb{R}
_{+,0}^{m}\rightarrow 
\mathbb{R}
$ is a strictly convex function on $%
\mathbb{R}
_{+,0}^{m}\times 
\mathbb{R}
_{+,0}^{m}$, with $%
\mathbb{R}
_{+,0}^{m}=\left[ 0,+\infty \right) \times ...\times \left[ 0,+\infty
\right) $, and a real positive constan $L_{0}>1$ exists such that 
\begin{equation}
\sum\limits_{\alpha =1}^{m}\left( \xi ^{\alpha }\right)
^{q}+\sum\limits_{\alpha =1}^{m}\left( s^{\alpha }\right) ^{q}\leq G\left(
s,\xi \right) \leq L_{0}\left[ \sum\limits_{\alpha =1}^{m}\left( \xi
^{\alpha }\right) ^{q}+\sum\limits_{\alpha =1}^{m}\left( s^{\alpha }\right)
^{q}+a\left( x\right) \right]
\end{equation}%
for every $s,\xi \in 
\mathbb{R}
_{+,0}^{m}$.

\item[H.2] $G\in C^{2}\left( 
\mathbb{R}
_{+,0}^{m}\times 
\mathbb{R}
_{+,0}^{m}\right) $ and a real positive constant $\vartheta >0$\ exists such
that%
\begin{equation}
\begin{tabular}{l}
$\sum\limits_{\alpha ,\beta =1}^{m}\partial _{s^{\beta }s^{\alpha }}G\left(
s,\xi \right) \lambda ^{\alpha }\lambda ^{\beta }+\sum\limits_{\alpha ,\beta
=1}^{m}\partial _{\xi ^{\beta }s^{\alpha }}G\left( s,\xi \right) \lambda
^{\alpha }\mu ^{\beta }$ \\ 
$+\sum\limits_{\alpha ,\beta =1}^{m}\partial _{s^{\beta }\xi ^{\alpha
}}G\left( s,\xi \right) \mu ^{\alpha }\lambda ^{\beta }+\sum\limits_{\alpha
,\beta =1}^{m}\partial _{\xi ^{\beta }\xi ^{\alpha }}G\left( s,\xi \right)
\mu ^{\alpha }\mu ^{\beta }$ \\ 
$\geq \vartheta \left( \left\vert s\right\vert ^{q-2}+\left\vert \xi
\right\vert ^{q-2}\right) \left( \left\vert \lambda \right\vert
^{2}+\left\vert \mu \right\vert ^{2}\right) $%
\end{tabular}%
\end{equation}%
for every $\left( s,\xi \right) \in 
\mathbb{R}
_{+,0}^{m}\times 
\mathbb{R}
_{+,0}^{m}$ and for every $\left( \lambda ,\mu \right) \in 
\mathbb{R}
^{m}\times 
\mathbb{R}
^{m}$. Moreover, a real positive costant $L_{1}>1$\ exists such that%
\begin{equation}
\begin{tabular}{l}
$0\leq \partial _{\xi ^{\alpha }}G(s,\xi )\leq L_{1}\left[ \left\vert
s\right\vert ^{q-1}+\left\vert \xi \right\vert ^{q-1}\right] $ \\ 
$0\leq \partial _{s^{\alpha }}G(s,\xi )\leq L_{1}\left[ \left\vert
s\right\vert ^{q-1}+\left\vert \xi \right\vert ^{q-1}\right] $ \\ 
$\left\vert \partial _{s^{\beta }s^{\alpha }}G\left( s,\xi \right)
\right\vert \leq L_{1}\left[ \left\vert s\right\vert ^{q-2}+\left\vert \xi
\right\vert ^{q-2}\right] $ \\ 
$\left\vert \partial _{s^{\beta }\xi ^{\alpha }}G\left( s,\xi \right)
\right\vert \leq L_{1}\left[ \left\vert s\right\vert ^{q-2}+\left\vert \xi
\right\vert ^{q-2}\right] $ \\ 
$\left\vert \partial _{\xi ^{\beta }\xi ^{\alpha }}G\left( s,\xi \right)
\right\vert \leq L_{1}\left[ \left\vert s\right\vert ^{q-2}+\left\vert \xi
\right\vert ^{q-2}\right] $%
\end{tabular}%
\end{equation}%
for $\alpha ,\beta =1,...,m$ and for every $\left( s,\xi \right) \in 
\mathbb{R}
_{+,0}^{m}\times 
\mathbb{R}
_{+,0}^{m}$.
\end{description}

Under these hypotheses we will obtain the following regularity theorem.

\begin{theorem}
If $u_{0}\in W^{1,q}\left( \Omega ,%
\mathbb{R}
^{m}\right) $, with $2\leq q<n$, $n\geq 2$ and $m\geq 1$, is a local
minimizer of the functional 1.1 and H.1 and H.2 hold then $\left\vert
u_{0}\right\vert \in L_{loc}^{\infty }\left( \Omega \right) $.
\end{theorem}

We will prove Theorem 1 by studying the boundary-bound behavior of the
following functional%
\begin{equation}
\mathcal{F}_{\varepsilon }\left( u,\Omega \right) =\int\limits_{\Omega
}\varepsilon \sum\limits_{\alpha =1}^{m}\left\vert \nabla u^{\alpha
}\right\vert ^{p}+G\left( \sqrt{\left\vert u^{1}\right\vert ^{2}+\varepsilon 
},...,\sqrt{\left\vert u^{m}\right\vert ^{2}+\varepsilon },\sqrt{\left\vert
\nabla u^{1}\right\vert ^{2}+\varepsilon },...,\sqrt{\left\vert \nabla
u^{m}\right\vert ^{2}+\varepsilon }\right) \,dx
\end{equation}

where $\Omega $ is a subset of $%
\mathbb{R}
^{n}$ and $u\in W^{1,p}\left( \Omega ,%
\mathbb{R}
^{m}\right) $\ with $n\geq 2$, $m\geq 1$ and $0<\varepsilon \leq 1$.

The study of this particular type of functionals was introduced and studied
in [9, 10, 28, 30, 35-39]. In particular the author in [30, 37-39] has
studied the Holder continuity of the minima of the functional $\mathcal{F}%
_{\varepsilon }\left( u,\Omega \right) $ defined by (1.5) under very general
hypotheses. In [37] the author considered the following class of functionals%
\begin{equation}
\mathcal{F}\left( u,\Omega \right) =\int\limits_{\Omega }\sum\limits_{\alpha
=1}^{m}\left\vert \nabla u^{\alpha }\right\vert ^{p}+G\left( x,u,\left\vert
\nabla u^{1}\right\vert ,...,\left\vert \nabla u^{m}\right\vert \right) \,dx
\end{equation}

taking on the following assumptions about the function $G$.

\begin{description}
\item[H.1.1] Let $\Omega $ be a bounded open subset of $%
\mathbb{R}
^{n}$ with $n\geq 2\ $and let $G:\Omega \times 
\mathbb{R}
^{m}\times 
\mathbb{R}
_{0,+}^{m}\rightarrow 
\mathbb{R}
$ be a Caratheodory function, where $%
\mathbb{R}
_{0,+}=\left[ 0,+\infty \right) $ $\ $and $%
\mathbb{R}
_{0,+}^{m}=%
\mathbb{R}
_{0,+}\times \cdots \times 
\mathbb{R}
_{0,+}$ with $m\geq 1$; we make the following growth conditions on $G$:\
there exists a constant $L>1$ such that%
\begin{equation*}
\sum\limits_{\alpha =1}^{m}\left\vert \xi ^{\alpha }\right\vert
^{q}-\sum\limits_{\alpha =1}^{m}\left\vert s^{\alpha }\right\vert
^{q}-a\left( x\right) \leq G\left( x,s^{1},...,s^{m},\left\vert \xi
^{1}\right\vert ,...,\left\vert \xi ^{m}\right\vert \right) \leq L\left[
\sum\limits_{\alpha =1}^{m}\left\vert \xi ^{\alpha }\right\vert
^{q}+\sum\limits_{\alpha =1}^{m}\left\vert s^{\alpha }\right\vert
^{q}+a\left( x\right) \right]
\end{equation*}%
for $\mathcal{L}^{n}$ a. e. $x\in \Omega $, for every $s^{\alpha }\in 
\mathbb{R}
$ and for every $\xi ^{\alpha }\in 
\mathbb{R}
$ with $\alpha =1,...,m$ and $m\geq 1$ and with $a\left( x\right) \in
L^{\sigma }\left( \Omega \right) $, $a(x)\geq 0$ for $\mathcal{L}^{n}$ a. e. 
$x\in \Omega $, $\sigma >\frac{n}{p}$, $1\leq q<\frac{p^{2}}{n}$ and $1<p<n$.
\end{description}

Assuming that the previous growth hypothesis H.1.1 holds, in [37] the author
proved the following regularity result.

\begin{theorem}
Let $\Omega $ be a bounded open subset of $%
\mathbb{R}
^{n}$ with $n\geq 2$; if $u\in W^{1,p}\left( \Omega ,%
\mathbb{R}
^{m}\right) $, with $m\geq 1$, is a local minimum of the functional 1.6 and $%
H.1.1$\ holds then $u^{\alpha }\in C_{loc}^{o,\beta _{0}}\left( \Omega
\right) $ for every $\alpha =1,...,m$, with $\beta _{0}\in \left( 0,1\right) 
$.
\end{theorem}

Hence we deduce that the minima of the functional (1.5) are locally Holder
continuous functions. In particular we are interested in the behavior at the
limit $\varepsilon \rightarrow 0$ of minima of the functional (1.5). In this
sense, the first interesting result on the local minima of the functional $%
F_{\varepsilon }\left( u,\Omega \right) $\ is the following Theorem with two
energy estimates of the Caccioppoli type.

\begin{theorem}
If $u_{\varepsilon }\in W^{1,p}\left( \Omega ,%
\mathbb{R}
^{m}\right) $, with $1<q<p<n$, $n\geq 2$ and $m\geq 1$, is a local minimizer
of the functional 1.5 and H.1 and H.2 hold then for every $\Sigma $ compact
subset of $\Omega $ there exist two positive constant $C_{1,\Sigma }$ and $%
C_{2,\Sigma ,\varepsilon }$ such that for every $x_{0}\in \Sigma $ and $%
0<\varrho <R<R_{0}=\frac{1}{4}\min \left\{ 1,dist\left( x_{0},\partial
\Sigma \right) \right\} $ the following Caccioppoli inequalities hold 
\begin{equation}
\begin{tabular}{l}
$\int\limits_{A_{k,\varrho }^{\alpha ,\varepsilon }}\varepsilon \left\vert
\nabla u_{\varepsilon }^{\alpha }\right\vert ^{p}\,+\left\vert \nabla
u_{\varepsilon }^{\alpha }\right\vert ^{q}\,dx$ \\ 
$\leq \frac{C_{1,\Sigma }}{\left( R-\varrho \right) ^{p}}\int%
\limits_{A_{k,R}^{\alpha ,\varepsilon }}\left( u_{\varepsilon }^{\alpha
}-k\right) ^{p}\,dx+C_{2,\Sigma ,\varepsilon }\left( 1+R^{-\epsilon
n}k^{p}\right) \left[ \mathcal{L}^{n}\left( A_{k,R}^{\alpha ,\varepsilon
}\right) \right] ^{1-\frac{p}{n}+\epsilon }$%
\end{tabular}%
\end{equation}%
and%
\begin{equation}
\begin{tabular}{l}
$\int\limits_{B_{k,\varrho }^{\alpha ,\varepsilon }}\varepsilon \left\vert
\nabla u_{\varepsilon }^{\alpha }\right\vert ^{p}\,+\left\vert \nabla
u_{\varepsilon }^{\alpha }\right\vert ^{q}\,dx$ \\ 
$\leq \frac{C_{1,\Sigma }}{\left( R-\varrho \right) ^{p}}\int%
\limits_{B_{k,R}^{\alpha ,\varepsilon }}\left( k-u_{\varepsilon }^{\alpha
}\right) ^{p}\,dx+C_{2,\Sigma ,\varepsilon }\left( 1+R^{-\epsilon
n}k^{p}\right) \left[ \mathcal{L}^{n}\left( B_{k,R}^{\alpha ,\varepsilon
}\right) \right] ^{1-\frac{p}{n}+\epsilon }$%
\end{tabular}%
\end{equation}%
for evrey $\alpha \in \left[ 1,...,m\right] $.
\end{theorem}

The inequalities of Caccioppoli (1.7) and (1.8) depend on the norm $%
\left\Vert u_{\varepsilon }\right\Vert _{W^{1,p}\left( \Sigma ,%
\mathbb{R}
^{m}\right) }$, where $\Sigma \subset \subset \Omega $ is a compact subset.
For us, it will be essential to get a uniform estimate on $\varepsilon \in
\left( 0,1\right) $ of the norm $\left\Vert u_{\varepsilon }\right\Vert
_{W^{1,p}\left( \Sigma ,%
\mathbb{R}
^{m}\right) }$. The following theorem of greater regularity of the gradient
of functions $u_{\varepsilon }$ is the first step to obtain this uniform
estimate.

\begin{theorem}
If $u_{\varepsilon }\in W^{1,p}\left( \Omega ,%
\mathbb{R}
^{m}\right) $, with $2\leq q<p<n$, $n\geq 2$ and $m\geq 1$, is a local
minimizer of the functional 1.5 and H.1 and H.2 hold then $u_{\varepsilon
}\in W_{loc}^{2,2}\left( \Omega ,%
\mathbb{R}
^{m}\right) $ and there exists a positive constant $C_{p,q}$ such that 
\begin{equation*}
\int\limits_{B_{s}\left( x_{0}\right) }\left[ \left\vert \nabla
u_{\varepsilon }\left( x\right) \right\vert ^{2}+\varepsilon \right] ^{\frac{%
q-2}{2}}\left\vert Hu_{\varepsilon }\right\vert ^{2}dx\leq \frac{C_{p,q}}{%
\left( t-s\right) ^{2}}\int\limits_{B_{t}\left( x_{0}\right) }\varepsilon
\left\vert \nabla u_{\varepsilon }\right\vert ^{p}+\left\vert u_{\varepsilon
}\right\vert ^{q}+\left\vert \nabla u_{\varepsilon }\right\vert ^{q}+1\,\,dx
\end{equation*}%
for every $0<s<t<R_{0}<\min \left\{ 1,dist\left( x_{0},\partial \Omega
\right) \right\} $.
\end{theorem}

From the previous theorem we deduce the following result which gives us a
uniform estimate on the norm $\left\Vert u_{\varepsilon }\right\Vert
_{W^{1,p}\left( \Sigma ,%
\mathbb{R}
^{m}\right) }$, where $\Sigma \subset \subset \Omega $ is a compact subset.

\begin{theorem}
\bigskip If $u_{\varepsilon }\in W^{1,p}\left( \Omega ,%
\mathbb{R}
^{m}\right) $, with $2\leq q<p<n$, $n\geq 2$ and $m\geq 1$, is a local
minimizer of the functional 1.5 and H.1 and H.2 hold then $\left\vert \nabla
u_{\varepsilon }\right\vert \in L_{loc}^{\frac{2^{\ast }}{2}q}\left( \Omega
\right) $.
\end{theorem}

\begin{theorem}
If $u_{\varepsilon }\in W^{1,p}\left( \Omega ,%
\mathbb{R}
^{m}\right) $, with $2\leq q<p<\min \left\{ \frac{2^{\ast }}{2}q,q^{\ast
}\right\} $, $n\geq 2$ and $m\geq 1$, is a local minimizer of the functional
1.5 and H.1 and H.2 hold then for every $x_{0}\in \Omega $ and for every
radius $r$, with $0<r<\frac{1}{4}dist\left( x_{0},\partial \Omega \right) $
there exists a positive constan $C_{p,q}$\ such that 
\begin{equation}
\left\Vert u_{\varepsilon }\right\Vert _{W^{1,p}\left( B_{r}\left(
x_{0}\right) \right) }<D_{p,q}
\end{equation}%
for every $\varepsilon \in \left( 0,1\right] $.
\end{theorem}

Theorem 6 is the heart of our proof, from these results, with a procedure of
passage to the limit, the following Caccioppoli inequalities are obtained.

\begin{theorem}
If $u_{0}\in W^{1,q}\left( \Omega ,%
\mathbb{R}
^{m}\right) $, with $2\leq q<\min \left\{ \frac{2^{\ast }}{2}q,q^{\ast
}\right\} $, $n\geq 2$ and $m\geq 1$, is a local minimizer of the functional
1.1 and H.1 and H.2 hold holds then%
\begin{equation}
\begin{tabular}{l}
$\int\limits_{A_{k,t\varrho y_{0}}^{\alpha ,0}}\left\vert \nabla
u_{0}^{\alpha }\right\vert ^{q}\,dx$ \\ 
$\leq \frac{C_{1,B_{r}\left( x_{0}\right) }}{\left( R-\varrho \right) ^{p}}%
\int\limits_{A_{k,R,y_{0}}^{\alpha ,0}}\left( u_{0}^{\alpha }-k\right)
^{p}\,dx+C_{2,B_{r}\left( x_{0}\right) }\left( 1+R^{-\epsilon n}k^{p}\right) %
\left[ \mathcal{L}^{n}\left( A_{k,R,y_{0}}^{\alpha ,0}\right) \right] ^{1-%
\frac{p}{n}+\epsilon }$%
\end{tabular}%
\end{equation}%
and%
\begin{equation}
\begin{tabular}{l}
$\int\limits_{B_{k,t\varrho y_{0}}^{\alpha ,0}}\left\vert \nabla
u_{0}^{\alpha }\right\vert ^{q}\,dx$ \\ 
$\leq \frac{C_{1,B_{r}\left( x_{0}\right) }}{\left( R-\varrho \right) ^{p}}%
\int\limits_{B_{k,R,y_{0}}^{\alpha ,0}}\left( k-u_{0}^{\alpha }\right)
^{p}\,dx+C_{2,B_{r}\left( x_{0}\right) }\left( 1+R^{-\epsilon n}k^{p}\right) %
\left[ \mathcal{L}^{n}\left( B_{k,R,y_{0}}^{\alpha ,0}\right) \right] ^{1-%
\frac{p}{n}+\epsilon }$%
\end{tabular}%
\end{equation}%
for every $\alpha =1,...,m$.
\end{theorem}

The Caccioppoli inequalities (1.10) and (1.11) are typical of scalar
variational problems with growths of the type q-p; in particular we will
refer to the ideas presented in [14, 29, 32-34] to obtain from the
inequalities (1.10) and (1.11) the boundedness theorem Theorem 1.

The results given in Theorem 1 and Theorem 2 are not obvious, as shown by
the examples developed by E. De Giorgi in [16] and by E. Giusti and M.
Miranda in [26]. In particular, in [16], E. De Giorgi proved that the
function $u:B_{1}\left( 0\right) \subset 
\mathbb{R}
^{n}\rightarrow 
\mathbb{R}
^{n}$ defined by%
\begin{equation*}
u^{\alpha }\left( x\right) =x_{\alpha }\left\vert x\right\vert ^{-\frac{n}{2}%
+\frac{n}{2\sqrt{\left( 2n-2\right) ^{2}+1}}}
\end{equation*}%
with $n>2$\ is a minimum of the functional%
\begin{equation*}
F_{E}\left( v,B_{1}\left( 0\right) \right) =\int\limits_{B_{1}\left(
0\right) }\sum\limits_{i,j=1}^{n}\sum\limits_{\alpha ,\beta =1}^{n}A_{\alpha
,\beta }^{i,j}\left( x\right) \partial _{i}v^{\alpha }\partial _{j}v^{\beta
}\,dx
\end{equation*}%
with 
\begin{equation*}
A_{\alpha ,\beta }^{i,j}\left( x\right) =\delta _{\alpha ,\beta }\delta
_{i,j}+\left[ \left( n-2\right) \delta _{\alpha ,i}+n\frac{x_{\alpha }x_{i}}{%
\left\vert x\right\vert ^{2}}\right] \left[ \left( n-2\right) \delta _{\beta
,j}+n\frac{x_{\beta }x_{j}}{\left\vert x\right\vert ^{2}}\right]
\end{equation*}%
among all functions taking the value x on the boundary of the unit ball $%
B_{1}\left( 0\right) $. While E. Giusti and M. Miranda in [26] demonstrated
that the function 
\begin{equation*}
u(x)=\frac{x}{\left\vert x\right\vert }
\end{equation*}%
minimizes the functional%
\begin{equation*}
F_{GM}\left( v,B_{1}\left( 0\right) \right) =\int\limits_{B_{1}\left(
0\right) }\sum\limits_{i,j=1}^{n}\sum\limits_{\alpha ,\beta =1}^{n}A_{\alpha
,\beta }^{i,j}\left( u\right) \partial _{i}v^{\alpha }\partial _{j}v^{\beta
}\,dx
\end{equation*}%
with%
\begin{equation*}
A_{\alpha ,\beta }^{i,j}\left( x\right) =\delta _{\alpha ,\beta }\delta
_{i,j}+\left[ \delta _{\alpha ,i}+\frac{4}{n-2}\frac{u_{\alpha }u_{i}}{%
1+\left\vert u\right\vert ^{2}}\right] \left[ \delta _{\beta ,j}+\frac{4}{n-2%
}\frac{u_{\beta }u_{j}}{1+\left\vert u\right\vert ^{2}}\right]
\end{equation*}%
among the functions $v\in W^{1,2}\left( B_{1}\left( 0\right) ;%
\mathbb{R}
^{n}\right) +x$ for $n$ sufficientely large. We observe that the functionals 
$F_{E}\left( v,B_{1}\left( 0\right) \right) $ and $F_{GM}\left(
v,B_{1}\left( 0\right) \right) $\ are convex functionolans. There are other
examples of non-bounded or non-continuous extremals for vector-type
variational problems, refer to [16, 22, 26, 27]. Starting from the works of
K. Uhlenbeck [55], P. Tolksdorf [53, 54] and E. Acerbi and N. Fusco [1]
during the eighties fundamental regularity results were obtained for
functionals of the type%
\begin{equation*}
\int\limits_{\Omega }\left\vert \nabla u\right\vert ^{p}\,dx
\end{equation*}%
The previous works cited [1, 53-55] were the starting point for the study of
the regularity of the minima of functionals of the type 
\begin{equation*}
\int\limits_{\Omega }f\left( \left\vert \nabla u\right\vert \right) \,dx
\end{equation*}%
where convexity and regularity hypotheses are assumed on the density $f$
which have been made more and more generic over time, we refer to [2, 4-7,
11-13, 17-21, 23, 40, 44, 45, 50, 51] for further details.The study of
functionals of type (1.5) starts from the results obtained in [9, 10]. In
particular in [9], Cupini, Focardi, Leonetti and Mascolo introduced the
following class of vectorial functionals%
\begin{equation*}
\int\limits_{\Omega }f\left( x,\nabla u\right) \,dx
\end{equation*}%
where $\Omega \subset 
\mathbb{R}
^{n},u:\Omega \rightarrow 
\mathbb{R}
^{m},n>1,m\geq 1$ and%
\begin{equation*}
f\left( x,\nabla u\right) =\sum\limits_{\alpha =1}^{m}F_{\alpha }\left(
x,\nabla u^{\alpha }\right) +G\left( x,\nabla u\right)
\end{equation*}%
where $F\alpha \times \ R^{n\times m}\rightarrow R$ is a Carath\'{e}odory
function satisfying the following standard growth condition%
\begin{equation*}
k_{1}\left\vert \xi ^{\alpha }\right\vert ^{p}-a\left( x\right) \leq
F_{\alpha }\left( x,\xi ^{\alpha }\right) \leq k_{2}\left\vert \xi ^{\alpha
}\right\vert ^{p}+a\left( x\right)
\end{equation*}%
for every $\xi ^{\alpha }\in 
\mathbb{R}
^{n}$ and for almost every $x\in \Omega $ ,where $k_{1}$ and $k_{2}$ are two
real positive constants, $p>1$ and $a\in L_{loc}^{\sigma }\left( \Omega
\right) $ is a non negative function. In [9],Cupini, Focardi, Leonetti and
Mascolo analyze two different types of hypotheses on the $G$ function. They
started by assuming that $G:\Omega \times \ R^{n\times m}\rightarrow R$ is a
Carath\'{e}odory rank one convex function satisfying the following growth
condition%
\begin{equation*}
|G(x,\xi )|\leq k_{3}|\xi |^{q}+b(x)
\end{equation*}%
for every $\xi \in 
\mathbb{R}
^{n\times m}$, for almost every $x\in \Omega $, here $k_{3}$ is a real
positive constant, $1\leq q<p$ and $b\in L_{loc}^{\sigma }\left( \Omega
\right) $ is a nonnegative function. Moreover Cupini, Focardi, Leonetti and
Mascolo in [9] study the case where $n\geq m\geq 3$, and $G:\Omega \times \ 
\mathbb{R}
^{n\times m}\rightarrow R$ is a Carath\'{e}odory function defined as%
\begin{equation*}
G(x,\xi )=\sum\limits_{\alpha =1}^{m}G_{\alpha }\left( x,\left( adj_{m-1}\xi
\right) ^{\alpha }\right)
\end{equation*}%
here $G\alpha :\Omega \times \ 
\mathbb{R}
^{\frac{m!}{n!\left( n-m\right) !}}\rightarrow 
\mathbb{R}
$ is a Carath\'{e}odory convex function satisfying the following growth
conditions%
\begin{equation*}
0\leq G_{\alpha }\left( x,\left( adj_{m-1}\xi \right) ^{\alpha }\right) \leq
k_{4}|\left( adj_{m-1}\xi \right) ^{\alpha }|^{r}+b(x)
\end{equation*}%
for every $\xi \in 
\mathbb{R}
^{n\times m}$, for almost every $x\in \Omega $, here $k_{3}$ is a real
positive constant, $1\leq r<p$ and $b\in L_{loc}^{\sigma }\left( \Omega
\right) $ is a non negative function. In both cases, by imposing appropriate
hypotheses on the parameters $q$ and $r$ , Cupini, Focardi, Leonetti and
Mascolo proved that the localminimizers of the vectorial functional (1.3)
are locally h\"{o}lder continuous functions. In [28] the author with M.
Randolfi demonstrated a regularity result for the minima of vector
functionals with anisotropic growths of the following type%
\begin{equation*}
\int\limits_{\Omega }\sum\limits_{\alpha =1}^{m}F_{\alpha }\left( x,\nabla
u^{\alpha }\right) +G\left( x,\nabla u\right) \,dx
\end{equation*}%
with%
\begin{equation*}
\sum\limits_{\alpha =1}^{m}\Phi _{i,\alpha }\left( \left\vert \xi
_{i}^{\alpha }\right\vert \right) \leq F_{\alpha }\left( x,\xi ^{\alpha
}\right) \leq L\left[ \bar{B}_{\alpha }^{\beta _{\alpha }}\left( \left\vert
\xi ^{\alpha }\right\vert \right) +a\left( x\right) \right]
\end{equation*}%
where $\Phi _{i,\alpha }$ are N functions belonging to the class $\triangle
_{2}^{m_{\alpha }}\cap \nabla _{2}^{r_{\alpha }}$, $\bar{B}_{\alpha }$\ is
the Sobolev function associated with $\Phi _{i,\alpha }$'s, $\beta _{\alpha
}\in \left( 0,1\right] $ and $a\in L_{loc}^{\sigma }\left( \Omega \right) $
is a non negative function; oppure with%
\begin{equation*}
\sum\limits_{\alpha =1}^{m}\Phi _{i,\alpha }\left( \left\vert \xi
_{i}^{\alpha }\right\vert \right) -a\left( x\right) \leq F_{\alpha }\left(
x,\xi ^{\alpha }\right) \leq L_{1}\left[ \sum\limits_{\alpha =1}^{m}\Phi
_{i,\alpha }\left( \left\vert \xi _{i}^{\alpha }\right\vert \right) +a\left(
x\right) \right]
\end{equation*}%
where $\Phi _{i,\alpha }$ are N-functions belonging to the class $\triangle
_{2}^{m_{\alpha }}\cap \nabla _{2}^{r_{\alpha }}$ and $a\in L_{loc}^{\sigma
}\left( \Omega \right) $ is a non negative function, moreover, appropriate
hypotheses are made on the density $G$, for more details we refer to [28],
in particular, using the techniques presented in [28, 29], the authors have
shown that the minima of the functional (1.3) are locally bounded functions
in the first case and locally h\"{o}lder continuous in the second case, we
refer to [28] for more details. In [30, 53-39] the author generalized the
previous results in the case of type functionals%
\begin{equation*}
\int\limits_{\Omega }f\left( x,u,\nabla u\right) \,dx
\end{equation*}%
where $\Omega \subset 
\mathbb{R}
^{n},u:\Omega \rightarrow 
\mathbb{R}
^{m},n>1,m\geq 1$ and%
\begin{equation*}
f\left( x,u,\nabla u\right) =\sum\limits_{\alpha =1}^{m}F_{\alpha }\left(
x,u,\nabla u^{\alpha }\right) +G\left( x,u,\nabla u\right)
\end{equation*}%
with appropriate hypotheses on density $G\left( x,s,\xi \right) $; in
particular, for the following article, the results given in [37] are of
particular relevance. The proof of Theorem 1 given in this article uses,
mixing them, various techniques, the results given in [9, 10, 26, 30] are
rephrased to obtain Caccioppoli estimates for the minima of the functional
(1.5) then we proceed with techniques similar to those presented in [4, 5,
6] to obtain uniform estimates on Sobolev norms of minima of functionals of
the type (1.5); finally, by properly rephrasing the results given in [14,
29, 32-34] the proof of Theorem1 is obtained. The author hopes, in some
forthcoming works, both to develop the case 1\TEXTsymbol{<}q\TEXTsymbol{<}2,
and to obtain the continuity of the minima of the functional (1.1). It would
also be interesting to understand whether analogous results are also valid
for densities with different growths, such as non-standard or anisotropic.

\section{\protect\bigskip Preliminary results}

In this section, for completeness, we will introduce some well-known lemmas
and some definitions concerning Sobolev spaces.

\subsection{Lemmata}

\begin{lemma}[Young Inequality]
Let $\varepsilon >0$, $a,b>0$ and $1<p,q<+\infty $ with $\frac{1}{p}+\frac{1%
}{q}=1$\ then it follows%
\begin{equation}  \label{4.1}
ab\leq \varepsilon \frac{a^{p}}{p}+\frac{b^{q}}{\varepsilon ^{\frac{q}{p}}q}
\end{equation}
\end{lemma}

\begin{lemma}[H\"{o}lder Inequality]
Assume $1\leq p,q\leq +\infty $ with $\frac{1}{p}+\frac{1}{q}=1$\ then if $%
u\in L^{p}\left( \Omega \right) $\ and $v\in L^{p}\left( \Omega \right) $\
it follows%
\begin{equation}  \label{4.2}
\int\limits_{\Omega }\left\vert uv\right\vert \,dx\leq \left(
\int\limits_{\Omega }\left\vert u\right\vert ^{p}\,dx\right) ^{\frac{1}{p}%
}\left( \int\limits_{\Omega }\left\vert v\right\vert ^{q}\,dx\right) ^{\frac{%
1}{q}}
\end{equation}
\end{lemma}

\begin{lemma}
\label{lemma3} Let $Z\left( t\right) $ be a nonnegative and bounded function
on the set $\left[ \varrho ,R\right] $; if for every $\varrho \leq t<s\leq R$
we get%
\begin{equation*}
Z\left( t\right) \leq \theta Z\left( s\right) +\frac{A}{\left( s-t\right)
^{\lambda }}+\frac{B}{\left( s-t\right) ^{\mu }}+C
\end{equation*}%
where $A,B,C\geq 0$, $\lambda >\mu >0$ and $0\leq \theta <1$ then it follows%
\begin{equation*}
Z\left( \varrho \right) \leq C\left( \theta ,\lambda \right) \left( \frac{A}{%
\left( R-\varrho \right) ^{\lambda }}+\frac{B}{\left( R-\varrho \right)
^{\mu }}+C\right)
\end{equation*}%
where $C\left( \theta ,\lambda \right) >0$ is a real constant depending only
on $\theta $ and $\lambda $.
\end{lemma}

\begin{lemma}
Let $\omega >0$ and $\left\{ x_{i}\right\} _{i\in 
\mathbb{N}
}$ be a sequence of positive real numbers, such that%
\begin{equation*}
x_{i+1}\leq \Gamma D^{i}x_{i+1}^{1+\omega }
\end{equation*}%
with $\Gamma >0$ and $D>1$. If $x_{0}\leq \Gamma ^{-\frac{1}{\omega }}D^{-%
\frac{1}{\omega ^{2}}}$ then 
\begin{equation*}
x_{i}\leq D^{-\frac{i}{\omega }}x_{0}
\end{equation*}%
and in particular we have%
\begin{equation*}
\lim_{i\rightarrow +\infty }x_{i}=0\text{.}
\end{equation*}
\end{lemma}

The previous lemmas are known to most readers, for more information refer to
[3, 8, 27].

\subsection{Sobolev Spaces}

For completeness we remember that if $\Omega $ is a open subset of $%
\mathbb{R}
^{N}$ and $u$\ is a Lebesgue measurable function then $L^{p}\left( \Omega
\right) $ is the set of the class of the Lebesgue measurable function such
that $\int\limits_{\Omega }\left\vert u\right\vert ^{p}\,dx<+\infty $ and $%
W^{1,p}\left( \Omega \right) $\ is the set of the function $u\in L^{p}\left(
\Omega \right) $ such that its waek derivate $\partial _{i}u\in L^{p}\left(
\Omega \right) $. The spaces $L^{p}\left( \Omega \right) $ and $%
W^{1,p}\left( \Omega \right) $ are Banach spaces with the respective norms 
\begin{equation}
\left\Vert u\right\Vert _{L^{p}\left( \Omega \right) }=\left(
\int\limits_{\Omega }\left\vert u\right\vert ^{p}\,dx\right) ^{\frac{1}{p}}
\end{equation}%
and%
\begin{equation}
\left\Vert u\right\Vert _{W^{1,p}\left( \Omega \right) }=\left\Vert
u\right\Vert _{L^{p}\left( \Omega \right) }+\sum\limits_{i=1}^{N}\left\Vert
\partial _{i}u\right\Vert _{L^{p}\left( \Omega \right) }
\end{equation}%
We say that the function $u:\Omega \subset 
\mathbb{R}
^{N}\rightarrow 
\mathbb{R}
^{n}$ belong \ \ \ in $W^{1,p}\left( \Omega ,%
\mathbb{R}
^{n}\right) $ if $u^{\alpha }\in W^{1,p}\left( \Omega \right) $ for every $%
\alpha =1,...,n$, where $u^{\alpha }$ is the $\alpha $ component of the
vector-valued function $u$; we end by remembering that $W^{1,p}\left( \Omega
,%
\mathbb{R}
^{n}\right) $ is a Banach space with the norm%
\begin{equation}
\left\Vert u\right\Vert _{W^{1,p}\left( \Omega ,%
\mathbb{R}
^{n}\right) }=\sum\limits_{\alpha =1}^{n}\left\Vert u^{\alpha }\right\Vert
_{W^{1,p}\left( \Omega \right) }
\end{equation}

In particular, the following basic results hold.

\begin{theorem}[Sobolev Inequality]
Let $\Omega $ be a open subset of $%
\mathbb{R}
^{N}$ if $u\in W_{0}^{1,p}\left( \Omega \right) $ with $1\leq p<N$ there
exists a real positive constant $C_{SN}$, depending only on $p$ and $N$,
such that%
\begin{equation}
\left\Vert u\right\Vert _{L^{p^{\ast }}\left( \Omega \right) }\leq
C_{SN}\left\Vert \nabla u\right\Vert _{L^{p}\left( \Omega \right) }
\end{equation}%
where $p^{\ast }=\frac{Np}{N-p}$.
\end{theorem}

\begin{theorem}
(Rellich-Sobolev Immersion Theorem) Let $\Omega $ be a open bounded subset
of $%
\mathbb{R}
^{N}$ with lipschitz boundary then if $u\in W^{1,p}\left( \Omega \right) $
with $1\leq p<N$ there exists a real positive constant $C_{IS}$, depending
only on $p$ and $N$, such that%
\begin{equation}
\left\Vert u\right\Vert _{L^{p^{\ast }}\left( \Omega \right) }\leq
C_{IS}\left\Vert u\right\Vert _{W^{1,p}\left( \Omega \right) }
\end{equation}%
where $p^{\ast }=\frac{Np}{N-p}$.
\end{theorem}

Now we report other known results which give an equivalent characterization
of Sobolev spaces.

\begin{definition}
Let $f\left( x\right) $ be a function defined in an open set $\Omega \subset 
\mathbb{R}
^{N}$ and let $h$ be a real number. We call the difference quotien of $f$
with respect to $x_{s}$ the function%
\begin{equation}
\triangle _{s,h}f\left( x\right) =\frac{f\left( x+he_{s}\right) -f\left(
x\right) }{h}
\end{equation}%
where $e_{s}$ denotes the direction of the $x_{s}$ axis.
\end{definition}

The function $\triangle _{s,h}f$ is defined in the set%
\begin{equation}
\triangle _{s,h}\Omega =\left\{ x\in \Omega :x+he_{s}\in \Omega \right\}
\end{equation}%
hence in the set%
\begin{equation}
\Omega _{\left\vert h\right\vert }=\left\{ x\in \Omega :dist\left(
x,\partial \Omega \right) >\left\vert h\right\vert \right\} .
\end{equation}

\begin{lemma}
If $f\in W^{1,p}\left( \Omega \right) $ then $\triangle _{s,h}f\in
W^{1,p}\left( \Omega _{\left\vert h\right\vert }\right) $ and 
\begin{equation}
\partial _{i}\left( \triangle _{s,h}f\right) =\triangle _{s,h}\left(
\partial _{i}f\right)
\end{equation}%
for every $i=1,...,N$.
\end{lemma}

\begin{lemma}
Let $f\in L^{p}\left( \Omega \right) $ and $g\in L^{\frac{p}{p-1}}\left(
\Omega \right) $ then if at least one of the functions $f$ or $g$ has
support contained in $\Omega _{\left\vert h\right\vert }$ then%
\begin{equation}
\int\limits_{\Omega }f\,\triangle _{s,h}g\,dx=-\int\limits_{\Omega
}g\,\triangle _{s,-h}f\,dx.
\end{equation}
\end{lemma}

\begin{lemma}
Let $f\in L^{p}\left( \Omega \right) $, $1<p<\infty $, and assume that there
exists a positive constant $K$ such that for $\left\vert h\right\vert $
small enough we have%
\begin{equation}
\int\limits_{\Omega _{\left\vert h\right\vert }}\left\vert \triangle
_{s,h}f\right\vert ^{p}\,dx\leq K^{p}.
\end{equation}%
Then, $\partial _{s}f\in L^{p}\left( \Omega \right) $, and 
\begin{equation}
\left\Vert \partial _{s}f\right\Vert _{L^{p}\left( \Omega \right) }\leq K.
\end{equation}%
Moreover, when $h\rightarrow 0$, $\triangle _{s,h}f\rightarrow \partial
_{s}f $ in $L^{p}\left( \Omega \right) $.
\end{lemma}

\begin{lemma}
Let $\xi $ and $\pi $ be two vectors in $%
\mathbb{R}
^{N}$ and let%
\begin{equation}
Z\left( t\right) =\left( 1+\left\vert \left( 1-t\right) \xi +t\pi
\right\vert ^{2}\right) ^{\frac{1}{2}}
\end{equation}%
for every $s>-1$ and $r>0$ there exists two constants $c_{1}\left(
r,s\right) $ and $c_{2}\left( r,s\right) $ such that%
\begin{equation}
c_{1}\left( r,s\right) \left( 1+\left\vert \xi \right\vert ^{2}+\left\vert
\pi \right\vert ^{2}\right) ^{\frac{s}{2}}\leq \int\limits_{0}^{1}\left(
1-t\right) ^{r}Z^{r}\left( t\right) \,dt\leq c_{2}\left( r,s\right) \left(
1+\left\vert \xi \right\vert ^{2}+\left\vert \pi \right\vert ^{2}\right) ^{%
\frac{s}{2}}
\end{equation}
\end{lemma}

To prove our regularity theorem, Theorem 1, we need to introduce some
fundamental regularity results introduced in the 70s that generalize the
famous Theorem of De Giorgi-Nash-Moser, refer to [15, 47, 48, 52], in
particular we refer to the results of [25, 41, 42] as presented in [25] and
[27].

\begin{definition}
Let $\Omega \subset 
\mathbb{R}
^{N}$ be a bounded open set and $v:\Omega \rightarrow 
\mathbb{R}
$, we say that $v\in W_{loc}^{1,p}\left( \Omega \right) $ belong to the De
Giorgi class $DG^{+}\left( \Omega ,p,\lambda ,\lambda _{\ast },\chi
,\varepsilon ,R_{0},k_{0}\right) $ with $p>1$, $\lambda >0$, $\lambda _{\ast
}>0$, $\chi >0$, $\varepsilon >0$, $R_{0}>0$ and $k_{0}\geq 0$ if%
\begin{equation*}
\int\limits_{A_{k,\varrho }}\left\vert \nabla v\right\vert ^{p}\,dx\leq 
\frac{\lambda }{\left( R-\varrho \right) ^{p}}\int\limits_{A_{k,R}}\left(
v-k\right) ^{p}\,dx+\lambda _{\ast }\left( \chi ^{p}+k^{p}R^{-N\varepsilon
}\right) \left\vert A_{k,R}\right\vert ^{1-\frac{p}{N}+\varepsilon }
\end{equation*}%
for all $k\geq k_{0}\geq 0$ and for all pair of balls $B_{\varrho }\left(
x_{0}\right) \subset B_{R}\left( x_{0}\right) \subset \subset \Omega $ with $%
0<\varrho <R<R_{0}$ and $A_{k,s}=B_{s}\left( x_{0}\right) \cap \left\{
v>k\right\} $ with $s>0$.
\end{definition}

\begin{definition}
Let $\Omega \subset 
\mathbb{R}
^{N}$ be a bounded open set and $v:\Omega \rightarrow 
\mathbb{R}
$, we say that $v\in W_{loc}^{1,p}\left( \Omega \right) $ belong to the De
Giorgi class $DG^{-}\left( \Omega ,p,\lambda ,\lambda _{\ast },\chi
,\varepsilon ,R_{0},k_{0}\right) $ with $p>1$, $\lambda >0$, $\lambda _{\ast
}>0$, $\chi >0$ and $k_{0}\geq 0$ if%
\begin{equation*}
\int\limits_{B_{k,\varrho }}\left\vert \nabla v\right\vert ^{p}\,dx\leq 
\frac{\lambda }{\left( R-\varrho \right) ^{p}}\int\limits_{B_{k,R}}\left(
k-v\right) ^{p}\,dx+\lambda _{\ast }\left( \chi ^{p}+\left\vert k\right\vert
^{p}R^{-N\varepsilon }\right) \left\vert B_{k,R}\right\vert ^{1-\frac{p}{N}%
+\varepsilon }
\end{equation*}%
for all $k\leq -k_{0}\leq 0$ and for all pair of balls $B_{\varrho }\left(
x_{0}\right) \subset B_{R}\left( x_{0}\right) \subset \subset \Omega $ with $%
0<\varrho <R<R_{0}$ and $B_{k,s}=B_{s}\left( x_{0}\right) \cap \left\{
v<k\right\} $ with $s>0$.
\end{definition}

\begin{definition}
We set $DG\left( \Omega ,p,\lambda ,\lambda _{\ast },\chi ,\varepsilon
,R_{0},k_{0}\right) =DG^{+}\left( \Omega ,p,\lambda ,\lambda _{\ast },\chi
,\varepsilon ,R_{0},k_{0}\right) \cap DG^{-}\left( \Omega ,p,\lambda
,\lambda _{\ast },\chi ,\varepsilon ,R_{0},k_{0}\right) $.
\end{definition}

\begin{theorem}
Let $v\in DG\left( \Omega ,p,\lambda ,\lambda _{\ast },\chi ,\varepsilon
,R_{0},k_{0}\right) $ and $\tau \in (0,1)$, then there exists a constant $%
C>1 $ depending only upon the data and not-dependent on $v$ and $x_{0}\in
\Omega $ such that for every pair of balls $B_{\tau \varrho }\left(
x_{0}\right) \subset B_{\varrho }\left( x_{0}\right) \subset \subset \Omega $
with $0<\varrho <R_{0}$ 
\begin{equation}
\left\Vert v\right\Vert _{L^{\infty }\left( B_{\tau \varrho }\left(
x_{0}\right) \right) }\leq \max \left\{ \lambda _{\ast }\varrho ^{\frac{%
N\varepsilon }{p}};\frac{C}{\left( 1-\tau \right) ^{\frac{N}{p}}}\left[ 
\frac{1}{\left\vert B_{\varrho }\left( x_{0}\right) \right\vert }%
\int\limits_{B_{\varrho }\left( x_{0}\right) }\left\vert v\right\vert
^{p}\,dx\right] ^{\frac{1}{p}}\right\} .
\end{equation}
\end{theorem}

For more details on De Giorgi's classes and for the proof of the Theorem 10
refer to [15, 25, 27, 41, 42, 47, 48, 52] while for more details on Sobolev
spaces, Theorem 8 and Theorem 9, Lemma 7 and Lemma 8 refer to [3, 8, 27].

\section{Proof of Teorem 3}

In this section we will give the proof of Theorem . \ Let $\Sigma \subset
\subset \Omega $ be a compact subset of $\Omega $, let $R_{0}$ be a positive
real number, let us consider $x_{0}\in \Sigma $, $0<\varrho \leq t<s\leq
R<\min \left( R_{0},1,\frac{dist\left( x_{0},\partial \Sigma \right) }{2}%
\right) $, $k>k_{0}\geq 0$\ and $\eta \in C_{c}^{\infty }\left( B_{s}\left(
x_{0}\right) \right) $ such that $0\leq \eta \leq 1$ on $B_{s}\left(
x_{0}\right) $, $\eta =1$ in $B_{t}\left( x_{0}\right) $ and $\left\vert
\nabla \eta \right\vert \leq \frac{2}{s-t}$ on $B_{s}\left( x_{0}\right) $,
moreover we we can observe that supp$\left( \nabla \eta \right) \subset
B_{s}\left( x_{0}\right) \backslash B_{t}\left( x_{0}\right) $. We define $%
\varphi _{\varepsilon }=-\eta ^{\lambda }w_{\varepsilon ,+}$ where%
\begin{equation*}
w_{\varepsilon ,+}=\left( 
\begin{tabular}{l}
$\left( u_{\varepsilon }^{1}-k\right) _{+}$ \\ 
$0$ \\ 
$\vdots $ \\ 
$0$%
\end{tabular}%
\right)
\end{equation*}%
with $\left( u_{\varepsilon }^{1}-k\right) _{+}=\max \left\{ u_{\varepsilon
}^{1}-k,0\right\} $ then, since $u_{\varepsilon }$ is a local sub-minimum of 
$\mathcal{F}_{\varepsilon }\left( u,\Omega \right) $, it follows 
\begin{equation}
\mathcal{F}\left( u_{\varepsilon },B_{s}\left( x_{0}\right) \right) \leq 
\mathcal{F}\left( u_{\varepsilon }+\varphi _{\varepsilon },B_{s}\left(
x_{0}\right) \right)
\end{equation}%
Using (3.1)\ we get%
\begin{equation}
\begin{tabular}{l}
$\int\limits_{\Omega }\varepsilon \left\vert \nabla u_{\varepsilon
}^{1}\right\vert ^{p}\,dx+\int\limits_{\Omega }\varepsilon
\sum\limits_{\alpha =2}^{m}\left\vert \nabla u_{\varepsilon }^{\alpha
}\right\vert ^{p}\,dx+\int\limits_{\Omega }G\left( \left\vert u_{\varepsilon
}^{1}\right\vert ,...,\left\vert u_{\varepsilon }^{m}\right\vert ,\left\vert
\nabla u_{\varepsilon }^{1}\right\vert ,...,\left\vert \nabla u_{\varepsilon
}^{m}\right\vert \right) \,dx$ \\ 
$\leq \int\limits_{\Omega }\varepsilon \left\vert \nabla u_{\varepsilon
}^{1}+\nabla \varphi _{\varepsilon }^{1}\right\vert
^{p}\,dx+\int\limits_{\Omega }\varepsilon \sum\limits_{\alpha
=2}^{m}\left\vert \nabla u_{\varepsilon }^{\alpha }\right\vert ^{p}\,dx$ \\ 
$+\int\limits_{\Omega }G\left( \left\vert u_{\varepsilon }^{1}+\varphi
_{\varepsilon }^{1}\right\vert ,\left\vert u_{\varepsilon }^{2}\right\vert
,...,\left\vert u_{\varepsilon }^{m}\right\vert ,\left\vert \nabla
u_{\varepsilon }^{1}+\nabla \varphi _{\varepsilon }^{1}\right\vert
,\left\vert \nabla u_{\varepsilon }^{2}\right\vert ,...,\left\vert \nabla
u_{\varepsilon }^{m}\right\vert \right) \,dx$%
\end{tabular}%
\end{equation}%
and%
\begin{equation}
\begin{tabular}{l}
$\int\limits_{B_{s}\left( x_{0}\right) }\varepsilon \left\vert \nabla
u_{\varepsilon }^{1}\right\vert ^{p}\,dx+\int\limits_{B_{s}\left(
x_{0}\right) }G\left( \left\vert u_{\varepsilon }^{1}\right\vert
,...,\left\vert u_{\varepsilon }^{m}\right\vert ,\left\vert \nabla
u_{\varepsilon }^{1}\right\vert ,...,\left\vert \nabla u_{\varepsilon
}^{m}\right\vert \right) \,dx$ \\ 
$+\int\limits_{\Omega \backslash B_{s}\left( x_{0}\right) }G\left(
\left\vert u_{\varepsilon }^{1}\right\vert ,...,\left\vert u_{\varepsilon
}^{m}\right\vert ,\left\vert \nabla u_{\varepsilon }^{1}\right\vert
,...,\left\vert \nabla u_{\varepsilon }^{m}\right\vert \right) \,dx$ \\ 
$\leq \int\limits_{B_{s}\left( x_{0}\right) }\varepsilon \left\vert \nabla
u_{\varepsilon }^{1}+\nabla \varphi _{\varepsilon }^{1}\right\vert
^{p}\,dx+\int\limits_{B_{s}\left( x_{0}\right) }G\left( \left\vert
u_{\varepsilon }^{1}+\varphi _{\varepsilon }^{1}\right\vert ,\left\vert
u_{\varepsilon }^{2}\right\vert ,...,\left\vert u_{\varepsilon
}^{m}\right\vert ,\left\vert \nabla u_{\varepsilon }^{1}+\nabla \varphi
_{\varepsilon }^{1}\right\vert ,\left\vert \nabla u_{\varepsilon
}^{2}\right\vert ,...,\left\vert \nabla u_{\varepsilon }^{m}\right\vert
\right) \,dx$ \\ 
$+\int\limits_{\Omega \backslash B_{s}\left( x_{0}\right) }G\left(
\left\vert u_{\varepsilon }^{1}\right\vert ,...,\left\vert u_{\varepsilon
}^{m}\right\vert ,\left\vert \nabla u_{\varepsilon }^{1}\right\vert
,...,\left\vert \nabla u_{\varepsilon }^{m}\right\vert \right) \,dx$%
\end{tabular}%
\end{equation}%
Using (3.3) and ipotesi H.1 we have%
\begin{equation*}
\begin{tabular}{l}
$\int\limits_{A_{k,\varrho }^{1,\varepsilon }}\varepsilon \left\vert \nabla
u_{\varepsilon }^{1}\right\vert ^{p}\,dx+\int\limits_{A_{k,\varrho
}^{1,\varepsilon }}\sum\limits_{\alpha =1}^{m}\left\vert \nabla
u_{\varepsilon }^{\alpha }\right\vert ^{q}\,dx$ \\ 
$\leq \int\limits_{A_{k,\varrho }^{1,\varepsilon }}\varepsilon \left( 1-\eta
^{p}\right) \left\vert \nabla u_{\varepsilon }^{1}\right\vert
^{p}\,dx+\varepsilon p^{p}\int\limits_{A_{k,\varrho }^{1,\varepsilon
}\backslash A_{k,t}^{1,\varepsilon }}\eta ^{p-1}\left\vert \nabla \eta
\right\vert ^{p}\left( u_{\varepsilon }^{1}-k\right) ^{p}\,dx+$ \\ 
$\int\limits_{A_{k,\varrho }^{1,\varepsilon }}L\left( 1-\eta ^{p}\right)
^{q}\left\vert \nabla u_{\varepsilon }^{1}\right\vert ^{q}+p^{q}\eta
^{\left( p-1\right) q}\left\vert \nabla \eta \right\vert ^{q}\left(
u_{\varepsilon }^{1}-k\right) ^{q}+L\sum\limits_{\alpha =2}^{m}\left\vert
\nabla u_{\varepsilon }^{\alpha }\right\vert ^{q}\,dx+$ \\ 
$\int\limits_{A_{k,\varrho }^{1,\varepsilon }}L\left\vert \left( 1-\eta
^{p}\right) \left( u_{\varepsilon }^{1}-k\right) +k\right\vert
^{q}+L\sum\limits_{\alpha =2}^{m}\left\vert u_{\varepsilon }^{\alpha
}\right\vert ^{q}+L\,a\left( x\right) \,dx$%
\end{tabular}%
\end{equation*}%
where $A_{k,\varrho }^{1,\varepsilon }=\left\{ u_{\varepsilon
}^{1}>k\right\} \cap B_{\varrho }\left( x_{0}\right) $. Since $1-\eta ^{p}=0$
on $A_{k,\varrho }^{1,\varepsilon }$ and $0<\varepsilon \leq 1$, it follows
that%
\begin{equation}
\begin{tabular}{l}
$\ \int\limits_{A_{k,\varrho }^{1,\varepsilon }}\varepsilon \left\vert
\nabla u_{\varepsilon }^{1}\right\vert ^{p}\,+\left\vert \nabla
u_{\varepsilon }^{1}\right\vert ^{q}\,dx$ \\ 
$\leq L\int\limits_{A_{k,\varrho }^{1,\varepsilon }\backslash
A_{k,t}^{1.\varepsilon }}\varepsilon \left\vert \nabla u_{\varepsilon
}^{1}\right\vert ^{p}+\left\vert \nabla u_{\varepsilon }^{1}\right\vert
^{q}\,dx+p^{p}\int\limits_{A_{k,\varrho }^{1,\varepsilon }\backslash
A_{k,t}^{1,\varepsilon }}\eta ^{p-1}\left\vert \nabla \eta \right\vert
^{p}\left( u_{\varepsilon }^{1}-k\right) ^{p}\,dx+$ \\ 
$\int\limits_{A_{k,\varrho }^{1,\varepsilon }}L\left\vert \left( 1-\eta
^{p}\right) \left( u_{\varepsilon }^{1}-k\right) +k\right\vert
^{q}+L\sum\limits_{\alpha =1}^{m}\left\vert u_{\varepsilon }^{\alpha
}\right\vert ^{q}+L\,a\left( x\right) \,dx+$ \\ 
$\int\limits_{A_{k,\varrho }^{1,\varepsilon }}p^{q}\eta ^{\left( p-1\right)
q}\left\vert \nabla \eta \right\vert ^{q}\left( u_{\varepsilon
}^{1}-k\right) ^{q}+L\sum\limits_{\alpha =2}^{m}\left\vert \nabla
u_{\varepsilon }^{\alpha }\right\vert ^{q}\,dx$%
\end{tabular}
\label{5.1}
\end{equation}%
We observe that we can estimate the term 
\begin{equation*}
L\int\limits_{A_{k,\varrho }^{1,\varepsilon }}\sum\limits_{\alpha
=1}^{m}\left\vert u_{\varepsilon }^{\alpha }\right\vert ^{q}+\,a\left(
x\right) \,dx
\end{equation*}%
using the H\"{o}lder inequality \ref{4.2} obtaining%
\begin{equation}
\begin{tabular}{l}
$L\int\limits_{A_{k,\varrho }^{1,\varepsilon }}\sum\limits_{\alpha
=1}^{m}\left\vert u_{\varepsilon }^{\alpha }\right\vert ^{q}+\,a\left(
x\right) \,dx$ \\ 
$\leq mL\left[ \mathcal{L}^{n}\left( A_{k,\varrho }^{1,\varepsilon }\right) %
\right] ^{1-\frac{q}{p}}\left[ \int\limits_{A_{k,\varrho }^{1,\varepsilon
}}\left\vert u_{\varepsilon }\right\vert ^{p}\,dx\right] ^{\frac{q}{p}}+$ \\ 
$L\left[ \mathcal{L}^{n}\left( A_{k,\varrho }^{1,\varepsilon }\right) \right]
^{1-\frac{1}{\sigma }}\left[ \int\limits_{A_{k,\varrho }^{1,\varepsilon
}}a^{\sigma }\,dx\right] ^{\frac{1}{\sigma }}$%
\end{tabular}
\label{5.2}
\end{equation}%
Similarly we can also estimate the term%
\begin{equation*}
\int\limits_{A_{k,\varrho }^{1,\varepsilon }}L\sum\limits_{\alpha
=2}^{m}\left\vert \nabla u_{\varepsilon }^{\alpha }\right\vert ^{q}\,dx
\end{equation*}%
using the H\"{o}lder inequality \ref{4.2} obtaining%
\begin{equation}
\begin{tabular}{l}
$\int\limits_{A_{k,\varrho }^{1,\varepsilon }}L\sum\limits_{\alpha
=2}^{m}\left\vert \nabla u_{\varepsilon }^{\alpha }\right\vert ^{q}\,dx$ \\ 
$\leq mL\left[ \mathcal{L}^{n}\left( A_{k,\varrho }^{1,\varepsilon }\right) %
\right] ^{1-\frac{q}{p}}\left[ \int\limits_{A_{k,\varrho }^{1,\varepsilon
}}\left\vert \nabla u_{\varepsilon }\right\vert ^{p}\,dx\right] ^{\frac{q}{p}%
}$%
\end{tabular}
\label{5.3}
\end{equation}%
Using (\ref{5.1}), (\ref{5.2}) and (\ref{5.3})\ we get%
\begin{equation}
\begin{tabular}{l}
$\ \int\limits_{A_{k,\varrho }^{1,\varepsilon }}\varepsilon \left\vert
\nabla u_{\varepsilon }^{1}\right\vert ^{p}\,+\left\vert \nabla
u_{\varepsilon }^{1}\right\vert ^{q}\,dx$ \\ 
$\leq L\int\limits_{A_{k,\varrho }^{1,\varepsilon }\backslash
A_{k,t}^{1}}\varepsilon \left\vert \nabla u_{\varepsilon }^{1}\right\vert
^{p}+\left\vert \nabla u_{\varepsilon }^{1}\right\vert
^{q}\,dx+p^{p}\int\limits_{A_{k,\varrho }^{1,\varepsilon }\backslash
A_{k,t}^{1}}\eta ^{p-1}\left\vert \nabla \eta \right\vert ^{p}\left(
u_{\varepsilon }^{1}-k\right) ^{p}\,dx+$ \\ 
$\int\limits_{A_{k,\varrho }^{1,\varepsilon }}p^{q}\eta ^{\left( p-1\right)
q}\left\vert \nabla \eta \right\vert ^{q}\left( u_{\varepsilon
}^{1}-k\right) ^{q}\,dx+\int\limits_{A_{k,\varrho }^{1,\varepsilon
}}L\left\vert \left( 1-\eta ^{p}\right) \left( u_{\varepsilon }^{1}-k\right)
+k\right\vert ^{q}\,dx+$ \\ 
$2mL\left[ \mathcal{L}^{n}\left( A_{k,\varrho }^{1,\varepsilon }\right) %
\right] ^{1-\frac{q}{p}}\left\Vert u_{\varepsilon }\right\Vert
_{W^{1,p}\left( A_{k,\varrho }^{1,\varepsilon }\right) }^{q}+L\left[ 
\mathcal{L}^{n}\left( A_{k,\varrho }^{1,\varepsilon }\right) \right] ^{1-%
\frac{1}{\sigma }}\left\Vert a\right\Vert _{L^{\sigma }\left( A_{k,\varrho
}^{1,\varepsilon }\right) }$%
\end{tabular}
\label{5.4}
\end{equation}%
Using the Young Inequality \ref{4.1} we have%
\begin{equation}
\begin{tabular}{l}
$\ \int\limits_{A_{k,\varrho }^{1,\varepsilon }}p^{q}\eta ^{\left(
p-1\right) q}\left\vert \nabla \eta \right\vert ^{q}\left( u_{\varepsilon
}^{1}-k\right) ^{q}\,dx$ \\ 
$\leq p^{q}\int\limits_{A_{k,\varrho }^{1,\varepsilon }}1+\left\vert \nabla
\eta \right\vert ^{p}\left( u_{\varepsilon }^{1}-k\right) ^{p}\,dx$ \\ 
$\leq p^{q}\mathcal{L}^{n}\left( A_{k,\varrho }^{1,\varepsilon }\right)
+p^{q}\int\limits_{A_{k,\varrho }^{1,\varepsilon }}\left\vert \nabla \eta
\right\vert ^{p}\left( u_{\varepsilon }^{1}-k\right) ^{p}\,dx$%
\end{tabular}
\label{5.5}
\end{equation}%
then using (\ref{5.4}) and (\ref{5.5}) it follows%
\begin{equation*}
\begin{tabular}{l}
$\ \int\limits_{A_{k,\varrho }^{1,\varepsilon }}\varepsilon \left\vert
\nabla u_{\varepsilon }^{1}\right\vert ^{p}\,+\left\vert \nabla
u_{\varepsilon }^{1}\right\vert ^{q}\,dx$ \\ 
$\leq L\int\limits_{A_{k,\varrho }^{1,\varepsilon }\backslash
A_{k,t}^{1}}\varepsilon \left\vert \nabla u_{\varepsilon }^{1}\right\vert
^{p}+\left\vert \nabla u_{\varepsilon }^{1}\right\vert ^{q}\,dx+\left(
p^{p}+p^{q}\right) \int\limits_{A_{k,\varrho }^{1,\varepsilon }\backslash
A_{k,t}^{1}}\eta ^{p-1}\left\vert \nabla \eta \right\vert ^{p}\left(
u_{\varepsilon }^{1}-k\right) ^{p}\,dx+$ \\ 
$p^{q}\mathcal{L}^{n}\left( A_{k,\varrho }^{1,\varepsilon }\right)
+\int\limits_{A_{k,\varrho }^{1,\varepsilon }}L\left\vert \left( 1-\eta
^{p}\right) \left( u_{\varepsilon }^{1}-k\right) +k\right\vert ^{q}\,dx+$ \\ 
$2mL\left[ \mathcal{L}^{n}\left( A_{k,\varrho }^{1,\varepsilon }\right) %
\right] ^{1-\frac{q}{p}}\left\Vert u_{\varepsilon }\right\Vert
_{W^{1,p}\left( A_{k,\varrho }^{1,\varepsilon }\right) }^{q}+L\left[ 
\mathcal{L}^{n}\left( A_{k,\varrho }^{1,\varepsilon }\right) \right] ^{1-%
\frac{1}{\sigma }}\left\Vert a\right\Vert _{L^{\sigma }\left( A_{k,\varrho
}^{1,\varepsilon }\right) }$%
\end{tabular}%
\end{equation*}%
and%
\begin{equation*}
\begin{tabular}{l}
$\ \int\limits_{A_{k,\varrho }^{1,\varepsilon }}\varepsilon \left\vert
\nabla u_{\varepsilon }^{1}\right\vert ^{p}\,+\left\vert \nabla
u_{\varepsilon }^{1}\right\vert ^{q}\,dx$ \\ 
$\leq L\int\limits_{A_{k,\varrho }^{1,\varepsilon }\backslash
A_{k,t}^{1}}\varepsilon \left\vert \nabla u_{\varepsilon }^{1}\right\vert
^{p}+\left\vert \nabla u_{\varepsilon }^{1}\right\vert ^{q}\,dx+\left(
p^{p}+p^{q}\right) \int\limits_{A_{k,\varrho }^{1,\varepsilon }\backslash
A_{k,t}^{1}}\eta ^{p-1}\left\vert \nabla \eta \right\vert ^{p}\left(
u_{\varepsilon }^{1}-k\right) ^{p}\,dx+$ \\ 
$\left( p^{q}+2^{p-1}L\,k^{p}\right) \mathcal{L}^{n}\left( A_{k,\varrho
}^{1,\varepsilon }\right) +2^{p-1}L\int\limits_{A_{k,s}^{1,\varepsilon
}\backslash A_{k,t}^{1,\varepsilon }}\left( u_{\varepsilon }^{1}-k\right)
^{p}\,dx+$ \\ 
$2mL\left[ \mathcal{L}^{n}\left( A_{k,\varrho }^{1,\varepsilon }\right) %
\right] ^{1-\frac{q}{p}}\left\Vert u_{\varepsilon }\right\Vert
_{W^{1,p}\left( A_{k,\varrho }^{1,\varepsilon }\right) }^{q}+L\left[ 
\mathcal{L}^{n}\left( A_{k,\varrho }^{1,\varepsilon }\right) \right] ^{1-%
\frac{1}{\sigma }}\left\Vert a\right\Vert _{L^{\sigma }\left( A_{k,\varrho
}^{1,\varepsilon }\right) }$%
\end{tabular}%
\end{equation*}%
Moreover, since $1\leq q<\frac{p^{2}}{n}$ and $\sigma >\frac{n}{p}$, we get%
\begin{equation*}
\begin{tabular}{l}
$\int\limits_{A_{k,\varrho }^{1,\varepsilon }}\varepsilon \left\vert \nabla
u_{\varepsilon }^{1}\right\vert ^{p}\,+\left\vert \nabla u_{\varepsilon
}^{1}\right\vert ^{q}\,dx$ \\ 
$\leq L\int\limits_{A_{k,\varrho }^{1,\varepsilon }\backslash
A_{k,t}^{1}}\varepsilon \left\vert \nabla u_{\varepsilon }^{1}\right\vert
^{p}+\left\vert \nabla u_{\varepsilon }^{1}\right\vert ^{q}\,dx+\frac{%
2^{p}\left( p^{p}+p^{q}\right) }{\left( s-t\right) ^{p}}\int%
\limits_{A_{k,s}^{1,\varepsilon }\backslash A_{k,t}^{1,\varepsilon }}\left(
u_{\varepsilon }^{1}-k\right) ^{p}\,dx+$ \\ 
$D_{\Sigma ,\varepsilon }\left( 1+k^{p}\right) \left[ \mathcal{L}^{n}\left(
A_{k,s}^{1,\varepsilon }\right) \right] ^{1-\frac{p}{n}+\epsilon }$%
\end{tabular}%
\end{equation*}%
where%
\begin{equation*}
D_{\Sigma ,\varepsilon }=2^{p-1}L\left( p^{q}+L+2m\left( L+1\right)
\left\Vert u_{\varepsilon }\right\Vert _{W^{1,p}\left( \Sigma \right)
}^{q}+\left( L+1\right) \left\Vert a\right\Vert _{L^{\sigma }\left( \Sigma
\right) }\right)
\end{equation*}%
then it follows%
\begin{equation*}
\begin{tabular}{l}
$\int\limits_{A_{k,t}^{1,\varepsilon }}\varepsilon \left\vert \nabla
u_{\varepsilon }^{1}\right\vert ^{p}\,+\left\vert \nabla u_{\varepsilon
}^{1}\right\vert ^{q}\,dx$ \\ 
$\leq \frac{L}{L+1}\int\limits_{A_{k,s}^{1,\varepsilon }}\varepsilon
\left\vert \nabla u_{\varepsilon }^{1}\right\vert ^{p}+\left\vert \nabla
u_{\varepsilon }^{1}\right\vert ^{q}\,dx+\frac{2^{p}\left(
p^{p}+p^{q}\right) }{\left( L+1\right) \left( s-t\right) ^{p}}%
\int\limits_{A_{k,s}^{1,\varepsilon }}\left( u_{\varepsilon }^{1}-k\right)
^{p}\,dx+$ \\ 
$\frac{D_{\Sigma ,\varepsilon }}{L+1}\left( 1+k^{p}\right) \left[ \mathcal{L}%
^{n}\left( A_{k,s}^{1,\varepsilon }\right) \right] ^{1-\frac{p}{n}%
+\varepsilon }$%
\end{tabular}%
\end{equation*}%
Using Lemma \ref{lemma3} we have the Caccioppoli Inequaity%
\begin{equation*}
\begin{tabular}{l}
$\int\limits_{A_{k,\varrho }^{1,\varepsilon }}\varepsilon \left\vert \nabla
u_{\varepsilon }^{1}\right\vert ^{p}\,+\left\vert \nabla u_{\varepsilon
}^{1}\right\vert ^{q}\,dx$ \\ 
$\leq \frac{C_{1,\Sigma }}{\left( R-\varrho \right) ^{p}}\int%
\limits_{A_{k,R}^{1,\varepsilon }}\left( u_{\varepsilon }^{1}-k\right)
^{p}\,dx+C_{2,\Sigma ,\varepsilon }\left( 1+R^{-\epsilon n}k^{p}\right) %
\left[ \mathcal{L}^{n}\left( A_{k,R}^{1,\varepsilon }\right) \right] ^{1-%
\frac{p}{n}+\epsilon }$%
\end{tabular}%
\end{equation*}%
Similarly you can proceed for $\alpha =2,...,m$ and we get $u_{\varepsilon
}^{\alpha }\in DG^{+}\left( \Omega ,p,\lambda ,\lambda _{\ast },\chi
,\varepsilon ,R_{0},k_{0}\right) $ for every $\alpha =1,...,m$, with $%
\lambda =C_{1,\Sigma }$, $\lambda _{\ast }=C_{2,\Sigma ,\varepsilon }$ and $%
\chi =1$. Since $-u$ is a minimizer of the integral functional 
\begin{equation*}
\mathcal{\tilde{F}}_{\varepsilon }\left( v,\Omega \right)
=\int\limits_{\Omega }\varepsilon \sum\limits_{\alpha =1}^{m}\left\vert
\nabla v^{\alpha }\right\vert ^{p}+G\left( \sqrt{\left\vert v^{1}\right\vert
^{2}+\varepsilon },...,\sqrt{\left\vert v^{m}\right\vert ^{2}+\varepsilon },%
\sqrt{\left\vert \nabla v^{1}\right\vert ^{2}+\varepsilon },...,\sqrt{%
\left\vert \nabla v^{m}\right\vert ^{2}+\varepsilon }\right) \,dx
\end{equation*}%
then $u_{\varepsilon }^{\alpha }\in DG^{-}\left( \Omega ,p,\lambda ,\lambda
_{\ast },\chi ,\varepsilon ,R_{0},k_{0}\right) $ for every $\alpha =1,...,m$%
, with $\lambda =C_{1,\Sigma }$, $\lambda _{\ast }=C_{2,\Sigma ,\varepsilon
} $ and $\chi =1$. It follows that $u_{\varepsilon }^{\alpha }\in DG\left(
\Omega ,p,\lambda ,\lambda _{\ast },\chi ,\varepsilon ,R_{0},k_{0}\right) $
for every $\alpha =1,...,m$, with $\lambda =C_{1,\Sigma }$, $\lambda _{\ast
}=C_{2,\Sigma }$ and $\chi =1$.

\section{Proof of Theorem 4}

Let us consider $u_{\varepsilon }\in W^{1,p}\left( \Omega ,%
\mathbb{R}
^{m}\right) $, a local minimizer of the functional 1.5, then%
\begin{equation*}
\begin{tabular}{l}
$0$ \\ 
$=\int\limits_{\Omega }\varepsilon \left\vert \nabla u_{\varepsilon
}^{\alpha }\right\vert ^{p-2}\nabla u_{\varepsilon }^{\alpha }\nabla \varphi
^{\alpha }dx$ \\ 
$+\int\limits_{\Omega }\partial _{s^{\alpha }}G\left( \sqrt{\left\vert
u_{\varepsilon }^{1}\right\vert ^{2}+\varepsilon },...,\sqrt{\left\vert
u_{\varepsilon }^{m}\right\vert ^{2}+\varepsilon },\sqrt{\left\vert \nabla
u_{\varepsilon }^{1}\right\vert ^{2}+\varepsilon },...,\sqrt{\left\vert
\nabla u_{\varepsilon }^{m}\right\vert ^{2}+\varepsilon }\right) \frac{%
u_{\varepsilon }^{\alpha }\varphi ^{\alpha }}{\sqrt{\left\vert
u_{\varepsilon }^{\alpha }\right\vert ^{2}+\varepsilon }}\,dx$ \\ 
$+\int\limits_{\Omega }\partial _{\xi ^{\alpha }}G\left( \sqrt{\left\vert
u_{\varepsilon }^{1}\right\vert ^{2}+\varepsilon },...,\sqrt{\left\vert
u_{\varepsilon }^{m}\right\vert ^{2}+\varepsilon },\sqrt{\left\vert \nabla
u_{\varepsilon }^{1}\right\vert ^{2}+\varepsilon },...,\sqrt{\left\vert
\nabla u_{\varepsilon }^{m}\right\vert ^{2}+\varepsilon }\right) \frac{%
\nabla u_{\varepsilon }^{\alpha }\nabla \varphi ^{\alpha }}{\sqrt{\left\vert
\nabla u_{\varepsilon }^{\alpha }\right\vert ^{2}+\varepsilon }}\,dx$%
\end{tabular}%
\end{equation*}%
for every $\alpha =1,...,m$. Fix $R_{0}>0$, let us take $\varphi ^{\alpha
}=\triangle _{-h,s}\left( \eta ^{2}\triangle _{h,s}u_{\varepsilon }^{\alpha
}\right) $ where $\eta \in C_{0}^{\infty }\left( B_{\varrho }\left(
x_{0}\right) \right) $ is a be a cut-of function, where $0<\tau <\varrho <$ $%
\frac{1}{2}\min \left\{ R_{0},dist\left( x_{0},\partial \Omega \right)
,1\right\} $, $0\leq \eta \leq 1$ on $B_{\varrho }\left( x_{0}\right) $, $%
\eta =1$ on $B_{\tau }\left( x_{0}\right) $, $\left\vert \nabla \eta
\right\vert \leq \frac{c}{\varrho -\tau }$ on $B_{\varrho }\left(
x_{0}\right) $\ and 
\begin{equation}
\triangle _{h,s}f\left( x\right) =\frac{f\left( x+h\,e_{s}\right) -f\left(
x\right) }{h}
\end{equation}%
for every $h\in \left( -\frac{\varrho }{4},\frac{\varrho }{4}\right)
\backslash \left\{ 0\right\} $ and for every $e_{s}\in 
\mathbb{R}
^{n}$ with $\left\vert e_{s}\right\vert =1$,then we get%
\begin{equation*}
\begin{tabular}{l}
$0$ \\ 
$=\int\limits_{\Omega }\varepsilon \left\vert \nabla u_{\varepsilon
}^{\alpha }\right\vert ^{p-2}\nabla u_{\varepsilon }^{\alpha }\nabla \left(
\triangle _{-h,s}\left( \eta ^{2}\triangle _{h,s}u_{\varepsilon }^{\alpha
}\right) \right) dx$ \\ 
$+\int\limits_{\Omega }\partial _{s^{\alpha }}G\left( \sqrt{\left\vert
u_{\varepsilon }^{1}\right\vert ^{2}+\varepsilon },...,\sqrt{\left\vert
u_{\varepsilon }^{m}\right\vert ^{2}+\varepsilon },\sqrt{\left\vert \nabla
u_{\varepsilon }^{1}\right\vert ^{2}+\varepsilon },...,\sqrt{\left\vert
\nabla u_{\varepsilon }^{m}\right\vert ^{2}+\varepsilon }\right) $ \\ 
$\frac{u_{\varepsilon }^{\alpha }\left( \triangle _{-h,s}\left( \eta
^{2}\triangle _{h,s}u_{\varepsilon }^{\alpha }\right) \right) }{\sqrt{%
\left\vert u_{\varepsilon }^{\alpha }\right\vert ^{2}+\varepsilon }}\,dx$ \\ 
$+\int\limits_{\Omega }\partial _{\xi ^{\alpha }}G\left( \sqrt{\left\vert
u_{\varepsilon }^{1}\right\vert ^{2}+\varepsilon },...,\sqrt{\left\vert
u_{\varepsilon }^{m}\right\vert ^{2}+\varepsilon },\sqrt{\left\vert \nabla
u_{\varepsilon }^{1}\right\vert ^{2}+\varepsilon },...,\sqrt{\left\vert
\nabla u_{\varepsilon }^{m}\right\vert ^{2}+\varepsilon }\right) $ \\ 
$\frac{\nabla u_{\varepsilon }^{\alpha }\nabla \left( \triangle
_{-h,s}\left( \eta ^{2}\triangle _{h,s}u_{\varepsilon }^{\alpha }\right)
\right) }{\sqrt{\left\vert \nabla u_{\varepsilon }^{\alpha }\right\vert
^{2}+\varepsilon }}\,dx$%
\end{tabular}%
\end{equation*}%
from Lemma 5 follows%
\begin{equation*}
\begin{tabular}{l}
$0$ \\ 
$=\int\limits_{\Omega }\varepsilon \left\vert \nabla u_{\varepsilon
}^{\alpha }\right\vert ^{p-2}\nabla u_{\varepsilon }^{\alpha }\left(
\triangle _{-h,s}\nabla \left( \eta ^{2}\triangle _{h,s}u_{\varepsilon
}^{\alpha }\right) \right) dx$ \\ 
$+\int\limits_{\Omega }\partial _{s^{\alpha }}G\left( \sqrt{\left\vert
u_{\varepsilon }^{1}\right\vert ^{2}+\varepsilon },...,\sqrt{\left\vert
u_{\varepsilon }^{m}\right\vert ^{2}+\varepsilon },\sqrt{\left\vert \nabla
u_{\varepsilon }^{1}\right\vert ^{2}+\varepsilon },...,\sqrt{\left\vert
\nabla u_{\varepsilon }^{m}\right\vert ^{2}+\varepsilon }\right) $ \\ 
$\frac{u_{\varepsilon }^{\alpha }\left( \triangle _{-h,s}\left( \eta
^{2}\triangle _{h,s}u_{\varepsilon }^{\alpha }\right) \right) }{\sqrt{%
\left\vert u_{\varepsilon }^{\alpha }\right\vert ^{2}+\varepsilon }}\,dx$ \\ 
$+\int\limits_{\Omega }\partial _{\xi ^{\alpha }}G\left( \sqrt{\left\vert
u_{\varepsilon }^{1}\right\vert ^{2}+\varepsilon },...,\sqrt{\left\vert
u_{\varepsilon }^{m}\right\vert ^{2}+\varepsilon },\sqrt{\left\vert \nabla
u_{\varepsilon }^{1}\right\vert ^{2}+\varepsilon },...,\sqrt{\left\vert
\nabla u_{\varepsilon }^{m}\right\vert ^{2}+\varepsilon }\right) $ \\ 
$\frac{\nabla u_{\varepsilon }^{\alpha }\left( \triangle _{-h,s}\nabla
\left( \eta ^{2}\triangle _{h,s}u_{\varepsilon }^{\alpha }\right) \right) }{%
\sqrt{\left\vert \nabla u_{\varepsilon }^{\alpha }\right\vert
^{2}+\varepsilon }}\,dx$%
\end{tabular}%
\end{equation*}%
from which from Lemma 6 we have%
\begin{equation}
\begin{tabular}{l}
$0$ \\ 
$=\int\limits_{\Omega }\varepsilon \,\triangle _{h,s}\left( \left\vert
\nabla u_{\varepsilon }^{\alpha }\right\vert ^{p-2}\nabla u_{\varepsilon
}^{\alpha }\right) \nabla \left( \eta ^{2}\triangle _{h,s}u_{\varepsilon
}^{\alpha }\right) dx$ \\ 
$+\int\limits_{\Omega }\triangle _{h,s}\left[ \partial _{s^{\alpha }}G\left( 
\sqrt{\left\vert u_{\varepsilon }^{1}\right\vert ^{2}+\varepsilon },...,%
\sqrt{\left\vert u_{\varepsilon }^{m}\right\vert ^{2}+\varepsilon },\sqrt{%
\left\vert \nabla u_{\varepsilon }^{1}\right\vert ^{2}+\varepsilon },...,%
\sqrt{\left\vert \nabla u_{\varepsilon }^{m}\right\vert ^{2}+\varepsilon }%
\right) \frac{u_{\varepsilon }^{\alpha }}{\sqrt{\left\vert u_{\varepsilon
}^{\alpha }\right\vert ^{2}+\varepsilon }}\right] $ \\ 
$\left( \eta ^{2}\triangle _{h,s}u_{\varepsilon }^{\alpha }\right) \,dx$ \\ 
$+\int\limits_{\Omega }\triangle _{h,s}\left[ \partial _{\xi ^{\alpha
}}G\left( \sqrt{\left\vert u_{\varepsilon }^{1}\right\vert ^{2}+\varepsilon }%
,...,\sqrt{\left\vert u_{\varepsilon }^{m}\right\vert ^{2}+\varepsilon },%
\sqrt{\left\vert \nabla u_{\varepsilon }^{1}\right\vert ^{2}+\varepsilon }%
,...,\sqrt{\left\vert \nabla u_{\varepsilon }^{m}\right\vert
^{2}+\varepsilon }\right) \frac{\nabla u_{\varepsilon }^{\alpha }}{\sqrt{%
\left\vert \nabla u_{\varepsilon }^{\alpha }\right\vert ^{2}+\varepsilon }}%
\right] $ \\ 
$\nabla \left( \eta ^{2}\triangle _{h,s}u_{\varepsilon }^{\alpha }\right)
\,dx$%
\end{tabular}%
\end{equation}%
Now let's analyze the following three elements%
\begin{equation}
\triangle _{h,s}\left( \left\vert \nabla u_{\varepsilon }^{\alpha
}\right\vert ^{p-2}\nabla u_{\varepsilon }^{\alpha }\right) 
\end{equation}%
\begin{equation}
\triangle _{h,s}\left[ \partial _{s^{\alpha }}G\left( \sqrt{\left\vert
u_{\varepsilon }^{1}\right\vert ^{2}+\varepsilon },...,\sqrt{\left\vert
u_{\varepsilon }^{m}\right\vert ^{2}+\varepsilon },\sqrt{\left\vert \nabla
u_{\varepsilon }^{1}\right\vert ^{2}+\varepsilon },...,\sqrt{\left\vert
\nabla u_{\varepsilon }^{m}\right\vert ^{2}+\varepsilon }\right) \frac{%
u_{\varepsilon }^{\alpha }}{\sqrt{\left\vert u_{\varepsilon }^{\alpha
}\right\vert ^{2}+\varepsilon }}\right] 
\end{equation}%
\begin{equation}
\triangle _{h,s}\left[ \partial _{\xi ^{\alpha }}G\left( \sqrt{\left\vert
u_{\varepsilon }^{1}\right\vert ^{2}+\varepsilon },...,\sqrt{\left\vert
u_{\varepsilon }^{m}\right\vert ^{2}+\varepsilon },\sqrt{\left\vert \nabla
u_{\varepsilon }^{1}\right\vert ^{2}+\varepsilon },...,\sqrt{\left\vert
\nabla u_{\varepsilon }^{m}\right\vert ^{2}+\varepsilon }\right) \frac{%
\nabla u_{\varepsilon }^{\alpha }}{\sqrt{\left\vert \nabla u_{\varepsilon
}^{\alpha }\right\vert ^{2}+\varepsilon }}\right] 
\end{equation}%
We begin by estimating the relation defined in (4.3)%
\begin{equation*}
\begin{tabular}{l}
$\triangle _{h,s}\left( \left\vert \nabla u_{\varepsilon }^{\alpha
}\right\vert ^{p-2}\nabla u_{\varepsilon }^{\alpha }\right) $ \\ 
$=\frac{\left\vert \nabla u_{\varepsilon }^{\alpha }\left( x+he_{s}\right)
\right\vert ^{p-2}\nabla u_{\varepsilon }^{\alpha }\left( x+he_{s}\right)
-\left\vert \nabla u_{\varepsilon }^{\alpha }\left( x\right) \right\vert
^{p-2}\nabla u_{\varepsilon }^{\alpha }\left( x\right) }{h}$ \\ 
$=\frac{1}{h}\left[ \left\vert \nabla u_{\varepsilon }^{\alpha }\left(
x\right) +th\triangle _{h,s}\nabla u_{\varepsilon }^{\alpha }\left( x\right)
\right\vert ^{p-2}\left( \nabla u_{\varepsilon }^{\alpha }\left( x\right)
+th\triangle _{h,s}\nabla u_{\varepsilon }^{\alpha }\left( x\right) \right) %
\right] _{0}^{1}$ \\ 
$=\frac{1}{h}\int\limits_{0}^{1}\frac{d}{dt}\left[ \left\vert \nabla
u_{\varepsilon }^{\alpha }\left( x\right) +th\triangle _{h,s}\nabla
u_{\varepsilon }^{\alpha }\left( x\right) \right\vert ^{p-2}\left( \nabla
u_{\varepsilon }^{\alpha }\left( x\right) +th\triangle _{h,s}\nabla
u_{\varepsilon }^{\alpha }\left( x\right) \right) \right] \,dt$ \\ 
$=\frac{1}{h}\int\limits_{0}^{1}\left( p-2\right) \left\vert \nabla
u_{\varepsilon }^{\alpha }\left( x\right) +th\triangle _{h,s}\nabla
u_{\varepsilon }^{\alpha }\left( x\right) \right\vert ^{p-3}\frac{\left(
\nabla u_{\varepsilon }^{\alpha }\left( x\right) +th\triangle _{h,s}\nabla
u_{\varepsilon }^{\alpha }\left( x\right) \right) ^{2}}{\left\vert \nabla
u_{\varepsilon }^{\alpha }\left( x\right) +th\triangle _{h,s}\nabla
u_{\varepsilon }^{\alpha }\left( x\right) \right\vert }h\triangle
_{h,s}\nabla u_{\varepsilon }^{\alpha }\left( x\right) $ \\ 
$+\left\vert \nabla u^{\alpha }\left( x\right) +th\triangle _{h,s}\nabla
u^{\alpha }\left( x\right) \right\vert ^{p-2}h\triangle _{h,s}\nabla
u^{\alpha }\left( x\right) \,dt$ \\
\end{tabular}%
\end{equation*}%
from which we get%
\begin{equation}
\begin{tabular}{l}
$\int\limits_{\Omega }\varepsilon \,\triangle _{h,s}\left( \left\vert \nabla
u_{\varepsilon }^{\alpha }\right\vert ^{p-2}\nabla u_{\varepsilon }^{\alpha
}\right) \nabla \left( \eta ^{2}\triangle _{h,s}u_{\varepsilon }^{\alpha
}\right) dx$ \\ 
$=\int\limits_{\Omega }\varepsilon \,\left( p-1\right) \,\triangle
_{h,s}\nabla u_{\varepsilon }^{\alpha }\left( x\right)
\int\limits_{0}^{1}\left\vert \nabla u_{\varepsilon }^{\alpha }\left(
x\right) +th\triangle _{h,s}\nabla u_{\varepsilon }^{\alpha }\left( x\right)
\right\vert ^{p-2}\,\,dt\,\nabla \left( \eta ^{2}\triangle
_{h,s}u_{\varepsilon }^{\alpha }\right) dx$ \\ 
$=\left( p-1\right) \int\limits_{\Omega }\varepsilon \,I_{1,\alpha }\,\left(
\triangle _{h,s}\left( \nabla u_{\varepsilon }^{\alpha }\right) \right)
^{2}\eta ^{2}dx$ \\ 
$+2\left( p-1\right) \int\limits_{\Omega }\varepsilon \,\triangle
_{h,s}\left( \nabla u_{\varepsilon }^{\alpha }\right) \left( \triangle
_{h,s}u_{\varepsilon }^{\alpha }\right) I_{1,\alpha }\,\eta \,\nabla \eta
\,dx$%
\end{tabular}%
\end{equation}

where%
\begin{equation}
I_{1,\alpha }=\int\limits_{0}^{1}\left\vert \nabla u_{\varepsilon }^{\alpha
}\left( x\right) +th\triangle _{h,s}\nabla u_{\varepsilon }^{\alpha }\left(
x\right) \right\vert ^{p-2}\,\,dt
\end{equation}%
Now we have to estimate the relation defined in (4.4) 
\begin{equation*}
\begin{tabular}{l}
$\triangle _{h,s}\left[ \partial _{s^{\alpha }}G\left( \sqrt{\left\vert
u_{\varepsilon }^{1}\right\vert ^{2}+\varepsilon },...,\sqrt{\left\vert
u_{\varepsilon }^{m}\right\vert ^{2}+\varepsilon },\sqrt{\left\vert \nabla
u_{\varepsilon }^{1}\right\vert ^{2}+\varepsilon },...,\sqrt{\left\vert
\nabla u_{\varepsilon }^{m}\right\vert ^{2}+\varepsilon }\right) \frac{%
u_{\varepsilon }^{\alpha }}{\sqrt{\left\vert u_{\varepsilon }^{\alpha
}\right\vert ^{2}+\varepsilon }}\right] $ \\ 
$=\frac{1}{h}\left[ \partial _{s^{\alpha }}G\left( \sqrt{\left\vert
u_{\varepsilon }^{1}+th\triangle _{h,s}u_{\varepsilon }^{1}\right\vert
^{2}+\varepsilon },...,\sqrt{\left\vert u^{m}+th\triangle
_{h,s}u_{\varepsilon }^{m}\right\vert ^{2}+\varepsilon },\right. \right. $
\\ 
$\left. \left. \sqrt{\left\vert \nabla u_{\varepsilon }^{1}+th\triangle
_{h,s}\nabla u_{\varepsilon }^{1}\right\vert ^{2}+\varepsilon },...,\sqrt{%
\left\vert \nabla u_{\varepsilon }^{m}+th\triangle _{h,s}\nabla
u_{\varepsilon }^{m}\right\vert ^{2}+\varepsilon }\right) \frac{%
u_{\varepsilon }^{\alpha }+th\triangle _{h,s}u_{\varepsilon }^{\alpha }}{%
\sqrt{\left\vert u^{\alpha }+th\triangle _{h,s}u^{\alpha }\right\vert
^{2}+\varepsilon }}\right] _{0}^{1}$ \\ 
$=\frac{1}{h}\int\limits_{0}^{1}\frac{d}{dt}\left( \partial _{s^{\alpha
}}G\left( \sqrt{\left\vert u_{\varepsilon }^{1}+th\triangle
_{h,s}u_{\varepsilon }^{1}\right\vert ^{2}+\varepsilon },...,\sqrt{%
\left\vert u_{\varepsilon }^{m}+th\triangle _{h,s}u_{\varepsilon
}^{m}\right\vert ^{2}+\varepsilon },\right. \right. $ \\
\end{tabular}%
\end{equation*}%
so long as $G\in C^{2}\left( 
\mathbb{R}
_{+,0}^{m}\times 
\mathbb{R}
_{+,0}^{m}\right) $ we obtain%
\begin{equation*}
\begin{tabular}{l}
$\triangle _{h,s}\left[ \partial _{s^{\alpha }}G\left( \sqrt{\left\vert
u_{\varepsilon }^{1}\right\vert ^{2}+\varepsilon },...,\sqrt{\left\vert
u_{\varepsilon }^{m}\right\vert ^{2}+\varepsilon },\sqrt{\left\vert \nabla
u_{\varepsilon }^{1}\right\vert ^{2}+\varepsilon },...,\sqrt{\left\vert
\nabla u_{\varepsilon }^{m}\right\vert ^{2}+\varepsilon }\right) \frac{%
u_{\varepsilon }^{\alpha }}{\sqrt{\left\vert u_{\varepsilon }^{\alpha
}\right\vert ^{2}+\varepsilon }}\right] $ \\ 
$=\frac{1}{h}\int\limits_{0}^{1}\sum\limits_{\beta =1}^{m}\partial
_{s^{\beta }s^{\alpha }}G\left( \sqrt{\left\vert u_{\varepsilon
}^{1}+th\triangle _{h,s}u_{\varepsilon }^{1}\right\vert ^{2}+\varepsilon }%
,...,\sqrt{\left\vert u_{\varepsilon }^{m}+th\triangle _{h,s}u_{\varepsilon
}^{m}\right\vert ^{2}+\varepsilon },\right. $ \\ 
$\left. \sqrt{\left\vert \nabla u_{\varepsilon }^{1}+th\triangle
_{h,s}\nabla u_{\varepsilon }^{1}\right\vert ^{2}+\varepsilon },...,\sqrt{%
\left\vert \nabla u_{\varepsilon }^{m}+th\triangle _{h,s}\nabla
u_{\varepsilon }^{m}\right\vert ^{2}+\varepsilon }\right) \frac{%
u_{\varepsilon }^{\beta }+th\triangle _{h,s}u_{\varepsilon }^{\beta }}{\sqrt{%
\left\vert u_{\varepsilon }^{\beta }+th\triangle _{h,s}u_{\varepsilon
}^{\beta }\right\vert ^{2}+\varepsilon }}$ \\ 
$\cdot \frac{u_{\varepsilon }^{\alpha }+th\triangle _{h,s}u_{\varepsilon
}^{\alpha }}{\sqrt{\left\vert u_{\varepsilon }^{\alpha }+th\triangle
_{h,s}u_{\varepsilon }^{\alpha }\right\vert ^{2}+\varepsilon }}h\triangle
_{h,s}u_{\varepsilon }^{\beta }$ \\ 
$+\sum\limits_{\beta =1}^{m}\partial _{\xi ^{\beta }s^{\alpha }}G\left( 
\sqrt{\left\vert u_{\varepsilon }^{1}+th\triangle _{h,s}u_{\varepsilon
}^{1}\right\vert ^{2}+\varepsilon },...,\sqrt{\left\vert u_{\varepsilon
}^{m}+th\triangle _{h,s}u_{\varepsilon }^{m}\right\vert ^{2}+\varepsilon }%
,\right. $ \\
\end{tabular}%
\end{equation*}%
which we can rewrite as follows%
\begin{equation}
\begin{tabular}{l}
$\triangle _{h,s}\left[ \partial _{s^{\alpha }}G\left( \sqrt{\left\vert
u_{\varepsilon }^{1}\right\vert ^{2}+\varepsilon },...,\sqrt{\left\vert
u_{\varepsilon }^{m}\right\vert ^{2}+\varepsilon },\sqrt{\left\vert \nabla
u_{\varepsilon }^{1}\right\vert ^{2}+\varepsilon },...,\sqrt{\left\vert
\nabla u_{\varepsilon }^{m}\right\vert ^{2}+\varepsilon }\right) \frac{%
u_{\varepsilon }^{\alpha }}{\sqrt{\left\vert u_{\varepsilon }^{\alpha
}\right\vert ^{2}+\varepsilon }}\right] $ \\ 
$=\sum\limits_{\beta =1}^{m}\int\limits_{0}^{1}\partial _{s^{\beta
}s^{\alpha }}G\left( x,t\right) \frac{u_{\varepsilon }^{\beta }+th\triangle
_{h,s}u_{\varepsilon }^{\beta }}{\sqrt{\left\vert u_{\varepsilon }^{\beta
}+th\triangle _{h,s}u_{\varepsilon }^{\beta }\right\vert ^{2}+\varepsilon }}%
\frac{u_{\varepsilon }^{\alpha }+th\triangle _{h,s}u_{\varepsilon }^{\alpha }%
}{\sqrt{\left\vert u_{\varepsilon }^{\alpha }+th\triangle
_{h,s}u_{\varepsilon }^{\alpha }\right\vert ^{2}+\varepsilon }}%
\,dt\,\triangle _{h,s}u_{\varepsilon }^{\beta }$ \\ 
$+\sum\limits_{\beta =1}^{m}\int\limits_{0}^{1}\partial _{\xi ^{\beta
}s^{\alpha }}G\left( x,t\right) \frac{\nabla u_{\varepsilon }^{\beta
}+th\triangle _{h,s}\nabla u_{\varepsilon }^{\beta }}{\sqrt{\left\vert
\nabla u_{\varepsilon }^{\beta }+th\triangle _{h,s}\nabla u_{\varepsilon
}^{\beta }\right\vert ^{2}+\varepsilon }}\frac{u_{\varepsilon }^{\alpha
}+th\triangle _{h,s}u_{\varepsilon }^{\alpha }}{\sqrt{\left\vert
u_{\varepsilon }^{\alpha }+th\triangle _{h,s}u_{\varepsilon }^{\alpha
}\right\vert ^{2}+\varepsilon }}\,dt\,\triangle _{h,s}\nabla u_{\varepsilon
}^{\beta }$ \\ 
$+\int\limits_{0}^{1}\frac{\varepsilon \,\partial _{s^{\alpha }}G(x,t)}{%
\left( \left\vert u_{\varepsilon }^{\alpha }+th\triangle
_{h,s}u_{\varepsilon }^{\alpha }\right\vert ^{2}+\varepsilon \right) ^{\frac{%
3}{2}}}\,dt\,\triangle _{h,s}u_{\varepsilon }^{\alpha }$%
\end{tabular}%
\end{equation}%
where, with an abuse of notation, we have denoted%
\begin{equation*}
\begin{tabular}{l}
$\partial _{s^{\beta }s^{\alpha }}G\left( x,t\right) =\partial _{s^{\beta
}s^{\alpha }}G\left( \sqrt{\left\vert s_{\varepsilon }^{1}\left( x,t\right)
\right\vert ^{2}+\varepsilon },...,\sqrt{\left\vert s_{\varepsilon
}^{m}\left( x,t\right) \right\vert ^{2}+\varepsilon },\sqrt{\left\vert \xi
_{\varepsilon }^{1}\left( x,t\right) \right\vert ^{2}+\varepsilon },...,%
\sqrt{\left\vert \xi _{\varepsilon }^{m}\left( x,t\right) \right\vert
^{2}+\varepsilon }\right) $ \\ 
\\ 
$\partial _{\xi ^{\beta }s^{\alpha }}G\left( x,t\right) =\partial _{\xi
^{\beta }s^{\alpha }}G\left( \sqrt{\left\vert s_{\varepsilon }^{1}\left(
x,t\right) \right\vert ^{2}+\varepsilon },...,\sqrt{\left\vert
s_{\varepsilon }^{m}\left( x,t\right) \right\vert ^{2}+\varepsilon },\sqrt{%
\left\vert \xi _{\varepsilon }^{1}\left( x,t\right) \right\vert
^{2}+\varepsilon },...,\sqrt{\left\vert \xi _{\varepsilon }^{m}\left(
x,t\right) \right\vert ^{2}+\varepsilon }\right) $ \\ 
\\ 
$\partial _{s^{\alpha }}G(x,t)=\partial _{s^{\alpha }}G\left( \sqrt{%
\left\vert s_{\varepsilon }^{1}\left( x,t\right) \right\vert
^{2}+\varepsilon },...,\sqrt{\left\vert s_{\varepsilon }^{m}\left(
x,t\right) \right\vert ^{2}+\varepsilon },\sqrt{\left\vert \xi _{\varepsilon
}^{1}\left( x,t\right) \right\vert ^{2}+\varepsilon },...,\sqrt{\left\vert
\xi _{\varepsilon }^{m}\left( x,t\right) \right\vert ^{2}+\varepsilon }%
\right) $%
\end{tabular}%
\end{equation*}%
where 
\begin{equation*}
\begin{tabular}{l}
$s_{\varepsilon }^{\alpha }\left( x,t\right) =u_{\varepsilon }^{\alpha
}+th\triangle _{h,s}u_{\varepsilon }^{\alpha }$ \\ 
\\ 
$\xi _{\varepsilon }^{\alpha }\left( x,t\right) =\nabla u_{\varepsilon
}^{\alpha }+th\triangle _{h,s}\nabla u_{\varepsilon }^{\alpha }$%
\end{tabular}%
\end{equation*}%
with $\alpha =1,...,m$. From the previous relation (4.8) we have%
\begin{equation*}
\begin{tabular}{l}
$\int\limits_{\Omega }\triangle _{h,s}\left[ \partial _{s^{\alpha }}G\left( 
\sqrt{\left\vert u_{\varepsilon }^{1}\right\vert ^{2}+\varepsilon },...,%
\sqrt{\left\vert u_{\varepsilon }^{m}\right\vert ^{2}+\varepsilon },\sqrt{%
\left\vert \nabla u_{\varepsilon }^{1}\right\vert ^{2}+\varepsilon },...,%
\sqrt{\left\vert \nabla u_{\varepsilon }^{m}\right\vert ^{2}+\varepsilon }%
\right) \frac{u_{\varepsilon }^{\alpha }}{\sqrt{\left\vert u_{\varepsilon
}^{\alpha }\right\vert ^{2}+\varepsilon }}\right] $ \\ 
$\cdot \left( \eta ^{2}\triangle _{h,s}u_{\varepsilon }^{\alpha }\right)
\,dx $ \\ 
$=\int\limits_{\Omega }\sum\limits_{\beta =1}^{m}\int\limits_{0}^{1}\partial
_{s^{\beta }s^{\alpha }}G\left( x,t\right) \frac{u_{\varepsilon }^{\beta
}+th\triangle _{h,s}u_{\varepsilon }^{\beta }}{\sqrt{\left\vert
u_{\varepsilon }^{\beta }+th\triangle _{h,s}u_{\varepsilon }^{\beta
}\right\vert ^{2}+\varepsilon }}\frac{u_{\varepsilon }^{\alpha }+th\triangle
_{h,s}u_{\varepsilon }^{\alpha }}{\sqrt{\left\vert u_{\varepsilon }^{\alpha
}+th\triangle _{h,s}u_{\varepsilon }^{\alpha }\right\vert ^{2}+\varepsilon }}%
\,dt$ \\ 
$\cdot \,\triangle _{h,s}u_{\varepsilon }^{\beta }\,\eta ^{2}\triangle
_{h,s}u_{\varepsilon }^{\alpha }\,dx$ \\ 
$+\int\limits_{\Omega }\sum\limits_{\beta =1}^{m}\int\limits_{0}^{1}\partial
_{\xi ^{\beta }s^{\alpha }}G\left( x,t\right) \frac{\nabla u_{\varepsilon
}^{\beta }+th\triangle _{h,s}\nabla u_{\varepsilon }^{\beta }}{\sqrt{%
\left\vert \nabla u_{\varepsilon }^{\beta }+th\triangle _{h,s}\nabla
u_{\varepsilon }^{\beta }\right\vert ^{2}+\varepsilon }}\frac{u_{\varepsilon
}^{\alpha }+th\triangle _{h,s}u_{\varepsilon }^{\alpha }}{\sqrt{\left\vert
u_{\varepsilon }^{\alpha }+th\triangle _{h,s}u_{\varepsilon }^{\alpha
}\right\vert ^{2}+\varepsilon }}\,dt\,$ \\ 
$\cdot \triangle _{h,s}\nabla u_{\varepsilon }^{\beta }\,\eta ^{2}\triangle
_{h,s}u_{\varepsilon }^{\alpha }\,dx$ \\
\end{tabular}%
\end{equation*}%
using Fubini's theorem we have%
\begin{equation}
\begin{tabular}{l}
$\int\limits_{\Omega }\triangle _{h,s}\left[ \partial _{s^{\alpha }}G\left( 
\sqrt{\left\vert u_{\varepsilon }^{1}\right\vert ^{2}+\varepsilon },...,%
\sqrt{\left\vert u_{\varepsilon }^{m}\right\vert ^{2}+\varepsilon },\sqrt{%
\left\vert \nabla u_{\varepsilon }^{1}\right\vert ^{2}+\varepsilon },...,%
\sqrt{\left\vert \nabla u_{\varepsilon }^{m}\right\vert ^{2}+\varepsilon }%
\right) \frac{u_{\varepsilon }^{\alpha }}{\sqrt{\left\vert u_{\varepsilon
}^{\alpha }\right\vert ^{2}+\varepsilon }}\right] $ \\ 
$\cdot \left( \eta ^{2}\triangle _{h,s}u_{\varepsilon }^{\alpha }\right)
\,dx $ \\ 
$=\int\limits_{0}^{1}\int\limits_{\Omega }\sum\limits_{\beta =1}^{m}\partial
_{s^{\beta }s^{\alpha }}G\left( x,t\right) \frac{u_{\varepsilon }^{\beta
}+th\triangle _{h,s}u_{\varepsilon }^{\beta }}{\sqrt{\left\vert
u_{\varepsilon }^{\beta }+th\triangle _{h,s}u_{\varepsilon }^{\beta
}\right\vert ^{2}+\varepsilon }}\frac{u_{\varepsilon }^{\alpha }+th\triangle
_{h,s}u_{\varepsilon }^{\alpha }}{\sqrt{\left\vert u_{\varepsilon }^{\alpha
}+th\triangle _{h,s}u_{\varepsilon }^{\alpha }\right\vert ^{2}+\varepsilon }}%
\,$ \\ 
$\cdot \,\triangle _{h,s}u_{\varepsilon }^{\beta }\,\triangle
_{h,s}u_{\varepsilon }^{\alpha }\,\eta ^{2}\,dx\,dt$ \\ 
$+\int\limits_{0}^{1}\int\limits_{\Omega }\sum\limits_{\beta =1}^{m}\partial
_{\xi ^{\beta }s^{\alpha }}G\left( x,t\right) \frac{\nabla u_{\varepsilon
}^{\beta }+th\triangle _{h,s}\nabla u_{\varepsilon }^{\beta }}{\sqrt{%
\left\vert \nabla u_{\varepsilon }^{\beta }+th\triangle _{h,s}\nabla
u_{\varepsilon }^{\beta }\right\vert ^{2}+\varepsilon }}\frac{u_{\varepsilon
}^{\alpha }+th\triangle _{h,s}u_{\varepsilon }^{\alpha }}{\sqrt{\left\vert
u_{\varepsilon }^{\alpha }+th\triangle _{h,s}u_{\varepsilon }^{\alpha
}\right\vert ^{2}+\varepsilon }}\,$ \\ 
$\cdot \,\triangle _{h,s}\nabla u_{\varepsilon }^{\beta }\,\triangle
_{h,s}u_{\varepsilon }^{\alpha }\,\eta ^{2}\,dx\,dt$ \\
\end{tabular}%
\end{equation}%
Now it remains to estimate the term defined by the relation in (4.5)%
\begin{equation*}
\begin{tabular}{l}
$\triangle _{h,s}\left[ \partial _{\xi ^{\alpha }}G\left( \sqrt{\left\vert
u_{\varepsilon }^{1}\right\vert ^{2}+\varepsilon },...,\sqrt{\left\vert
u_{\varepsilon }^{m}\right\vert ^{2}+\varepsilon },\sqrt{\left\vert \nabla
u_{\varepsilon }^{1}\right\vert ^{2}+\varepsilon },...,\sqrt{\left\vert
\nabla u_{\varepsilon }^{m}\right\vert ^{2}+\varepsilon }\right) \frac{%
\nabla u_{\varepsilon }^{\alpha }}{\sqrt{\left\vert \nabla u_{\varepsilon
}^{\alpha }\right\vert ^{2}+\varepsilon }}\right] $ \\ 
$=\frac{1}{h}\left[ \partial _{\xi ^{\alpha }}G\left( \sqrt{\left\vert
u_{\varepsilon }^{1}+th\triangle _{h,s}u_{\varepsilon }^{1}\right\vert
^{2}+\varepsilon },...,\sqrt{\left\vert u_{\varepsilon }^{m}+th\triangle
_{h,s}u_{\varepsilon }^{m}\right\vert ^{2}+\varepsilon }\right. \right. $ \\ 
$\left. \left. \sqrt{\left\vert \nabla u_{\varepsilon }^{1}+th\triangle
_{h,s}\nabla u_{\varepsilon }^{1}\right\vert ^{2}+\varepsilon },...,\sqrt{%
\left\vert \nabla u_{\varepsilon }^{m}+th\triangle _{h,s}\nabla
u_{\varepsilon }^{m}\right\vert ^{2}+\varepsilon }\right) \frac{\nabla
u_{\varepsilon }^{\alpha }+th\triangle _{h,s}\nabla u_{\varepsilon }^{\alpha
}}{\sqrt{\left\vert \nabla u_{\varepsilon }^{\alpha }+th\triangle
_{h,s}\nabla u_{\varepsilon }^{\alpha }\right\vert ^{2}+\varepsilon }}\right]
_{0}^{1}$ \\ 
$=\frac{1}{h}\int\limits_{0}^{1}\frac{d}{dt}\left( \partial _{\xi ^{\alpha
}}G\left( \sqrt{\left\vert u_{\varepsilon }^{1}+th\triangle
_{h,s}u_{\varepsilon }^{1}\right\vert ^{2}+\varepsilon },...,\sqrt{%
\left\vert u_{\varepsilon }^{m}+th\triangle _{h,s}u_{\varepsilon
}^{m}\right\vert ^{2}+\varepsilon }\right. \right. $ \\
\end{tabular}%
\end{equation*}%
since $G\in C^{2}\left( 
\mathbb{R}
_{+,0}^{m}\times 
\mathbb{R}
_{+,0}^{m}\right) $ we obtain%
\begin{equation}
\begin{tabular}{l}
$\triangle _{h,s}\left[ \partial _{\xi ^{\alpha }}G\left( \sqrt{\left\vert
u_{\varepsilon }^{1}\right\vert ^{2}+\varepsilon },...,\sqrt{\left\vert
u_{\varepsilon }^{m}\right\vert ^{2}+\varepsilon },\sqrt{\left\vert \nabla
u_{\varepsilon }^{1}\right\vert ^{2}+\varepsilon },...,\sqrt{\left\vert
\nabla u_{\varepsilon }^{m}\right\vert ^{2}+\varepsilon }\right) \frac{%
\nabla u_{\varepsilon }^{\alpha }}{\sqrt{\left\vert \nabla u_{\varepsilon
}^{\alpha }\right\vert ^{2}+\varepsilon }}\right] $ \\ 
$=\frac{1}{h}\int\limits_{0}^{1}\sum\limits_{\beta =1}^{m}\partial
_{s^{\beta }\xi ^{\alpha }}G\left( x,t\right) \frac{u_{\varepsilon }^{\beta
}+th\triangle _{h,s}u_{\varepsilon }^{\beta }}{\sqrt{\left\vert u^{\beta
}+th\triangle _{h,s}u^{\beta }\right\vert ^{2}+\varepsilon }}\frac{\nabla
u_{\varepsilon }^{\alpha }+th\triangle _{h,s}\nabla u_{\varepsilon }^{\alpha
}}{\sqrt{\left\vert \nabla u^{\alpha }+th\triangle _{h,s}\nabla u^{\alpha
}\right\vert ^{2}+\varepsilon }}\,\,h\triangle _{h,s}u_{\varepsilon }^{\beta
}$ \\ 
$+\sum\limits_{\beta =1}^{m}\partial _{\xi ^{\beta }\xi ^{\alpha }}G\left(
x,t\right) \frac{\nabla u_{\varepsilon }^{\beta }+th\triangle _{h,s}\nabla
u_{\varepsilon }^{\beta }}{\sqrt{\left\vert \nabla u_{\varepsilon }^{\beta
}+th\triangle _{h,s}\nabla u_{\varepsilon }^{\beta }\right\vert
^{2}+\varepsilon }}\frac{\nabla u_{\varepsilon }^{\alpha }+th\triangle
_{h,s}\nabla u_{\varepsilon }^{\alpha }}{\sqrt{\left\vert \nabla
u_{\varepsilon }^{\alpha }+th\triangle _{h,s}\nabla u_{\varepsilon }^{\alpha
}\right\vert ^{2}+\varepsilon }}\,\,h\triangle _{h,s}\nabla u_{\varepsilon
}^{\beta }$ \\ 
$+\frac{\varepsilon \,\partial _{\xi ^{\alpha }}G(x,t)}{\left( \left\vert
\nabla u_{\varepsilon }^{\alpha }+th\triangle _{h,s}\nabla u_{\varepsilon
}^{\alpha }\right\vert ^{2}+\varepsilon \right) ^{\frac{3}{2}}}%
\,h\,\triangle _{h,s}\nabla u_{\varepsilon }^{\alpha }\,dt$%
\end{tabular}%
\end{equation}%
where we denoted%
\begin{equation*}
\begin{tabular}{l}
$\partial _{s^{\beta }\xi ^{\alpha }}G\left( x,t\right) =\partial _{s^{\beta
}\xi ^{\alpha }}G\left( \sqrt{\left\vert s_{\varepsilon }^{1}\left(
x,t\right) \right\vert ^{2}+\varepsilon },...,\sqrt{\left\vert
s_{\varepsilon }^{m}\left( x,t\right) \right\vert ^{2}+\varepsilon },\sqrt{%
\left\vert \xi _{\varepsilon }^{1}\left( x,t\right) \right\vert
^{2}+\varepsilon },...,\sqrt{\left\vert \xi _{\varepsilon }^{m}\left(
x,t\right) \right\vert ^{2}+\varepsilon }\right) $ \\ 
\\ 
$\partial _{\xi ^{\beta }\xi ^{\alpha }}G\left( x,t\right) =\partial _{\xi
^{\beta }\xi ^{\alpha }}G\left( \sqrt{\left\vert s_{\varepsilon }^{1}\left(
x,t\right) \right\vert ^{2}+\varepsilon },...,\sqrt{\left\vert
s_{\varepsilon }^{m}\left( x,t\right) \right\vert ^{2}+\varepsilon },\sqrt{%
\left\vert \xi _{\varepsilon }^{1}\left( x,t\right) \right\vert
^{2}+\varepsilon },...,\sqrt{\left\vert \xi _{\varepsilon }^{m}\left(
x,t\right) \right\vert ^{2}+\varepsilon }\right) $ \\ 
\\ 
$\partial _{\xi ^{\alpha }}G(x,t)=\partial _{\xi ^{\alpha }}G\left( \sqrt{%
\left\vert s_{\varepsilon }^{1}\left( x,t\right) \right\vert
^{2}+\varepsilon },...,\sqrt{\left\vert s_{\varepsilon }^{m}\left(
x,t\right) \right\vert ^{2}+\varepsilon },\sqrt{\left\vert \xi _{\varepsilon
}^{1}\left( x,t\right) \right\vert ^{2}+\varepsilon },...,\sqrt{\left\vert
\xi _{\varepsilon }^{m}\left( x,t\right) \right\vert ^{2}+\varepsilon }%
\right) $%
\end{tabular}%
\end{equation*}%
where 
\begin{equation*}
\begin{tabular}{l}
$s_{\varepsilon }^{\alpha }\left( x,t\right) =u_{\varepsilon }^{\alpha
}+th\triangle _{h,s}u_{\varepsilon }^{\alpha }$ \\ 
\\ 
$\xi _{\varepsilon }^{\alpha }\left( x,t\right) =\nabla u_{\varepsilon
}^{\alpha }+th\triangle _{h,s}\nabla u_{\varepsilon }^{\alpha }$%
\end{tabular}%
\end{equation*}%
with $\alpha =1,...,m$. Using (4.10) we obtain%
\begin{equation*}
\begin{tabular}{l}
$\triangle _{h,s}\left[ \partial _{\xi ^{\alpha }}G\left( \sqrt{\left\vert
u_{\varepsilon }^{1}\right\vert ^{2}+\varepsilon },...,\sqrt{\left\vert
u_{\varepsilon }^{m}\right\vert ^{2}+\varepsilon },\sqrt{\left\vert \nabla
u_{\varepsilon }^{1}\right\vert ^{2}+\varepsilon },...,\sqrt{\left\vert
\nabla u_{\varepsilon }^{m}\right\vert ^{2}+\varepsilon }\right) \frac{%
\nabla u_{\varepsilon }^{\alpha }}{\sqrt{\left\vert \nabla u_{\varepsilon
}^{\alpha }\right\vert ^{2}+\varepsilon }}\right] $ \\ 
$=\int\limits_{0}^{1}\sum\limits_{\beta =1}^{m}\partial _{s^{\beta }\xi
^{\alpha }}G\left( x,t\right) \frac{u_{\varepsilon }^{\beta }+th\triangle
_{h,s}u_{\varepsilon }^{\beta }}{\sqrt{\left\vert u^{\beta }+th\triangle
_{h,s}u^{\beta }\right\vert ^{2}+\varepsilon }}\frac{\nabla u_{\varepsilon
}^{\alpha }+th\triangle _{h,s}\nabla u_{\varepsilon }^{\alpha }}{\sqrt{%
\left\vert \nabla u^{\alpha }+th\triangle _{h,s}\nabla u^{\alpha
}\right\vert ^{2}+\varepsilon }}\,\,\triangle _{h,s}u_{\varepsilon }^{\beta
}\,dt$ \\ 
$+\int\limits_{0}^{1}\sum\limits_{\beta =1}^{m}\partial _{\xi ^{\beta }\xi
^{\alpha }}G\left( x,t\right) \frac{\nabla u_{\varepsilon }^{\beta
}+th\triangle _{h,s}\nabla u_{\varepsilon }^{\beta }}{\sqrt{\left\vert
\nabla u_{\varepsilon }^{\beta }+th\triangle _{h,s}\nabla u_{\varepsilon
}^{\beta }\right\vert ^{2}+\varepsilon }}\frac{\nabla u_{\varepsilon
}^{\alpha }+th\triangle _{h,s}\nabla u_{\varepsilon }^{\alpha }}{\sqrt{%
\left\vert \nabla u_{\varepsilon }^{\alpha }+th\triangle _{h,s}\nabla
u_{\varepsilon }^{\alpha }\right\vert ^{2}+\varepsilon }}\,\,\triangle
_{h,s}\nabla u_{\varepsilon }^{\beta }\,dt$ \\ 
$+\int\limits_{0}^{1}\frac{\varepsilon \,\partial _{\xi ^{\alpha }}G(x,t)}{%
\left( \left\vert \nabla u_{\varepsilon }^{\alpha }+th\triangle _{h,s}\nabla
u_{\varepsilon }^{\alpha }\right\vert ^{2}+\varepsilon \right) ^{\frac{3}{2}}%
}\,\,\triangle _{h,s}\nabla u_{\varepsilon }^{\alpha }\,dt$%
\end{tabular}%
\end{equation*}%
it follows that%
\begin{equation*}
\begin{tabular}{l}
$\int\limits_{\Omega }\triangle _{h,s}\left[ \partial _{\xi ^{\alpha
}}G\left( \sqrt{\left\vert u_{\varepsilon }^{1}\right\vert ^{2}+\varepsilon }%
,...,\sqrt{\left\vert u_{\varepsilon }^{m}\right\vert ^{2}+\varepsilon },%
\sqrt{\left\vert \nabla u_{\varepsilon }^{1}\right\vert ^{2}+\varepsilon }%
,...,\sqrt{\left\vert \nabla u_{\varepsilon }^{m}\right\vert
^{2}+\varepsilon }\right) \frac{\nabla u_{\varepsilon }^{\alpha }}{\sqrt{%
\left\vert \nabla u_{\varepsilon }^{\alpha }\right\vert ^{2}+\varepsilon }}%
\right] $ \\ 
$\cdot \nabla \left( \eta ^{2}\triangle _{h,s}u_{\varepsilon }^{\alpha
}\right) \,dx$ \\ 
$=\int\limits_{\Omega }\,\int\limits_{0}^{1}\sum\limits_{\beta
=1}^{m}\partial _{s^{\beta }\xi ^{\alpha }}G\left( x,t\right) \frac{%
u_{\varepsilon }^{\beta }+th\triangle _{h,s}u_{\varepsilon }^{\beta }}{\sqrt{%
\left\vert u^{\beta }+th\triangle _{h,s}u^{\beta }\right\vert
^{2}+\varepsilon }}\frac{\nabla u_{\varepsilon }^{\alpha }+th\triangle
_{h,s}\nabla u_{\varepsilon }^{\alpha }}{\sqrt{\left\vert \nabla u^{\alpha
}+th\triangle _{h,s}\nabla u^{\alpha }\right\vert ^{2}+\varepsilon }}%
\,\,\triangle _{h,s}u_{\varepsilon }^{\beta }\,dt$ \\ 
$\cdot \nabla \left( \eta ^{2}\triangle _{h,s}u_{\varepsilon }^{\alpha
}\right) \,dx$ \\ 
$+\int\limits_{\Omega }\,\int\limits_{0}^{1}\sum\limits_{\beta
=1}^{m}\partial _{\xi ^{\beta }\xi ^{\alpha }}G\left( x,t\right) \frac{%
\nabla u_{\varepsilon }^{\beta }+th\triangle _{h,s}\nabla u_{\varepsilon
}^{\beta }}{\sqrt{\left\vert \nabla u_{\varepsilon }^{\beta }+th\triangle
_{h,s}\nabla u_{\varepsilon }^{\beta }\right\vert ^{2}+\varepsilon }}\frac{%
\nabla u_{\varepsilon }^{\alpha }+th\triangle _{h,s}\nabla u_{\varepsilon
}^{\alpha }}{\sqrt{\left\vert \nabla u_{\varepsilon }^{\alpha }+th\triangle
_{h,s}\nabla u_{\varepsilon }^{\alpha }\right\vert ^{2}+\varepsilon }}%
\,\,\triangle _{h,s}\nabla u_{\varepsilon }^{\beta }\,dt$ \\ 
$\cdot \nabla \left( \eta ^{2}\triangle _{h,s}u_{\varepsilon }^{\alpha
}\right) \,dx$ \\
\end{tabular}%
\end{equation*}%
Since $\nabla \left( \eta ^{2}\triangle _{h,s}u_{\varepsilon }^{\alpha
}\right) =\eta ^{2}\triangle _{h,s}\nabla u_{\varepsilon }^{\alpha }+2\eta
\nabla \eta \triangle _{h,s}u_{\varepsilon }^{\alpha }$ then we get%
\begin{equation*}
\begin{tabular}{l}
$\int\limits_{\Omega }\triangle _{h,s}\left[ \partial _{\xi ^{\alpha
}}G\left( \sqrt{\left\vert u_{\varepsilon }^{1}\right\vert ^{2}+\varepsilon }%
,...,\sqrt{\left\vert u_{\varepsilon }^{m}\right\vert ^{2}+\varepsilon },%
\sqrt{\left\vert \nabla u_{\varepsilon }^{1}\right\vert ^{2}+\varepsilon }%
,...,\sqrt{\left\vert \nabla u_{\varepsilon }^{m}\right\vert
^{2}+\varepsilon }\right) \frac{\nabla u_{\varepsilon }^{\alpha }}{\sqrt{%
\left\vert \nabla u_{\varepsilon }^{\alpha }\right\vert ^{2}+\varepsilon }}%
\right] $ \\ 
$\cdot \nabla \left( \eta ^{2}\triangle _{h,s}u_{\varepsilon }^{\alpha
}\right) \,dx$ \\ 
$=\int\limits_{\Omega }\,\int\limits_{0}^{1}\sum\limits_{\beta
=1}^{m}\partial _{s^{\beta }\xi ^{\alpha }}G\left( x,t\right) \frac{%
u_{\varepsilon }^{\beta }+th\triangle _{h,s}u_{\varepsilon }^{\beta }}{\sqrt{%
\left\vert u^{\beta }+th\triangle _{h,s}u^{\beta }\right\vert
^{2}+\varepsilon }}\frac{\nabla u_{\varepsilon }^{\alpha }+th\triangle
_{h,s}\nabla u_{\varepsilon }^{\alpha }}{\sqrt{\left\vert \nabla u^{\alpha
}+th\triangle _{h,s}\nabla u^{\alpha }\right\vert ^{2}+\varepsilon }}%
\,\,\triangle _{h,s}u_{\varepsilon }^{\beta }\,dt$ \\ 
$\cdot \eta ^{2}\triangle _{h,s}\nabla u_{\varepsilon }^{\alpha }\,dx$ \\ 
$+2\int\limits_{\Omega }\,\int\limits_{0}^{1}\sum\limits_{\beta
=1}^{m}\partial _{s^{\beta }\xi ^{\alpha }}G\left( x,t\right) \frac{%
u_{\varepsilon }^{\beta }+th\triangle _{h,s}u_{\varepsilon }^{\beta }}{\sqrt{%
\left\vert u^{\beta }+th\triangle _{h,s}u^{\beta }\right\vert
^{2}+\varepsilon }}\frac{\nabla u_{\varepsilon }^{\alpha }+th\triangle
_{h,s}\nabla u_{\varepsilon }^{\alpha }}{\sqrt{\left\vert \nabla u^{\alpha
}+th\triangle _{h,s}\nabla u^{\alpha }\right\vert ^{2}+\varepsilon }}%
\,\,\triangle _{h,s}u_{\varepsilon }^{\beta }\,dt$ \\ 
$\cdot \eta \nabla \eta \triangle _{h,s}u_{\varepsilon }^{\alpha }\,dx$ \\
\end{tabular}%
\end{equation*}

using Fubini's theorem we have%
\begin{equation}
\begin{tabular}{l}
$\int\limits_{\Omega }\triangle _{h,s}\left[ \partial _{\xi ^{\alpha
}}G\left( \sqrt{\left\vert u_{\varepsilon }^{1}\right\vert ^{2}+\varepsilon }%
,...,\sqrt{\left\vert u_{\varepsilon }^{m}\right\vert ^{2}+\varepsilon },%
\sqrt{\left\vert \nabla u_{\varepsilon }^{1}\right\vert ^{2}+\varepsilon }%
,...,\sqrt{\left\vert \nabla u_{\varepsilon }^{m}\right\vert
^{2}+\varepsilon }\right) \frac{\nabla u_{\varepsilon }^{\alpha }}{\sqrt{%
\left\vert \nabla u_{\varepsilon }^{\alpha }\right\vert ^{2}+\varepsilon }}%
\right] $ \\ 
$\cdot \nabla \left( \eta ^{2}\triangle _{h,s}u_{\varepsilon }^{\alpha
}\right) \,dx$ \\ 
$=\,\int\limits_{0}^{1}\,\int\limits_{\Omega }\,\sum\limits_{\beta
=1}^{m}\partial _{s^{\beta }\xi ^{\alpha }}G\left( x,t\right) \frac{%
u_{\varepsilon }^{\beta }+th\triangle _{h,s}u_{\varepsilon }^{\beta }}{\sqrt{%
\left\vert u^{\beta }+th\triangle _{h,s}u^{\beta }\right\vert
^{2}+\varepsilon }}\frac{\nabla u_{\varepsilon }^{\alpha }+th\triangle
_{h,s}\nabla u_{\varepsilon }^{\alpha }}{\sqrt{\left\vert \nabla u^{\alpha
}+th\triangle _{h,s}\nabla u^{\alpha }\right\vert ^{2}+\varepsilon }}%
\,\,\triangle _{h,s}u_{\varepsilon }^{\beta }\,$ \\ 
$\cdot \eta ^{2}\triangle _{h,s}\nabla u_{\varepsilon }^{\alpha }\,dx\,dt$
\\ 
$+2\,\int\limits_{0}^{1}\,\int\limits_{\Omega }\,\sum\limits_{\beta
=1}^{m}\partial _{s^{\beta }\xi ^{\alpha }}G\left( x,t\right) \frac{%
u_{\varepsilon }^{\beta }+th\triangle _{h,s}u_{\varepsilon }^{\beta }}{\sqrt{%
\left\vert u^{\beta }+th\triangle _{h,s}u^{\beta }\right\vert
^{2}+\varepsilon }}\frac{\nabla u_{\varepsilon }^{\alpha }+th\triangle
_{h,s}\nabla u_{\varepsilon }^{\alpha }}{\sqrt{\left\vert \nabla u^{\alpha
}+th\triangle _{h,s}\nabla u^{\alpha }\right\vert ^{2}+\varepsilon }}%
\,\,\triangle _{h,s}u_{\varepsilon }^{\beta }\,$ \\ 
$\cdot \eta \nabla \eta \triangle _{h,s}u_{\varepsilon }^{\alpha }\,dx\,dt$
\\
\end{tabular}%
\end{equation}%
Using (4.2), (4.6), (4.9) and (4.11) we get%
\begin{equation}
\begin{tabular}{l}
$0$ \\ 
$=\left( p-1\right) \int\limits_{\Omega }\varepsilon \,I_{1,\alpha }\,\left(
\triangle _{h,s}\left( \nabla u_{\varepsilon }^{\alpha }\right) \right)
^{2}\eta ^{2}dx$ \\ 
$+2\left( p-1\right) \int\limits_{\Omega }\varepsilon \,\triangle
_{h,s}\left( \nabla u_{\varepsilon }^{\alpha }\right) \left( \triangle
_{h,s}u_{\varepsilon }^{\alpha }\right) I_{1,\alpha }\,\eta \,\nabla \eta
\,dx$ \\ 
$+\int\limits_{0}^{1}\int\limits_{\Omega }\sum\limits_{\beta =1}^{m}\partial
_{s^{\beta }s^{\alpha }}G\left( x,t\right) \frac{u_{\varepsilon }^{\beta
}+th\triangle _{h,s}u_{\varepsilon }^{\beta }}{\sqrt{\left\vert
u_{\varepsilon }^{\beta }+th\triangle _{h,s}u_{\varepsilon }^{\beta
}\right\vert ^{2}+\varepsilon }}\frac{u_{\varepsilon }^{\alpha }+th\triangle
_{h,s}u_{\varepsilon }^{\alpha }}{\sqrt{\left\vert u_{\varepsilon }^{\alpha
}+th\triangle _{h,s}u_{\varepsilon }^{\alpha }\right\vert ^{2}+\varepsilon }}%
\,$ \\ 
$\cdot \,\triangle _{h,s}u_{\varepsilon }^{\beta }\,\triangle
_{h,s}u_{\varepsilon }^{\alpha }\,\eta ^{2}\,dx\,dt$ \\ 
$+\int\limits_{0}^{1}\int\limits_{\Omega }\sum\limits_{\beta =1}^{m}\partial
_{\xi ^{\beta }s^{\alpha }}G\left( x,t\right) \frac{\nabla u_{\varepsilon
}^{\beta }+th\triangle _{h,s}\nabla u_{\varepsilon }^{\beta }}{\sqrt{%
\left\vert \nabla u_{\varepsilon }^{\beta }+th\triangle _{h,s}\nabla
u_{\varepsilon }^{\beta }\right\vert ^{2}+\varepsilon }}\frac{u_{\varepsilon
}^{\alpha }+th\triangle _{h,s}u_{\varepsilon }^{\alpha }}{\sqrt{\left\vert
u_{\varepsilon }^{\alpha }+th\triangle _{h,s}u_{\varepsilon }^{\alpha
}\right\vert ^{2}+\varepsilon }}\,$ \\ 
$\cdot \,\triangle _{h,s}\nabla u_{\varepsilon }^{\beta }\,\triangle
_{h,s}u_{\varepsilon }^{\alpha }\,\eta ^{2}\,dx\,dt$ \\
\end{tabular}%
\end{equation}%
The previous equation (4.12) can be rewritten as follows%
\begin{equation}
A_{\alpha }=-B_{\alpha }
\end{equation}%
where%
\begin{equation}
\begin{tabular}{l}
$A_{\alpha }$ \\ 
$=\left( p-1\right) \int\limits_{\Omega }\varepsilon \,I_{1,\alpha }\,\left(
\triangle _{h,s}\left( \nabla u_{\varepsilon }^{\alpha }\right) \right)
^{2}\eta ^{2}dx$ \\ 
$+\int\limits_{0}^{1}\int\limits_{\Omega }\sum\limits_{\beta =1}^{m}\partial
_{s^{\beta }s^{\alpha }}G\left( x,t\right) \frac{u_{\varepsilon }^{\beta
}+th\triangle _{h,s}u_{\varepsilon }^{\beta }}{\sqrt{\left\vert
u_{\varepsilon }^{\beta }+th\triangle _{h,s}u_{\varepsilon }^{\beta
}\right\vert ^{2}+\varepsilon }}\frac{u_{\varepsilon }^{\alpha }+th\triangle
_{h,s}u_{\varepsilon }^{\alpha }}{\sqrt{\left\vert u_{\varepsilon }^{\alpha
}+th\triangle _{h,s}u_{\varepsilon }^{\alpha }\right\vert ^{2}+\varepsilon }}%
\,$ \\ 
$\cdot \,\triangle _{h,s}u_{\varepsilon }^{\beta }\,\triangle
_{h,s}u_{\varepsilon }^{\alpha }\,\eta ^{2}\,dx\,dt$ \\ 
$+\int\limits_{0}^{1}\int\limits_{\Omega }\sum\limits_{\beta =1}^{m}\partial
_{\xi ^{\beta }s^{\alpha }}G\left( x,t\right) \frac{\nabla u_{\varepsilon
}^{\beta }+th\triangle _{h,s}\nabla u_{\varepsilon }^{\beta }}{\sqrt{%
\left\vert \nabla u_{\varepsilon }^{\beta }+th\triangle _{h,s}\nabla
u_{\varepsilon }^{\beta }\right\vert ^{2}+\varepsilon }}\frac{u_{\varepsilon
}^{\alpha }+th\triangle _{h,s}u_{\varepsilon }^{\alpha }}{\sqrt{\left\vert
u_{\varepsilon }^{\alpha }+th\triangle _{h,s}u_{\varepsilon }^{\alpha
}\right\vert ^{2}+\varepsilon }}\,$ \\ 
$\cdot \,\triangle _{h,s}\nabla u_{\varepsilon }^{\beta }\,\triangle
_{h,s}u_{\varepsilon }^{\alpha }\,\eta ^{2}\,dx\,dt$ \\
\end{tabular}%
\end{equation}%
and%
\begin{equation}
\begin{tabular}{l}
$B_{\alpha }$ \\ 
$=2\left( p-1\right) \int\limits_{\Omega }\varepsilon \,\triangle
_{h,s}\left( \nabla u_{\varepsilon }^{\alpha }\right) \left( \triangle
_{h,s}u_{\varepsilon }^{\alpha }\right) I_{1,\alpha }\,\eta \,\nabla \eta
\,dx$ \\ 
$+2\,\int\limits_{0}^{1}\,\int\limits_{\Omega }\,\sum\limits_{\beta
=1}^{m}\partial _{s^{\beta }\xi ^{\alpha }}G\left( x,t\right) \frac{%
u_{\varepsilon }^{\beta }+th\triangle _{h,s}u_{\varepsilon }^{\beta }}{\sqrt{%
\left\vert u^{\beta }+th\triangle _{h,s}u^{\beta }\right\vert
^{2}+\varepsilon }}\frac{\nabla u_{\varepsilon }^{\alpha }+th\triangle
_{h,s}\nabla u_{\varepsilon }^{\alpha }}{\sqrt{\left\vert \nabla u^{\alpha
}+th\triangle _{h,s}\nabla u^{\alpha }\right\vert ^{2}+\varepsilon }}%
\,\,\triangle _{h,s}u_{\varepsilon }^{\beta }\,$ \\ 
$\cdot \eta \nabla \eta \triangle _{h,s}u_{\varepsilon }^{\alpha }\,dx\,dt$
\\ 
$+2\int\limits_{0}^{1}\int\limits_{\Omega }\sum\limits_{\beta
=1}^{m}\partial _{\xi ^{\beta }\xi ^{\alpha }}G\left( x,t\right) \frac{%
\nabla u_{\varepsilon }^{\beta }+th\triangle _{h,s}\nabla u_{\varepsilon
}^{\beta }}{\sqrt{\left\vert \nabla u_{\varepsilon }^{\beta }+th\triangle
_{h,s}\nabla u_{\varepsilon }^{\beta }\right\vert ^{2}+\varepsilon }}\frac{%
\nabla u_{\varepsilon }^{\alpha }+th\triangle _{h,s}\nabla u_{\varepsilon
}^{\alpha }}{\sqrt{\left\vert \nabla u_{\varepsilon }^{\alpha }+th\triangle
_{h,s}\nabla u_{\varepsilon }^{\alpha }\right\vert ^{2}+\varepsilon }}%
\,\,\triangle _{h,s}\nabla u_{\varepsilon }^{\beta }\,$ \\ 
$\cdot \eta \nabla \eta \triangle _{h,s}u_{\varepsilon }^{\alpha }\,dx\,dt$
\\ 
$+2\int\limits_{0}^{1}\,\int\limits_{\Omega }\,\frac{\varepsilon \,\partial
_{\xi ^{\alpha }}G(x,t)}{\left( \left\vert \nabla u_{\varepsilon }^{\alpha
}+th\triangle _{h,s}\nabla u_{\varepsilon }^{\alpha }\right\vert
^{2}+\varepsilon \right) ^{\frac{3}{2}}}\,\,\triangle _{h,s}\nabla
u_{\varepsilon }^{\alpha }\,\,$ \\
\end{tabular}%
\end{equation}%
for every $\alpha =1,...,m$. Adding for $\alpha =1,...,m$ we get%
\begin{equation}
A=-B
\end{equation}%
where%
\begin{equation*}
A=\sum\limits_{\alpha =1}^{m}A_{\alpha }
\end{equation*}%
and%
\begin{equation*}
B=\sum\limits_{\alpha =1}^{m}B_{\alpha }
\end{equation*}

Let us consider%
\begin{equation*}
\begin{tabular}{l}
$A$ \\ 
$=\left( p-1\right) \int\limits_{\Omega }\varepsilon \,\,\sum\limits_{\alpha
=1}^{m}I_{1,\alpha }\left( \triangle _{h,s}\left( \nabla u_{\varepsilon
}^{\alpha }\right) \right) ^{2}\eta ^{2}dx$ \\ 
$+\int\limits_{0}^{1}\int\limits_{\Omega }\sum\limits_{\alpha ,\beta
=1}^{m}\partial _{s^{\beta }s^{\alpha }}G\left( x,t\right) \frac{%
u_{\varepsilon }^{\beta }+th\triangle _{h,s}u_{\varepsilon }^{\beta }}{\sqrt{%
\left\vert u_{\varepsilon }^{\beta }+th\triangle _{h,s}u_{\varepsilon
}^{\beta }\right\vert ^{2}+\varepsilon }}\frac{u_{\varepsilon }^{\alpha
}+th\triangle _{h,s}u_{\varepsilon }^{\alpha }}{\sqrt{\left\vert
u_{\varepsilon }^{\alpha }+th\triangle _{h,s}u_{\varepsilon }^{\alpha
}\right\vert ^{2}+\varepsilon }}\,\,$ \\ 
$\cdot \triangle _{h,s}u_{\varepsilon }^{\beta }\,\triangle
_{h,s}u_{\varepsilon }^{\alpha }\,\eta ^{2}\,dx\,dt$ \\ 
$+\int\limits_{0}^{1}\int\limits_{\Omega }\sum\limits_{\alpha ,\beta
=1}^{m}\partial _{\xi ^{\beta }s^{\alpha }}G\left( x,t\right) \frac{\nabla
u_{\varepsilon }^{\beta }+th\triangle _{h,s}\nabla u_{\varepsilon }^{\beta }%
}{\sqrt{\left\vert \nabla u_{\varepsilon }^{\beta }+th\triangle _{h,s}\nabla
u_{\varepsilon }^{\beta }\right\vert ^{2}+\varepsilon }}\frac{u_{\varepsilon
}^{\alpha }+th\triangle _{h,s}u_{\varepsilon }^{\alpha }}{\sqrt{\left\vert
u_{\varepsilon }^{\alpha }+th\triangle _{h,s}u_{\varepsilon }^{\alpha
}\right\vert ^{2}+\varepsilon }}\,\,$ \\ 
$\cdot \triangle _{h,s}\nabla u_{\varepsilon }^{\beta }\,\triangle
_{h,s}u_{\varepsilon }^{\alpha }\,\eta ^{2}\,dx\,dt$ \\
\end{tabular}%
\end{equation*}%
since, $I_{1,\alpha }=\int\limits_{0}^{1}\left\vert \nabla u_{\varepsilon
}^{\alpha }\left( x\right) +th\triangle _{h,s}\nabla u_{\varepsilon
}^{\alpha }\left( x\right) \right\vert ^{p-2}\,\,dt\,\geq 0$, then%
\begin{equation*}
\left( p-1\right) \int\limits_{\Omega }\varepsilon \,\,\sum\limits_{\alpha
=1}^{m}I_{1,\alpha }\left( \triangle _{h,s}\left( \nabla u_{\varepsilon
}^{\alpha }\right) \right) ^{2}\eta ^{2}dx\geq 0
\end{equation*}%
Moreover using hypothesis H.2 we have $\partial _{s^{\alpha }}G(x,t)\geq 0$
and $\partial _{\xi ^{\alpha }}G(x,t)\geq 0$ for every $\alpha =1,...,m$ then%
\begin{equation*}
\int\limits_{0}^{1}\int\limits_{\Omega }\sum\limits_{\alpha =1}^{m}\frac{%
\varepsilon \,\partial _{s^{\alpha }}G(x,t)}{\left( \left\vert
u_{\varepsilon }^{\alpha }+th\triangle _{h,s}u_{\varepsilon }^{\alpha
}\right\vert ^{2}+\varepsilon \right) ^{\frac{3}{2}}}\,\,\left( \triangle
_{h,s}u_{\varepsilon }^{\alpha }\right) ^{2}\,\,\eta ^{2}\,dx\,dt\geq 0
\end{equation*}%
and%
\begin{equation*}
\int\limits_{0}^{1}\int\limits_{\Omega }\sum\limits_{\alpha =1}^{m}\frac{%
\varepsilon \,\partial _{\xi ^{\alpha }}G(x,t)}{\left( \left\vert \nabla
u_{\varepsilon }^{\alpha }+th\triangle _{h,s}\nabla u_{\varepsilon }^{\alpha
}\right\vert ^{2}+\varepsilon \right) ^{\frac{3}{2}}}\,\,\left( \triangle
_{h,s}\nabla u_{\varepsilon }^{\alpha }\right) ^{2}\,\eta ^{2}\,dx\,dt\,\geq
0
\end{equation*}%
Using these last observations we deduce that%
\begin{equation*}
\begin{tabular}{l}
$A\geq \left( p-1\right) \int\limits_{\Omega }\varepsilon
\,\,\sum\limits_{\alpha =1}^{m}I_{1,\alpha }\left( \triangle _{h,s}\left(
\nabla u_{\varepsilon }^{\alpha }\right) \right) ^{2}\eta ^{2}dx$ \\ 
$+\int\limits_{0}^{1}\int\limits_{\Omega }\sum\limits_{\alpha ,\beta
=1}^{m}\partial _{s^{\beta }s^{\alpha }}G\left( x,t\right) \frac{%
u_{\varepsilon }^{\beta }+th\triangle _{h,s}u_{\varepsilon }^{\beta }}{\sqrt{%
\left\vert u_{\varepsilon }^{\beta }+th\triangle _{h,s}u_{\varepsilon
}^{\beta }\right\vert ^{2}+\varepsilon }}\frac{u_{\varepsilon }^{\alpha
}+th\triangle _{h,s}u_{\varepsilon }^{\alpha }}{\sqrt{\left\vert
u_{\varepsilon }^{\alpha }+th\triangle _{h,s}u_{\varepsilon }^{\alpha
}\right\vert ^{2}+\varepsilon }}\,\,$ \\ 
$\cdot \triangle _{h,s}u_{\varepsilon }^{\beta }\,\triangle
_{h,s}u_{\varepsilon }^{\alpha }\,\eta ^{2}\,dx\,dt$ \\ 
$+\int\limits_{0}^{1}\int\limits_{\Omega }\sum\limits_{\alpha ,\beta
=1}^{m}\partial _{\xi ^{\beta }s^{\alpha }}G\left( x,t\right) \frac{\nabla
u_{\varepsilon }^{\beta }+th\triangle _{h,s}\nabla u_{\varepsilon }^{\beta }%
}{\sqrt{\left\vert \nabla u_{\varepsilon }^{\beta }+th\triangle _{h,s}\nabla
u_{\varepsilon }^{\beta }\right\vert ^{2}+\varepsilon }}\frac{u_{\varepsilon
}^{\alpha }+th\triangle _{h,s}u_{\varepsilon }^{\alpha }}{\sqrt{\left\vert
u_{\varepsilon }^{\alpha }+th\triangle _{h,s}u_{\varepsilon }^{\alpha
}\right\vert ^{2}+\varepsilon }}\,\,$ \\ 
$\cdot \triangle _{h,s}\nabla u_{\varepsilon }^{\beta }\,\triangle
_{h,s}u_{\varepsilon }^{\alpha }\,\eta ^{2}\,dx\,dt$ \\ 
$+\int\limits_{0}^{1}\int\limits_{\Omega }\sum\limits_{\alpha ,\beta
=1}^{m}\partial _{s^{\beta }\xi ^{\alpha }}G\left( x,t\right) \frac{%
u_{\varepsilon }^{\beta }+th\triangle _{h,s}u_{\varepsilon }^{\beta }}{\sqrt{%
\left\vert u_{\varepsilon }^{\beta }+th\triangle _{h,s}u_{\varepsilon
}^{\beta }\right\vert ^{2}+\varepsilon }}\frac{\nabla u_{\varepsilon
}^{\alpha }+th\triangle _{h,s}\nabla u_{\varepsilon }^{\alpha }}{\sqrt{%
\left\vert \nabla u_{\varepsilon }^{\alpha }+th\triangle _{h,s}\nabla
u_{\varepsilon }^{\alpha }\right\vert ^{2}+\varepsilon }}\,\,$ \\ 
$\cdot \triangle _{h,s}u_{\varepsilon }^{\beta }\,\,\triangle _{h,s}\nabla
u_{\varepsilon }^{\alpha }\,\eta ^{2}\,dx\,dt$ \\ 
$+\int\limits_{0}^{1}\int\limits_{\Omega }\sum\limits_{\alpha ,\beta
=1}^{m}\partial _{\xi ^{\beta }\xi ^{\alpha }}G\left( x,t\right) \frac{%
\nabla u_{\varepsilon }^{\beta }+th\triangle _{h,s}\nabla u_{\varepsilon
}^{\beta }}{\sqrt{\left\vert \nabla u_{\varepsilon }^{\beta }+th\triangle
_{h,s}\nabla u_{\varepsilon }^{\beta }\right\vert ^{2}+\varepsilon }}\frac{%
\nabla u^{\alpha }+th\triangle _{h,s}\nabla u_{\varepsilon }^{\alpha }}{%
\sqrt{\left\vert \nabla u_{\varepsilon }^{\alpha }+th\triangle _{h,s}\nabla
u_{\varepsilon }^{\alpha }\right\vert ^{2}+\varepsilon }}\,\,$ \\ 
$\cdot \triangle _{h,s}\nabla u_{\varepsilon }^{\beta }\,\triangle
_{h,s}\nabla u_{\varepsilon }^{\alpha }\,\eta ^{2}\,dx\,dt$%
\end{tabular}%
\end{equation*}%
Using hypothesis H.2 we get%
\begin{equation*}
A\geq \left( p-1\right) \int\limits_{\Omega }\varepsilon
\,\sum\limits_{\alpha =1}^{m}I_{1,\alpha }\left( \triangle _{h,s}\left(
\nabla u_{\varepsilon }^{\alpha }\right) \right) ^{2}\eta ^{2}dx+\vartheta
\int\limits_{\Omega }\Pi _{1}\left( x\right) \left( \left\vert \triangle
_{h,s}\left( \nabla u_{\varepsilon }\right) \right\vert ^{2}+\left\vert
\triangle _{h,s}u_{\varepsilon }\right\vert ^{2}\right) \eta ^{2}dx
\end{equation*}%
where 
\begin{equation*}
\Pi _{1}\left( x\right) =\int\limits_{0}^{1}\left( \left\vert u_{\varepsilon
}\left( x\right) +th\triangle _{h,s}u_{\varepsilon }\left( x\right)
\right\vert ^{2}+\varepsilon \right) ^{\frac{q-2}{2}}+\left( \left\vert
\nabla u_{\varepsilon }\left( x\right) +th\triangle _{h,s}\nabla
u_{\varepsilon }\left( x\right) \right\vert ^{2}+\varepsilon \right) ^{\frac{%
q-2}{2}}\,dt
\end{equation*}%
Moreover using Lemma 8 we get%
\begin{equation}
\Pi _{1}\left( x\right) \geq \nu \left( W_{1,\varepsilon
}^{q-2}+W_{2,\varepsilon }^{q-2}\right)
\end{equation}%
where 
\begin{equation}
W_{1,\varepsilon }=\left( \left\vert u_{\varepsilon }\left( x\right)
\right\vert ^{2}+\left\vert u_{\varepsilon }\left( x+he_{s}\right)
\right\vert ^{2}+\varepsilon \right) ^{\frac{1}{2}}
\end{equation}%
and%
\begin{equation}
W_{2,\varepsilon }=\left( \left\vert \nabla u_{\varepsilon }\left( x\right)
\right\vert ^{2}+\left\vert \nabla u_{\varepsilon }\left( x+he_{s}\right)
\right\vert ^{2}+\varepsilon \right) ^{\frac{q-2}{2}}\,
\end{equation}%
then it follows%
\begin{equation}
A\geq \left( p-1\right) \int\limits_{\Omega }\varepsilon
\,\sum\limits_{\alpha =1}^{m}I_{1,\alpha }\left( \triangle _{h,s}\left(
\nabla u_{\varepsilon }^{\alpha }\right) \right) ^{2}\eta ^{2}dx+\vartheta
\nu \int\limits_{\Omega }\left( W_{1,\varepsilon }^{q-2}+W_{2,\varepsilon
}^{q-2}\right) \left( \left\vert \triangle _{h,s}\left( \nabla
u_{\varepsilon }\right) \right\vert ^{2}+\left\vert \triangle
_{h,s}u_{\varepsilon }\right\vert ^{2}\right) \eta ^{2}dx
\end{equation}%
Recalling that 
\begin{equation*}
B\leq \left\vert B\right\vert \leq \sum\limits_{\alpha =1}^{m}\left\vert
B_{\alpha }\right\vert
\end{equation*}%
it follows%
\begin{equation}
\begin{tabular}{l}
$\left\vert B\right\vert $ \\ 
$\leq 2\left( p-1\right) \int\limits_{\Omega }\varepsilon
\sum\limits_{\alpha =1}^{m}\,\left\vert \triangle _{h,s}\left( \nabla
u_{\varepsilon }^{\alpha }\right) \right\vert \left\vert \triangle
_{h,s}u_{\varepsilon }^{\alpha }\right\vert I_{1,\alpha }\,\eta \,\left\vert
\nabla \eta \right\vert \,dx$ \\ 
$+2\int\limits_{0}^{1}\int\limits_{\Omega }\sum\limits_{\alpha ,\beta
=1}^{m}\left\vert \partial _{s^{\beta }\xi ^{\alpha }}G\left( x,t\right)
\right\vert \frac{\left\vert u_{\varepsilon }^{\beta }+th\triangle
_{h,s}u_{\varepsilon }^{\beta }\right\vert }{\sqrt{\left\vert u_{\varepsilon
}^{\beta }+th\triangle _{h,s}u_{\varepsilon }^{\beta }\right\vert
^{2}+\varepsilon }}\frac{\left\vert \nabla u_{\varepsilon }^{\alpha
}+th\triangle _{h,s}\nabla u_{\varepsilon }^{\alpha }\right\vert }{\sqrt{%
\left\vert \nabla u_{\varepsilon }^{\alpha }+th\triangle _{h,s}\nabla
u_{\varepsilon }^{\alpha }\right\vert ^{2}+\varepsilon }}\,$ \\ 
$\,\cdot \left\vert \triangle _{h,s}u_{\varepsilon }^{\beta }\right\vert
\,\,\left\vert \triangle _{h,s}u_{\varepsilon }^{\alpha }\right\vert \,\eta
\,\left\vert \nabla \eta \right\vert \,dx\,dt\,$ \\ 
$+2\int\limits_{0}^{1}\int\limits_{\Omega }\sum\limits_{\alpha ,\beta
=1}^{m}\left\vert \partial _{\xi ^{\beta }\xi ^{\alpha }}G\left( x,t\right)
\right\vert \frac{\left\vert \nabla u_{\varepsilon }^{\beta }+th\triangle
_{h,s}\nabla u_{\varepsilon }^{\beta }\right\vert }{\sqrt{\left\vert \nabla
u_{\varepsilon }^{\beta }+th\triangle _{h,s}\nabla u_{\varepsilon }^{\beta
}\right\vert ^{2}+\varepsilon }}\frac{\left\vert \nabla u_{\varepsilon
}^{\alpha }+th\triangle _{h,s}\nabla u_{\varepsilon }^{\alpha }\right\vert }{%
\sqrt{\left\vert \nabla u_{\varepsilon }^{\alpha }+th\triangle _{h,s}\nabla
u_{\varepsilon }^{\alpha }\right\vert ^{2}+\varepsilon }}\,\,$ \\ 
$\cdot \left\vert \triangle _{h,s}\nabla u^{\beta }\right\vert
\,\,\left\vert \triangle _{h,s}u^{\alpha }\right\vert \,\eta \,\left\vert
\nabla \eta \right\vert \,dx\,dt$ \\ 
$+2\int\limits_{0}^{1}\int\limits_{\Omega }\sum\limits_{\alpha =1}^{m}\frac{%
\varepsilon \,\partial _{\xi ^{\alpha }}G(x,t)}{\left( \left\vert \nabla
u_{\varepsilon }^{\alpha }+th\triangle _{h,s}\nabla u_{\varepsilon }^{\alpha
}\right\vert ^{2}+\varepsilon \right) ^{\frac{3}{2}}}\,\,\left\vert
\triangle _{h,s}\nabla u_{\varepsilon }^{\alpha }\right\vert \,\left\vert
\triangle _{h,s}u_{\varepsilon }^{\alpha }\right\vert \,\eta \,\left\vert
\nabla \eta \right\vert \,dx\,dt$%
\end{tabular}%
\end{equation}%
Now we have to estimate the various terms of (4.21), we start by estimating
the first term%
\begin{equation*}
2\left( p-1\right) \int\limits_{\Omega }\varepsilon \,\sum\limits_{\alpha
=1}^{m}\left\vert \triangle _{h,s}\left( \nabla u_{\varepsilon }^{\alpha
}\right) \right\vert \left\vert \triangle _{h,s}u_{\varepsilon }^{\alpha
}\right\vert I_{1,\alpha }\,\eta \,\left\vert \nabla \eta \right\vert \,dx
\end{equation*}%
using the $\delta $-Young inequality we obtain%
\begin{equation}
\begin{tabular}{l}
$2\left( p-1\right) \int\limits_{\Omega }\varepsilon \,\sum\limits_{\alpha
=1}^{m}\left\vert \triangle _{h,s}\left( \nabla u_{\varepsilon }^{\alpha
}\right) \right\vert \left\vert \triangle _{h,s}u_{\varepsilon }^{\alpha
}\right\vert I_{1,\alpha }\,\eta \,\left\vert \nabla \eta \right\vert \,dx$
\\ 
$\leq 2\left( p-1\right) \int\limits_{\Omega }\varepsilon
\sum\limits_{\alpha =1}^{m}I_{1,\alpha }\,\left[ \frac{\delta }{2}%
\,\left\vert \triangle _{h,s}\left( \nabla u_{\varepsilon }^{\alpha }\right)
\right\vert ^{2}\eta ^{2}+\frac{1}{2\delta }\left\vert \triangle
_{h,s}u_{\varepsilon }^{\alpha }\right\vert ^{2}\,\left\vert \nabla \eta
\right\vert ^{2}\right] \,dx$ \\ 
$\leq \delta \left( p-1\right) \int\limits_{\Omega }\varepsilon
\sum\limits_{\alpha =1}^{m}I_{1,\alpha }\,\,\left\vert \triangle
_{h,s}\left( \nabla u_{\varepsilon }^{\alpha }\right) \right\vert ^{2}\eta
^{2}\,dx$ \\ 
$+\frac{\left( p-1\right) }{\delta }\int\limits_{\Omega }\varepsilon
\sum\limits_{\alpha =1}^{m}I_{1,\alpha }\,\left\vert \triangle
_{h,s}u_{\varepsilon }^{\alpha }\right\vert ^{2}\,\left\vert \nabla \eta
\right\vert ^{2}\,dx$%
\end{tabular}%
\end{equation}%
where $\delta >0$\ is a real number that we fix later. Now let us consider
the second term%
\begin{eqnarray*}
&&p\int\limits_{0}^{1}\int\limits_{\Omega }\sum\limits_{\alpha ,\beta
=1}^{m}\left\vert \partial _{s^{\beta }\xi ^{\alpha }}G\left( x,t\right)
\right\vert \frac{\left\vert u_{\varepsilon }^{\beta }+th\triangle
_{h,s}u_{\varepsilon }^{\beta }\right\vert }{\sqrt{\left\vert u_{\varepsilon
}^{\beta }+th\triangle _{h,s}u_{\varepsilon }^{\beta }\right\vert
^{2}+\varepsilon }}\frac{\left\vert \nabla u_{\varepsilon }^{\alpha
}+th\triangle _{h,s}\nabla u_{\varepsilon }^{\alpha }\right\vert }{\sqrt{%
\left\vert \nabla u_{\varepsilon }^{\alpha }+th\triangle _{h,s}\nabla
u_{\varepsilon }^{\alpha }\right\vert ^{2}+\varepsilon }}\,\, \\
&&\cdot \left\vert \triangle _{h,s}u_{\varepsilon }^{\beta }\right\vert
\,\,\left\vert \triangle _{h,s}u_{\varepsilon }^{\alpha }\right\vert \,\eta
\,\left\vert \nabla \eta \right\vert \,dx\,dt
\end{eqnarray*}%
then using hypothesis H.2 we get%
\begin{equation*}
\begin{tabular}{l}
$\Theta _{1}=p\int\limits_{0}^{1}\int\limits_{\Omega }\sum\limits_{\alpha
,\beta =1}^{m}\left\vert \partial _{s^{\beta }\xi ^{\alpha }}G\left(
x,t\right) \right\vert \frac{\left\vert u_{\varepsilon }^{\beta
}+th\triangle _{h,s}u_{\varepsilon }^{\beta }\right\vert }{\sqrt{\left\vert
u_{\varepsilon }^{\beta }+th\triangle _{h,s}u_{\varepsilon }^{\beta
}\right\vert ^{2}+\varepsilon }}$ \\ 
$\cdot \frac{\left\vert \nabla u_{\varepsilon }^{\alpha }+th\triangle
_{h,s}\nabla u_{\varepsilon }{}^{\alpha }\right\vert }{\sqrt{\left\vert
\nabla u_{\varepsilon }^{\alpha }+th\triangle _{h,s}\nabla u_{\varepsilon
}^{\alpha }\right\vert ^{2}+\varepsilon }}\left\vert \triangle
_{h,s}u_{\varepsilon }^{\beta }\right\vert \,\,\left\vert \triangle
_{h,s}u_{\varepsilon }^{\alpha }\right\vert \,\eta \,\left\vert \nabla \eta
\right\vert \,dx\,dt\,$ \\ 
$\leq pLm^{2}\int\limits_{0}^{1}\int\limits_{\Omega }\left[ \left(
\left\vert u_{\varepsilon }+th\triangle _{h,s}u_{\varepsilon }\right\vert
^{2}+\varepsilon \right) ^{\frac{q-2}{2}}+\left( \left\vert \nabla
u_{\varepsilon }+th\triangle _{h,s}\nabla u_{\varepsilon }\right\vert
^{2}+\varepsilon \right) ^{\frac{q-2}{2}}\right] \,$ \\ 
$\left\vert \triangle _{h,s}u_{\varepsilon }\right\vert ^{2}\,\,\eta
\,\left\vert \nabla \eta \right\vert \,dx\,dt$%
\end{tabular}%
\end{equation*}%
Since by Lemma 8%
\begin{equation}
\int\limits_{0}^{1}\left[ \left( \left\vert u_{\varepsilon }+th\triangle
_{h,s}u_{\varepsilon }\right\vert ^{2}+\varepsilon \right) ^{\frac{q-2}{2}%
}+\left( \left\vert \nabla u_{\varepsilon }+th\triangle _{h,s}\nabla
u_{\varepsilon }\right\vert ^{2}+\varepsilon \right) ^{\frac{q-2}{2}}\right]
\,\,dt\leq W_{1,\varepsilon }^{q-2}+W_{2,\varepsilon }^{q-2}
\end{equation}%
it follows%
\begin{equation}
\Theta _{1}\leq pLm^{2}\int\limits_{\Omega }\left[ W_{1,\varepsilon
}^{q-2}+W_{2,\varepsilon }^{q-2}\right] \left\vert \triangle
_{h,s}u_{\varepsilon }\right\vert ^{2}\,\,\eta \,\left\vert \nabla \eta
\right\vert \,dx
\end{equation}%
Now we consider the term%
\begin{eqnarray*}
&&p\int\limits_{0}^{1}\int\limits_{\Omega }\sum\limits_{\alpha ,\beta
=1}^{m}\left\vert \partial _{\xi ^{\beta }\xi ^{\alpha }}G\left( x,t\right)
\right\vert \frac{\left\vert \nabla u_{\varepsilon }^{\beta }+th\triangle
_{h,s}\nabla u_{\varepsilon }^{\beta }\right\vert }{\sqrt{\left\vert \nabla
u_{\varepsilon }^{\beta }+th\triangle _{h,s}\nabla u_{\varepsilon }^{\beta
}\right\vert ^{2}+\varepsilon }}\frac{\left\vert \nabla u_{\varepsilon
}^{\alpha }+th\triangle _{h,s}\nabla u_{\varepsilon }^{\alpha }\right\vert }{%
\sqrt{\left\vert \nabla u_{\varepsilon }^{\alpha }+th\triangle _{h,s}\nabla
u_{\varepsilon }^{\alpha }\right\vert ^{2}+\varepsilon }}\,\, \\
&&\left\vert \triangle _{h,s}\nabla u_{\varepsilon }^{\beta }\right\vert
\,\,\left\vert \triangle _{h,s}u_{\varepsilon }^{\alpha }\right\vert \,\eta
\,\left\vert \nabla \eta \right\vert \,dx\,dt
\end{eqnarray*}%
then using hypothesis H.2 and $\tilde{\delta}$-Young inequality we obtain%
\begin{equation*}
\begin{tabular}{l}
$\Theta _{2}=p\int\limits_{0}^{1}\int\limits_{\Omega }\sum\limits_{\alpha
,\beta =1}^{m}\left\vert \partial _{\xi ^{\beta }\xi ^{\alpha }}G\left(
x,t\right) \right\vert \frac{\left\vert \nabla u_{\varepsilon }^{\beta
}+th\triangle _{h,s}\nabla u_{\varepsilon }^{\beta }\right\vert }{\sqrt{%
\left\vert \nabla u_{\varepsilon }^{\beta }+th\triangle _{h,s}\nabla
u_{\varepsilon }^{\beta }\right\vert ^{2}+\varepsilon }}\frac{\left\vert
\nabla u_{\varepsilon }^{\alpha }+th\triangle _{h,s}\nabla u_{\varepsilon
}^{\alpha }\right\vert }{\sqrt{\left\vert \nabla u_{\varepsilon }^{\alpha
}+th\triangle _{h,s}\nabla u_{\varepsilon }^{\alpha }\right\vert
^{2}+\varepsilon }}\,\,$ \\ 
$\left\vert \triangle _{h,s}\nabla u_{\varepsilon }^{\beta }\right\vert
\,\,\left\vert \triangle _{h,s}u_{\varepsilon }^{\alpha }\right\vert \,\eta
\,\left\vert \nabla \eta \right\vert \,dx\,dt$ \\ 
$\leq pL\int\limits_{0}^{1}\int\limits_{\Omega }\sum\limits_{\alpha ,\beta
=1}^{m}\left[ \left( \left\vert u+th\triangle _{h,s}u_{\varepsilon
}\right\vert ^{2}+\varepsilon \right) ^{\frac{q-2}{2}}+\left( \left\vert
\nabla u+th\triangle _{h,s}\nabla u_{\varepsilon }\right\vert
^{2}+\varepsilon \right) ^{\frac{q-2}{2}}\right] \,\,\,$ \\ 
$\left\vert \triangle _{h,s}\nabla u_{\varepsilon }^{\beta }\right\vert
\,\,\left\vert \triangle _{h,s}u_{\varepsilon }^{\alpha }\right\vert \,\eta
\,\left\vert \nabla \eta \right\vert \,dx\,dt$ \\ 
$\leq \frac{pLm^{2}}{2}\int\limits_{0}^{1}\int\limits_{\Omega }\left[ \left(
\left\vert u+th\triangle _{h,s}u_{\varepsilon }\right\vert ^{2}+\varepsilon
\right) ^{\frac{q-2}{2}}+\left( \left\vert \nabla u+th\triangle _{h,s}\nabla
u_{\varepsilon }\right\vert ^{2}+\varepsilon \right) ^{\frac{q-2}{2}}\right]
\,\,\,$ \\ 
$\left[ \tilde{\delta}\left\vert \triangle _{h,s}\nabla u_{\varepsilon
}\right\vert ^{2}\,\eta ^{2}+\frac{\,\left\vert \triangle
_{h,s}u_{\varepsilon }\right\vert ^{2}\left\vert \nabla \eta \right\vert ^{2}%
}{\tilde{\delta}}\right] \,\,\,dx\,dt$%
\end{tabular}%
\end{equation*}%
where $\tilde{\delta}>0$\ is a real positive number that we fix later. Using
Lemma 8 it follows%
\begin{equation}
\begin{tabular}{l}
$\Theta _{2}\leq \tilde{\delta}\frac{pLm^{2}}{2}\int\limits_{\Omega }\left[
W_{1,\varepsilon }^{q-2}+W_{2,\varepsilon }^{q-2}\right] \left\vert
\triangle _{h,s}\nabla u_{\varepsilon }\right\vert ^{2}\,\eta ^{2}\,dx\,\,$
\\ 
$+\frac{pLm^{2}}{2\tilde{\delta}}\int\limits_{\Omega }\left[
W_{1,\varepsilon }^{q-2}+W_{2,\varepsilon }^{q-2}\right] \left\vert
\triangle _{h,s}u_{\varepsilon }\right\vert ^{2}\left\vert \nabla \eta
\right\vert ^{2}\,dx$%
\end{tabular}%
\end{equation}

Finally we consider the term%
\begin{equation}
\left\vert \int\limits_{0}^{1}\int\limits_{\Omega }\sum\limits_{\alpha
=1}^{m}\frac{2\varepsilon \,\partial _{\xi ^{\alpha }}G(x,t)}{\left(
\left\vert \nabla u_{\varepsilon }^{\alpha }+th\triangle _{h,s}\nabla
u_{\varepsilon }^{\alpha }\right\vert ^{2}+\varepsilon \right) ^{\frac{3}{2}}%
}\,\,\triangle _{h,s}\nabla u_{\varepsilon }^{\alpha }\,\triangle
_{h,s}u_{\varepsilon }^{\alpha }\,\eta \,\nabla \eta \,dx\,dt\right\vert 
\end{equation}%
using the Fubini's\ Theorem it follows%
\begin{equation}
\begin{tabular}{l}
$\Theta _{3}$ \\ 
$=\left\vert \int\limits_{0}^{1}\int\limits_{\Omega }\sum\limits_{\alpha
=1}^{m}\frac{2\varepsilon \,\partial _{\xi ^{\alpha }}G(x,t)}{\left(
\left\vert \nabla u_{\varepsilon }^{\alpha }+th\triangle _{h,s}\nabla
u_{\varepsilon }^{\alpha }\right\vert ^{2}+\varepsilon \right) ^{\frac{3}{2}}%
}\,\,\triangle _{h,s}\nabla u_{\varepsilon }^{\alpha }\,\triangle
_{h,s}u_{\varepsilon }^{\alpha }\,\eta \,\nabla \eta \,dx\,dt\right\vert $
\\ 
$=2\varepsilon \left\vert \int\limits_{\Omega }\sum\limits_{\alpha
=1}^{m}\triangle _{h,s}u_{\varepsilon }^{\alpha }\,\eta \,\nabla \eta
\,\,\int\limits_{0}^{1}\,\partial _{\xi ^{\alpha }}G(x,t)\,\,\frac{\triangle
_{h,s}\nabla u_{\varepsilon }^{\alpha }}{\left( \left\vert \nabla
u_{\varepsilon }^{\alpha }+th\triangle _{h,s}\nabla u_{\varepsilon }^{\alpha
}\right\vert ^{2}+\varepsilon \right) ^{\frac{3}{2}}}\,dt\,dx\right\vert $
\\ 
$=2\varepsilon \left\vert \int\limits_{\Omega }\sum\limits_{\alpha
=1}^{m}\triangle _{h,s}u_{\varepsilon }^{\alpha }\,\eta \,\nabla \eta
\,\,\int\limits_{0}^{1}\,\partial _{\xi ^{\alpha }}G(x,t)\,\frac{\nabla
u_{\varepsilon }^{\alpha }+th\triangle _{h,s}\nabla u_{\varepsilon }^{\alpha
}-th\triangle _{h,s}\nabla u_{\varepsilon }^{\alpha }}{\left( \left\vert
\nabla u_{\varepsilon }^{\alpha }+th\triangle _{h,s}\nabla u_{\varepsilon
}^{\alpha }\right\vert ^{2}+\varepsilon \right) ^{\frac{3}{2}}}%
\,dt\,dx\right\vert $%
\end{tabular}%
\end{equation}%
then%
\begin{equation}
\begin{tabular}{l}
$\Theta _{3}$ \\ 
$\leq 2\varepsilon \sum\limits_{\alpha =1}^{m}\left\vert \int\limits_{\Omega
}\triangle _{h,s}u_{\varepsilon }^{\alpha }\,\eta \,\nabla \eta
\,\,\int\limits_{0}^{1}\,\partial _{\xi ^{\alpha }}G(x,t)\,\frac{\nabla
u_{\varepsilon }^{\alpha }+th\triangle _{h,s}\nabla u_{\varepsilon }^{\alpha
}}{\left( \left\vert \nabla u_{\varepsilon }^{\alpha }+th\triangle
_{h,s}\nabla u_{\varepsilon }^{\alpha }\right\vert ^{2}+\varepsilon \right)
^{\frac{3}{2}}}\,dt\,dx\right\vert $ \\ 
$+2\varepsilon \sum\limits_{\alpha =1}^{m}\left\vert \int\limits_{\Omega
}\triangle _{h,s}u_{\varepsilon }^{\alpha }\,\eta \,\nabla \eta
\,\,\int\limits_{0}^{1}\,\partial _{\xi ^{\alpha }}G(x,t)\,\frac{th\triangle
_{h,s}\nabla u_{\varepsilon }^{\alpha }}{\left( \left\vert \nabla
u_{\varepsilon }^{\alpha }+th\triangle _{h,s}\nabla u_{\varepsilon }^{\alpha
}\right\vert ^{2}+\varepsilon \right) ^{\frac{3}{2}}}\,dt\,dx\right\vert $%
\end{tabular}%
\end{equation}%
Now we have to estimate the two integrals of (4.28), so let us start by
considering the following integral%
\begin{equation}
\varepsilon \left\vert \int\limits_{\Omega }\triangle _{h,s}u_{\varepsilon
}^{\alpha }\,\eta \,\nabla \eta \,\,\int\limits_{0}^{1}\,\partial _{\xi
^{\alpha }}G(x,t)\,\frac{\nabla u_{\varepsilon }^{\alpha }+th\triangle
_{h,s}\nabla u_{\varepsilon }^{\alpha }}{\left( \left\vert \nabla
u_{\varepsilon }^{\alpha }+th\triangle _{h,s}\nabla u_{\varepsilon }^{\alpha
}\right\vert ^{2}+\varepsilon \right) ^{\frac{3}{2}}}\,dt\,dx\right\vert 
\end{equation}%
then we get%
\begin{equation}
\begin{tabular}{l}
$2\varepsilon \left\vert \int\limits_{\Omega }\triangle _{h,s}u_{\varepsilon
}^{\alpha }\,\eta \,\nabla \eta \,\,\int\limits_{0}^{1}\,\partial _{\xi
^{\alpha }}G(x,t)\,\frac{\nabla u_{\varepsilon }^{\alpha }+th\triangle
_{h,s}\nabla u_{\varepsilon }^{\alpha }}{\left( \left\vert \nabla
u_{\varepsilon }^{\alpha }+th\triangle _{h,s}\nabla u_{\varepsilon }^{\alpha
}\right\vert ^{2}+\varepsilon \right) ^{\frac{3}{2}}}\,dt\,dx\right\vert $
\\ 
$\leq 2\varepsilon \int\limits_{\Omega }\left\vert \triangle
_{h,s}u_{\varepsilon }^{\alpha }\right\vert \,\eta \,\left\vert \nabla \eta
\right\vert \,\,\int\limits_{0}^{1}\,\left\vert \partial _{\xi ^{\alpha
}}G(x,t)\right\vert \,\frac{\left\vert \nabla u_{\varepsilon }^{\alpha
}+th\triangle _{h,s}\nabla u_{\varepsilon }^{\alpha }\right\vert }{\left(
\left\vert \nabla u_{\varepsilon }^{\alpha }+th\triangle _{h,s}\nabla
u_{\varepsilon }^{\alpha }\right\vert ^{2}+\varepsilon \right) ^{\frac{3}{2}}%
}\,dt\,dx$ \\ 
$\leq 2\varepsilon \int\limits_{\Omega }\left\vert \triangle
_{h,s}u_{\varepsilon }^{\alpha }\right\vert \,\eta \,\left\vert \nabla \eta
\right\vert \,\,\int\limits_{0}^{1}\,\,\frac{\left\vert \partial _{\xi
^{\alpha }}G(x,t)\right\vert \,}{\left( \left\vert \nabla u_{\varepsilon
}^{\alpha }+th\triangle _{h,s}\nabla u_{\varepsilon }^{\alpha }\right\vert
^{2}+\varepsilon \right) ^{\frac{1}{2}}}\,dt\,dx$ \\ 
$\leq 2\sqrt{\varepsilon }\int\limits_{\Omega }\left\vert \triangle
_{h,s}u_{\varepsilon }^{\alpha }\right\vert \,\eta \,\left\vert \nabla \eta
\right\vert \,\,\int\limits_{0}^{1}\,\,\left\vert \partial _{\xi ^{\alpha
}}G(x,t)\right\vert \,\,dt\,dx$ \\ 
$\leq 2\sqrt{\varepsilon }Lm\int\limits_{\Omega }\left\vert \triangle
_{h,s}u_{\varepsilon }^{\alpha }\right\vert \,\eta \,\left\vert \nabla \eta
\right\vert \,\,\int\limits_{0}^{1}\left( \left\vert \nabla u_{\varepsilon
}+th\triangle _{h,s}\nabla u_{\varepsilon }\right\vert ^{2}+\varepsilon
\right) ^{\frac{q-1}{2}}+\left( \left\vert u_{\varepsilon }+th\triangle
_{h,s}u_{\varepsilon }\right\vert ^{2}+\varepsilon \right) ^{\frac{q-1}{2}%
}\,\,dt\,dx$%
\end{tabular}%
\end{equation}%
Since%
\begin{equation}
\begin{tabular}{l}
$\int\limits_{0}^{1}\left( \left\vert \nabla u_{\varepsilon }+th\triangle
_{h,s}\nabla u_{\varepsilon }\right\vert ^{2}+\varepsilon \right) ^{\frac{q-1%
}{2}}+\left( \left\vert u_{\varepsilon }+th\triangle _{h,s}u_{\varepsilon
}\right\vert ^{2}+\varepsilon \right) ^{\frac{q-1}{2}}\,\,dt\,$ \\ 
$=\int\limits_{0}^{1}\left( \left\vert \nabla u_{\varepsilon }+th\triangle
_{h,s}\nabla u_{\varepsilon }\right\vert ^{2}+\varepsilon \right) ^{\frac{q-2%
}{2}}\left( \left\vert \nabla u_{\varepsilon }+th\triangle _{h,s}\nabla
u_{\varepsilon }\right\vert ^{2}+\varepsilon \right) ^{\frac{1}{2}}\,\,dt\,$
\\ 
$+\int\limits_{0}^{1}\left( \left\vert u_{\varepsilon }+th\triangle
_{h,s}u_{\varepsilon }\right\vert ^{2}+\varepsilon \right) ^{\frac{q-1}{2}%
}\left( \left\vert u_{\varepsilon }+th\triangle _{h,s}u_{\varepsilon
}\right\vert ^{2}+\varepsilon \right) ^{\frac{1}{2}}\,\,dt\,$ \\ 
$\leq \int\limits_{0}^{1}\left( \left\vert \nabla u_{\varepsilon
}+th\triangle _{h,s}\nabla u_{\varepsilon }\right\vert ^{2}+\varepsilon
\right) ^{\frac{q-2}{2}}\left( \left\vert \nabla u_{\varepsilon
}+th\triangle _{h,s}\nabla u_{\varepsilon }\right\vert +1\right) \,\,dt$ \\ 
$+\int\limits_{0}^{1}\left( \left\vert u_{\varepsilon }+th\triangle
_{h,s}u_{\varepsilon }\right\vert ^{2}+\varepsilon \right) ^{\frac{q-2}{2}%
}\left( \left\vert u_{\varepsilon }+th\triangle _{h,s}u_{\varepsilon
}\right\vert +1\right) \,\,dt$ \\ 
$\leq \int\limits_{0}^{1}\left( \left\vert \nabla u_{\varepsilon
}+th\triangle _{h,s}\nabla u_{\varepsilon }\right\vert ^{2}+\varepsilon
\right) ^{\frac{q-2}{2}}\left( \left\vert \nabla u_{\varepsilon }\right\vert
+\left\vert \triangle _{h,s}\nabla u_{\varepsilon }\right\vert +1\right)
\,\,dt$ \\ 
$+\int\limits_{0}^{1}\left( \left\vert u_{\varepsilon }+th\triangle
_{h,s}u_{\varepsilon }\right\vert ^{2}+\varepsilon \right) ^{\frac{q-2}{2}%
}\left( \left\vert u_{\varepsilon }\right\vert +\left\vert \triangle
_{h,s}u_{\varepsilon }\right\vert +1\right) \,\,dt$ \\ 
$\leq \left( \left\vert \nabla u_{\varepsilon }\right\vert +\left\vert
\triangle _{h,s}\nabla u_{\varepsilon }\right\vert +1\right)
\int\limits_{0}^{1}\left( \left\vert \nabla u_{\varepsilon }+th\triangle
_{h,s}\nabla u_{\varepsilon }\right\vert ^{2}+\varepsilon \right) ^{\frac{q-2%
}{2}}\,\,dt$ \\ 
$+\left( \left\vert u_{\varepsilon }\right\vert +\left\vert \triangle
_{h,s}u_{\varepsilon }\right\vert +1\right) \int\limits_{0}^{1}\left(
\left\vert u_{\varepsilon }+th\triangle _{h,s}u_{\varepsilon }\right\vert
^{2}+\varepsilon \right) ^{\frac{q-2}{2}}\,\,dt$%
\end{tabular}%
\end{equation}%
then, using Lemma 8, it follows%
\begin{equation}
\begin{tabular}{l}
$2\varepsilon \sum\limits_{\alpha =1}^{m}\left\vert \int\limits_{\Omega
}\triangle _{h,s}u_{\varepsilon }^{\alpha }\,\eta \,\nabla \eta
\,\,\int\limits_{0}^{1}\,\partial _{\xi ^{\alpha }}G(x,t)\,\frac{\nabla
u_{\varepsilon }^{\alpha }+th\triangle _{h,s}\nabla u_{\varepsilon }^{\alpha
}}{\left( \left\vert \nabla u_{\varepsilon }^{\alpha }+th\triangle
_{h,s}\nabla u_{\varepsilon }^{\alpha }\right\vert ^{2}+\varepsilon \right)
^{\frac{3}{2}}}\,dt\,dx\right\vert $ \\ 
$\leq 2\sqrt{\varepsilon }Lm^{2}\int\limits_{\Omega }\left\vert \triangle
_{h,s}u_{\varepsilon }\right\vert \,\eta \,\left\vert \nabla \eta
\right\vert \left[ \left( \left\vert \nabla u_{\varepsilon }\right\vert
+\left\vert \triangle _{h,s}\nabla u_{\varepsilon }\right\vert +1\right)
\,\,W_{2,\varepsilon }^{q-2}+\left( \left\vert u_{\varepsilon }\right\vert
+\left\vert \triangle _{h,s}u_{\varepsilon }\right\vert +1\right)
W_{1,\varepsilon }^{q-2}\right] \,dx$ \\ 
$\leq 2Lm^{2}\int\limits_{\Omega }\left\vert \triangle _{h,s}u_{\varepsilon
}\right\vert \,\eta \,\left\vert \nabla \eta \right\vert \,\,\left(
\left\vert u_{\varepsilon }\right\vert +\left\vert \triangle
_{h,s}u_{\varepsilon }\right\vert +1\right) W_{1,\varepsilon }^{q-2}\,dx$ \\ 
$+2Lm^{2}\int\limits_{\Omega }\left\vert \triangle _{h,s}u_{\varepsilon
}\right\vert \,\eta \,\left\vert \nabla \eta \right\vert \,\,\left(
\left\vert \nabla u_{\varepsilon }\right\vert +\left\vert \triangle
_{h,s}\nabla u_{\varepsilon }\right\vert +1\right) \,\,W_{2,\varepsilon
}^{q-2}\,dx$%
\end{tabular}%
\end{equation}%
Let us consider the following integral%
\begin{equation}
2Lm^{2}\int\limits_{\Omega }\left\vert \triangle _{h,s}u_{\varepsilon
}\right\vert \,\eta \,\left\vert \nabla \eta \right\vert \,\,\left(
\left\vert u_{\varepsilon }\right\vert +\left\vert \triangle
_{h,s}u_{\varepsilon }\right\vert +1\right) W_{1,\varepsilon }^{q-2}\,dx
\end{equation}%
using the Young Inequality we get%
\begin{equation}
\begin{tabular}{l}
$2Lm^{2}\int\limits_{\Omega }\left\vert \triangle _{h,s}u_{\varepsilon
}\right\vert \,\eta \,\left\vert \nabla \eta \right\vert \,\,\left(
\left\vert u_{\varepsilon }\right\vert +\left\vert \triangle
_{h,s}u_{\varepsilon }\right\vert +1\right) W_{1,\varepsilon }^{q-2}\,dx$ \\ 
$\leq 2Lm^{2}\int\limits_{\Omega }\left\vert \triangle _{h,s}u_{\varepsilon
}\right\vert \,\eta \,\left\vert \nabla \eta \right\vert \,\,\left\vert
u_{\varepsilon }\right\vert W_{1,\varepsilon }^{q-2}\,dx$ \\ 
$+2Lm^{2}\int\limits_{\Omega }\left\vert \triangle _{h,s}u_{\varepsilon
}\right\vert ^{2}\,\eta \,\left\vert \nabla \eta \right\vert
\,\,W_{1,\varepsilon }^{q-2}\,dx$ \\ 
$+2Lm^{2}\int\limits_{\Omega }\left\vert \triangle _{h,s}u_{\varepsilon
}\right\vert \,\eta \,\left\vert \nabla \eta \right\vert
\,\,W_{1,\varepsilon }^{q-2}\,dx$ \\ 
$\leq Lm^{2}\int\limits_{\Omega }\left\vert \triangle _{h,s}u_{\varepsilon
}\right\vert ^{2}\,\,\left\vert \nabla \eta \right\vert
^{2}\,W_{1,\varepsilon }^{q-2}\,dx+Lm^{2}\int\limits_{\Omega }\eta
\,^{2}\,\left\vert u_{\varepsilon }\right\vert ^{2}W_{1,\varepsilon
}^{q-2}\,dx$ \\ 
$+2Lm^{2}\int\limits_{\Omega }\left\vert \triangle _{h,s}u_{\varepsilon
}\right\vert ^{2}\,\eta \,\left\vert \nabla \eta \right\vert
\,\,W_{1,\varepsilon }^{q-2}\,dx$ \\ 
$+Lm^{2}\int\limits_{\Omega }\left\vert \triangle _{h,s}u_{\varepsilon
}\right\vert ^{2}\,\,\left\vert \nabla \eta \right\vert
^{2}\,\,W_{1,\varepsilon }^{q-2}\,dx+Lm^{2}\int\limits_{\Omega }\,\eta
\,^{2}\,W_{1,\varepsilon }^{q-2}\,dx$%
\end{tabular}%
\end{equation}%
Let us consider the integral%
\begin{equation}
2Lm^{2}\int\limits_{\Omega }\left\vert \triangle _{h,s}u_{\varepsilon
}\right\vert \,\eta \,\left\vert \nabla \eta \right\vert \,\,\left(
\left\vert \nabla u_{\varepsilon }\right\vert +\left\vert \triangle
_{h,s}\nabla u_{\varepsilon }\right\vert +1\right) \,\,W_{2,\varepsilon
}^{q-2}\,dx
\end{equation}%
then, using the Young Inequality we obtain%
\begin{equation}
\begin{tabular}{l}
$2Lm^{2}\int\limits_{\Omega }\left\vert \triangle _{h,s}u_{\varepsilon
}\right\vert \,\eta \,\left\vert \nabla \eta \right\vert \,\,\left(
\left\vert \nabla u_{\varepsilon }\right\vert +\left\vert \triangle
_{h,s}\nabla u_{\varepsilon }\right\vert +1\right) \,\,W_{2,\varepsilon
}^{q-2}\,dx$ \\ 
$\leq 2Lm^{2}\int\limits_{\Omega }\left\vert \triangle _{h,s}u_{\varepsilon
}\right\vert \,\eta \,\left\vert \nabla \eta \right\vert \,\,\left\vert
\nabla u_{\varepsilon }\right\vert \,\,W_{2,\varepsilon }^{q-2}\,dx$ \\ 
$+2Lm^{2}\int\limits_{\Omega }\left\vert \triangle _{h,s}u_{\varepsilon
}\right\vert \,\eta \,\left\vert \nabla \eta \right\vert \,\,\left\vert
\triangle _{h,s}\nabla u_{\varepsilon }\right\vert \,\,W_{2,\varepsilon
}^{q-2}\,dx$ \\ 
$+2Lm^{2}\int\limits_{\Omega }\left\vert \triangle _{h,s}u_{\varepsilon
}\right\vert \,\eta \,\left\vert \nabla \eta \right\vert
\,\,\,W_{2,\varepsilon }^{q-2}\,dx$ \\ 
$\leq Lm^{2}\int\limits_{\Omega }\left\vert \triangle _{h,s}u_{\varepsilon
}\right\vert ^{2}\,\,\left\vert \nabla \eta \right\vert
^{2}\,\,\,\,W_{2,\varepsilon }^{q-2}\,dx$ \\ 
$+Lm^{2}\int\limits_{\Omega }\eta ^{2}\,\,\left\vert \nabla u_{\varepsilon
}\right\vert ^{2}\,\,W_{2,\varepsilon }^{q-2}\,dx$ \\ 
$+\tilde{\delta}_{2}Lm^{2}\int\limits_{\Omega }\,\eta ^{2}\,\,\,\left\vert
\triangle _{h,s}\nabla u_{\varepsilon }\right\vert ^{2}\,\,W_{2,\varepsilon
}^{q-2}\,dx$ \\ 
$+\frac{Lm^{2}}{\tilde{\delta}_{2}}\int\limits_{\Omega }\left\vert \triangle
_{h,s}u_{\varepsilon }\right\vert ^{2}\,\left\vert \nabla \eta \right\vert
^{2}\,\,W_{2,\varepsilon }^{q-2}\,dx$ \\ 
$+Lm^{2}\int\limits_{\Omega }\left\vert \triangle _{h,s}u_{\varepsilon
}\right\vert \,^{2}\,\left\vert \nabla \eta \right\vert
^{2}\,\,\,W_{2,\varepsilon }^{q-2}\,dx$ \\ 
$+Lm^{2}\int\limits_{\Omega }\,\eta \,^{2}\,\,W_{2,\varepsilon }^{q-2}\,dx$%
\end{tabular}%
\end{equation}%
Using (4.32), (4.34) and (4.36) we have%
\begin{equation}
\begin{tabular}{l}
$2\varepsilon \sum\limits_{\alpha =1}^{m}\left\vert \int\limits_{\Omega
}\triangle _{h,s}u_{\varepsilon }^{\alpha }\,\eta \,\nabla \eta
\,\,\int\limits_{0}^{1}\,\partial _{\xi ^{\alpha }}G(x,t)\,\frac{\nabla
u_{\varepsilon }^{\alpha }+th\triangle _{h,s}\nabla u_{\varepsilon }^{\alpha
}}{\left( \left\vert \nabla u_{\varepsilon }^{\alpha }+th\triangle
_{h,s}\nabla u_{\varepsilon }^{\alpha }\right\vert ^{2}+\varepsilon \right)
^{\frac{3}{2}}}\,dt\,dx\right\vert $ \\ 
$\leq Lm^{2}\int\limits_{\Omega }\left\vert \triangle _{h,s}u_{\varepsilon
}\right\vert ^{2}\,\,\left\vert \nabla \eta \right\vert
^{2}\,W_{1,\varepsilon }^{q-2}\,dx+Lm^{2}\int\limits_{\Omega }\eta
\,^{2}\,\left\vert u_{\varepsilon }\right\vert ^{2}W_{1,\varepsilon
}^{q-2}\,dx$ \\ 
$+2Lm^{2}\int\limits_{\Omega }\left\vert \triangle _{h,s}u_{\varepsilon
}\right\vert ^{2}\,\eta \,\left\vert \nabla \eta \right\vert
\,\,W_{1,\varepsilon }^{q-2}\,dx$ \\ 
$+Lm^{2}\int\limits_{\Omega }\left\vert \triangle _{h,s}u_{\varepsilon
}\right\vert ^{2}\,\,\left\vert \nabla \eta \right\vert
^{2}\,\,W_{1,\varepsilon }^{q-2}\,dx+Lm^{2}\int\limits_{\Omega }\,\eta
\,^{2}\,W_{1,\varepsilon }^{q-2}\,dx$ \\ 
$+Lm^{2}\int\limits_{\Omega }\left\vert \triangle _{h,s}u_{\varepsilon
}\right\vert ^{2}\,\,\left\vert \nabla \eta \right\vert
^{2}\,\,\,\,W_{2,\varepsilon }^{q-2}\,dx$ \\ 
$+Lm^{2}\int\limits_{\Omega }\eta ^{2}\,\,\left\vert \nabla u_{\varepsilon
}\right\vert ^{2}\,\,W_{2,\varepsilon }^{q-2}\,dx$ \\ 
$+\delta Lm^{2}\int\limits_{\Omega }\,\eta ^{2}\,\,\,\left\vert \triangle
_{h,s}\nabla u_{\varepsilon }\right\vert ^{2}\,\,W_{2,\varepsilon }^{q-2}\,dx
$ \\ 
$+\frac{Lm^{2}}{\delta }\int\limits_{\Omega }\left\vert \triangle
_{h,s}u_{\varepsilon }\right\vert ^{2}\,\left\vert \nabla \eta \right\vert
^{2}\,\,W_{2,\varepsilon }^{q-2}\,dx$ \\ 
$+Lm^{2}\int\limits_{\Omega }\left\vert \triangle _{h,s}u_{\varepsilon
}\right\vert \,^{2}\,\left\vert \nabla \eta \right\vert
^{2}\,\,\,W_{2,\varepsilon }^{q-2}\,dx$ \\ 
$+Lm^{2}\int\limits_{\Omega }\,\eta \,^{2}\,\,W_{2,\varepsilon }^{q-2}\,dx$%
\end{tabular}%
\end{equation}

Let us consider the following second integral of (4.28) 
\begin{equation}
2\varepsilon \left\vert \sum\limits_{\alpha =1}^{m}\int\limits_{\Omega
}\triangle _{h,s}u_{\varepsilon }^{\alpha }\,\eta \,\nabla \eta
\,\,\int\limits_{0}^{1}\,\partial _{\xi ^{\alpha }}G(x,t)\,\frac{th\triangle
_{h,s}\nabla u_{\varepsilon }^{\alpha }}{\left( \left\vert \nabla
u_{\varepsilon }^{\alpha }+th\triangle _{h,s}\nabla u_{\varepsilon }^{\alpha
}\right\vert ^{2}+\varepsilon \right) ^{\frac{3}{2}}}\,dt\,dx\right\vert 
\end{equation}%
then, integrating by parts, we obtain%
\begin{equation}
\begin{tabular}{l}
$=2\varepsilon \left\vert \sum\limits_{\alpha =1}^{m}\int\limits_{\Omega
}\triangle _{h,s}u_{\varepsilon }^{\alpha }\,\eta \,\nabla \eta \,\,\left[
\partial _{\xi ^{\alpha }}G(x,1)\frac{\nabla u_{\varepsilon }^{\alpha
}+h\triangle _{h,s}\nabla u_{\varepsilon }^{\alpha }}{\left( \left\vert
\nabla u_{\varepsilon }^{\alpha }+h\triangle _{h,s}\nabla u_{\varepsilon
}^{\alpha }\right\vert ^{2}+\varepsilon \right) ^{\frac{1}{2}}}-\partial
_{\xi ^{\alpha }}G(x,0)\frac{\nabla u_{\varepsilon }^{\alpha }}{\left(
\left\vert \nabla u_{\varepsilon }^{\alpha }\right\vert ^{2}+\varepsilon
\right) ^{\frac{1}{2}}}\right] \,dx-\right. $ \\ 
$-\left. 2\varepsilon \sum\limits_{\alpha =1}^{m}\int\limits_{\Omega
}\triangle _{h,s}u_{\varepsilon }^{\alpha }\,\eta \,\nabla \eta
\,\,\,\int\limits_{0}^{1}\,\frac{d}{dt}\,\left( \partial _{\xi ^{\alpha
}}G(x,t)\right) \,\frac{\nabla u_{\varepsilon }^{\alpha }+th\triangle
_{h,s}\nabla u_{\varepsilon }^{\alpha }}{\left( \left\vert \nabla
u_{\varepsilon }^{\alpha }+th\triangle _{h,s}\nabla u_{\varepsilon }^{\alpha
}\right\vert ^{2}+\varepsilon \right) ^{\frac{1}{2}}}\,dt\,dx\right\vert $
\\ 
$\leq 2\varepsilon \sum\limits_{\alpha =1}^{m}\left\vert \int\limits_{\Omega
}\triangle _{h,s}u_{\varepsilon }^{\alpha }\,\eta \,\nabla \eta \,\,\partial
_{\xi ^{\alpha }}G(x,1)\frac{\nabla u_{\varepsilon }^{\alpha }+h\triangle
_{h,s}\nabla u_{\varepsilon }^{\alpha }}{\left( \left\vert \nabla
u_{\varepsilon }^{\alpha }+h\triangle _{h,s}\nabla u_{\varepsilon }^{\alpha
}\right\vert ^{2}+\varepsilon \right) ^{\frac{1}{2}}}\,dx\right\vert
+2\varepsilon \left\vert \int\limits_{\Omega }\triangle _{h,s}u_{\varepsilon
}^{\alpha }\,\eta \,\nabla \eta \,\,\partial _{\xi ^{\alpha }}G(x,0)\frac{%
\nabla u_{\varepsilon }^{\alpha }}{\left( \left\vert \nabla u_{\varepsilon
}^{\alpha }\right\vert ^{2}+\varepsilon \right) ^{\frac{1}{2}}}%
\,dx\right\vert $ \\
\end{tabular}%
\end{equation}%
Since%
\begin{equation}
\begin{tabular}{l}
$\frac{d}{dt}\,\left( \partial _{\xi ^{\alpha }}G(x,t)\right) $ \\ 
$=\sum\limits_{\beta =1}^{m}\,\partial _{s^{\beta }\xi ^{\alpha }}G(x,t)%
\frac{u_{\varepsilon }^{\beta }+th\triangle _{h,s}u_{\varepsilon }^{\beta }}{%
\left( \left\vert u_{\varepsilon }^{\beta }+th\triangle _{h,s}u_{\varepsilon
}^{\beta }\right\vert ^{2}+\varepsilon \right) ^{\frac{1}{2}}}h\triangle
_{h,s}u_{\varepsilon }^{\beta }$ \\ 
$+\sum\limits_{\beta =1}^{m}\,\partial _{\xi ^{\beta }\xi ^{\alpha }}G(x,t)%
\frac{\nabla u_{\varepsilon }^{\beta }+th\triangle _{h,s}\nabla
u_{\varepsilon }^{\beta }}{\left( \left\vert \nabla u_{\varepsilon }^{\beta
}+th\triangle _{h,s}\nabla u_{\varepsilon }^{\beta }\right\vert
^{2}+\varepsilon \right) ^{\frac{1}{2}}}h\triangle _{h,s}\nabla
u_{\varepsilon }^{\beta }$%
\end{tabular}%
\end{equation}%
then it follows%
\begin{equation}
\begin{tabular}{l}
$2\varepsilon \sum\limits_{\alpha =1}^{m}\left\vert \int\limits_{\Omega
}\triangle _{h,s}u_{\varepsilon }^{\alpha }\,\eta \,\nabla \eta
\,\,\int\limits_{0}^{1}\,\frac{d}{dt}\,\left( \partial _{\xi ^{\alpha
}}G(x,t)\right) \,\frac{\nabla u_{\varepsilon }^{\alpha }+th\triangle
_{h,s}\nabla u_{\varepsilon }^{\alpha }}{\left( \left\vert \nabla
u_{\varepsilon }^{\alpha }+th\triangle _{h,s}\nabla u_{\varepsilon }^{\alpha
}\right\vert ^{2}+\varepsilon \right) ^{\frac{1}{2}}}\,dt\,dx\right\vert $
\\ 
$\leq 2\varepsilon \sum\limits_{\alpha =1}^{m}\left\vert \int\limits_{\Omega
}\triangle _{h,s}u_{\varepsilon }^{\alpha }\,\eta \,\nabla \eta
\,\,\int\limits_{0}^{1}\,\sum\limits_{\beta =1}^{m}\,\partial _{s^{\beta
}\xi ^{\alpha }}G(x,t)\frac{u_{\varepsilon }^{\beta }+th\triangle
_{h,s}u_{\varepsilon }^{\beta }}{\left( \left\vert u_{\varepsilon }^{\beta
}+th\triangle _{h,s}u_{\varepsilon }^{\beta }\right\vert ^{2}+\varepsilon
\right) ^{\frac{1}{2}}}h\triangle _{h,s}u_{\varepsilon }^{\beta }\,\frac{%
\nabla u_{\varepsilon }^{\alpha }+th\triangle _{h,s}\nabla u_{\varepsilon
}^{\alpha }}{\left( \left\vert \nabla u_{\varepsilon }^{\alpha }+th\triangle
_{h,s}\nabla u_{\varepsilon }^{\alpha }\right\vert ^{2}+\varepsilon \right)
^{\frac{1}{2}}}\,dt\,dx\right\vert $ \\
\end{tabular}%
\end{equation}%
Moreover, by hypothesis H.2, we get%
\begin{equation}
\begin{tabular}{l}
$2\varepsilon \sum\limits_{\alpha =1}^{m}\left\vert \int\limits_{\Omega
}\triangle _{h,s}u_{\varepsilon }^{\alpha }\,\eta \,\nabla \eta
\,\,\int\limits_{0}^{1}\,\frac{d}{dt}\,\left( \partial _{\xi ^{\alpha
}}G(x,t)\right) \,\frac{\nabla u_{\varepsilon }^{\alpha }+th\triangle
_{h,s}\nabla u_{\varepsilon }^{\alpha }}{\left( \left\vert \nabla
u_{\varepsilon }^{\alpha }+th\triangle _{h,s}\nabla u_{\varepsilon }^{\alpha
}\right\vert ^{2}+\varepsilon \right) ^{\frac{1}{2}}}\,dt\,dx\right\vert $
\\ 
$\leq 2\varepsilon \left\vert h\right\vert \sum\limits_{\alpha
=1}^{m}\int\limits_{\Omega }\left\vert \triangle _{h,s}u_{\varepsilon
}^{\alpha }\right\vert \,\eta \,\left\vert \nabla \eta \right\vert
\,\,\int\limits_{0}^{1}\,\sum\limits_{\beta =1}^{m}\,\left\vert \partial
_{s^{\beta }\xi ^{\alpha }}G(x,t)\right\vert \left\vert \triangle
_{h,s}u_{\varepsilon }^{\beta }\right\vert \,\,dt\,dx$ \\ 
$+2\varepsilon \left\vert h\right\vert \sum\limits_{\alpha
=1}^{m}\int\limits_{\Omega }\left\vert \triangle _{h,s}u_{\varepsilon
}^{\alpha }\,\right\vert \eta \,\left\vert \nabla \eta \right\vert
\,\,\int\limits_{0}^{1}\,\sum\limits_{\beta =1}^{m}\,\left\vert \partial
_{\xi ^{\beta }\xi ^{\alpha }}G(x,t)\right\vert \left\vert \triangle
_{h,s}\nabla u_{\varepsilon }^{\beta }\right\vert \,\,dt\,dx$ \\ 
$\leq 2\varepsilon \left\vert h\right\vert Lm^{q-2}\sum\limits_{\alpha
,\beta =1}^{m}\int\limits_{\Omega }\left\vert \triangle _{h,s}u_{\varepsilon
}^{\alpha }\right\vert \,\eta \,\left\vert \nabla \eta \right\vert \,$ \\ 
$\cdot \int\limits_{0}^{1}\left[ \left( \left\vert u_{\varepsilon
}+th\triangle _{h,s}u_{\varepsilon }\right\vert ^{2}+\varepsilon \right) ^{%
\frac{q-e}{2}}+\left( \left\vert \nabla u_{\varepsilon }+th\triangle
_{h,s}\nabla u_{\varepsilon }\right\vert ^{2}+\varepsilon \right) ^{\frac{q-2%
}{2}}\right] \left\vert \triangle _{h,s}u_{\varepsilon }^{\beta }\right\vert
\,\,dt\,dx$ \\
\end{tabular}%
\end{equation}%
Using () and Lemma 8, since $\left\vert h\right\vert <1$ and $0<\varepsilon
<1$, it follows%
\begin{equation}
\begin{tabular}{l}
$2\varepsilon \sum\limits_{\alpha =1}^{m}\left\vert \int\limits_{\Omega
}\triangle _{h,s}u_{\varepsilon }^{\alpha }\,\eta \,\nabla \eta
\,\,\int\limits_{0}^{1}\,\frac{d}{dt}\,\left( \partial _{\xi ^{\alpha
}}G(x,t)\right) \,\frac{\nabla u_{\varepsilon }^{\alpha }+th\triangle
_{h,s}\nabla u_{\varepsilon }^{\alpha }}{\left( \left\vert \nabla
u_{\varepsilon }^{\alpha }+th\triangle _{h,s}\nabla u_{\varepsilon }^{\alpha
}\right\vert ^{2}+\varepsilon \right) ^{\frac{1}{2}}}\,dt\,dx\right\vert $
\\ 
$\leq 4Lm^{q}\int\limits_{\Omega }\left\vert \triangle _{h,s}u_{\varepsilon
}\right\vert ^{2}\,\eta \,\left\vert \nabla \eta \right\vert \,\,\left[
W_{1,\varepsilon }^{q-2}+W_{2,\varepsilon }^{q-2}\right] \,\,dx$ \\ 
$+4Lm^{q}\int\limits_{\Omega }\left\vert \triangle _{h,s}u_{\varepsilon
}\right\vert \,\eta \,\left\vert \nabla \eta \right\vert \left\vert
\triangle _{h,s}\nabla u_{\varepsilon }\right\vert \,\,\left[
W_{1,\varepsilon }^{q-2}+W_{2,\varepsilon }^{q-2}\right] \,\,dx$%
\end{tabular}%
\end{equation}%
then by Young Inequality we obtain%
\begin{equation}
\begin{tabular}{l}
$2\varepsilon \sum\limits_{\alpha =1}^{m}\left\vert \int\limits_{\Omega
}\triangle _{h,s}u_{\varepsilon }^{\alpha }\,\eta \,\nabla \eta
\,\,\int\limits_{0}^{1}\,\frac{d}{dt}\,\left( \partial _{\xi ^{\alpha
}}G(x,t)\right) \,\frac{\nabla u_{\varepsilon }^{\alpha }+th\triangle
_{h,s}\nabla u_{\varepsilon }^{\alpha }}{\left( \left\vert \nabla
u_{\varepsilon }^{\alpha }+th\triangle _{h,s}\nabla u_{\varepsilon }^{\alpha
}\right\vert ^{2}+\varepsilon \right) ^{\frac{1}{2}}}\,dt\,dx\right\vert $
\\ 
$\leq 2Lm^{q}\int\limits_{\Omega }\left\vert \triangle _{h,s}u_{\varepsilon
}\right\vert ^{2}\,\eta ^{2}\,\,\left[ W_{1,\varepsilon
}^{q-2}+W_{2,\varepsilon }^{q-2}\right] \,\,dx$ \\ 
$+2Lm^{q}\int\limits_{\Omega }\left\vert \triangle _{h,s}u_{\varepsilon
}\right\vert ^{2}\,\left\vert \nabla \eta \right\vert ^{2}\,\,\left[
W_{1,\varepsilon }^{q-2}+W_{2,\varepsilon }^{q-2}\right] \,\,dx$ \\ 
$+2\tilde{\delta}_{3}Lm^{q}\int\limits_{\Omega }\,\eta ^{2}\,\left\vert
\triangle _{h,s}\nabla u_{\varepsilon }\right\vert ^{2}\,\,\left[
W_{1,\varepsilon }^{q-2}+W_{2,\varepsilon }^{q-2}\right] \,\,dx$ \\ 
$+\frac{2Lm^{q}}{\tilde{\delta}_{3}}\int\limits_{\Omega }\left\vert
\triangle _{h,s}u_{\varepsilon }\right\vert ^{2}\,\,\left\vert \nabla \eta
\right\vert ^{2}\,\,\left[ W_{1,\varepsilon }^{q-2}+W_{2,\varepsilon }^{q-2}%
\right] \,\,dx$%
\end{tabular}%
\end{equation}%
Using (4.37) and (4.44)\ we have%
\begin{equation}
\begin{tabular}{l}
$\Theta _{3}$ \\ 
$=\left\vert \int\limits_{0}^{1}\int\limits_{\Omega }\sum\limits_{\alpha
=1}^{m}\frac{2\varepsilon \,\partial _{\xi ^{\alpha }}G(x,t)}{\left(
\left\vert \nabla u_{\varepsilon }^{\alpha }+th\triangle _{h,s}\nabla
u_{\varepsilon }^{\alpha }\right\vert ^{2}+\varepsilon \right) ^{\frac{3}{2}}%
}\,\,\triangle _{h,s}\nabla u_{\varepsilon }^{\alpha }\,\triangle
_{h,s}u_{\varepsilon }^{\alpha }\,\eta \,\nabla \eta \,dx\,dt\right\vert $
\\ 
$\leq \left( \tilde{\delta}_{2}Lm^{2}+2\tilde{\delta}_{3}Lm^{q}\right)
\int\limits_{\Omega }\,\eta ^{2}\,\left\vert \triangle _{h,s}\nabla
u_{\varepsilon }\right\vert ^{2}\,\,\left[ W_{1,\varepsilon
}^{q-2}+W_{2,\varepsilon }^{q-2}\right] \,\,dx$ \\ 
$+\left( 4Lm^{q}+\frac{Lm^{2}}{\tilde{\delta}_{2}}+\frac{2Lm^{q}}{\tilde{%
\delta}_{3}}\right) \int\limits_{\Omega }\left\vert \triangle
_{h,s}u_{\varepsilon }\right\vert ^{2}\,\,\left\vert \nabla \eta \right\vert
^{2}\,\left[ W_{1,\varepsilon }^{q-2}+W_{2,\varepsilon }^{q-2}\right] \,dx$
\\ 
$+2Lm^{q}\int\limits_{\Omega }\left\vert \triangle _{h,s}u_{\varepsilon
}\right\vert ^{2}\,\eta ^{2}\,\,\left[ W_{1,\varepsilon
}^{q-2}+W_{2,\varepsilon }^{q-2}\right] \,\,dx$ \\ 
$+2Lm^{2}\int\limits_{\Omega }\left\vert \triangle _{h,s}u_{\varepsilon
}\right\vert ^{2}\,\eta \,\left\vert \nabla \eta \right\vert
\,\,W_{1,\varepsilon }^{q-2}\,dx$ \\ 
$+Lm^{2}\int\limits_{\Omega }\,\eta \,^{2}\,\left[ W_{1,\varepsilon
}^{q-2}+W_{2,\varepsilon }^{q-2}\right] \,dx$ \\ 
$+Lm^{2}\int\limits_{\Omega }\eta \,^{2}\,\left\vert u_{\varepsilon
}\right\vert ^{2}W_{1,\varepsilon }^{q-2}\,dx$ \\ 
$+Lm^{2}\int\limits_{\Omega }\eta ^{2}\,\,\left\vert \nabla u_{\varepsilon
}\right\vert ^{2}\,\,W_{2,\varepsilon }^{q-2}\,dx$%
\end{tabular}%
\end{equation}%
Using (4.16), (4.20), (4.21), (4.22), (4.24), (4.25) and (4.45) we get

\begin{equation}
\begin{tabular}{l}
$\left( p-1\right) \int\limits_{\Omega }\varepsilon \,\sum\limits_{\alpha
=1}^{m}I_{1,\alpha }\left( \triangle _{h,s}\left( \nabla u_{\varepsilon
}^{\alpha }\right) \right) ^{2}\eta ^{2}dx+\nu \vartheta \int\limits_{\Omega
}\left[ W_{1,\varepsilon }^{q-2}+W_{2,\varepsilon }^{q-2}\right] \left(
\left\vert \triangle _{h,s}\left( \nabla u_{\varepsilon }\right) \right\vert
^{2}+\left\vert \triangle _{h,s}u_{\varepsilon }\right\vert ^{2}\right) \eta
^{2}dx$ \\ 
$\leq \delta \left( p-1\right) \int\limits_{\Omega }\varepsilon
\sum\limits_{\alpha =1}^{m}I_{1,\alpha }\,\,\left\vert \triangle
_{h,s}\left( \nabla u_{\varepsilon }^{\alpha }\right) \right\vert ^{2}\eta
^{2}\,dx+\frac{\left( p-1\right) }{\delta }\int\limits_{\Omega }\varepsilon
\sum\limits_{\alpha =1}^{m}I_{1,\alpha }\,\left\vert \triangle
_{h,s}u_{\varepsilon }^{\alpha }\right\vert ^{2}\,\left\vert \nabla \eta
\right\vert ^{2}\,dx$ \\ 
$+\left( pLm^{2}+2Lm^{2}\right) \int\limits_{\Omega }\left[ W_{1,\varepsilon
}^{q-2}+W_{2,\varepsilon }^{q-2}\right] \left\vert \triangle
_{h,s}u_{\varepsilon }\right\vert ^{2}\,\,\eta \,\left\vert \nabla \eta
\right\vert \,dx$ \\ 
$+\left( \tilde{\delta}\frac{pLm^{2}}{2}+\tilde{\delta}_{2}Lm^{2}+2\tilde{%
\delta}_{3}Lm^{q}\right) \int\limits_{\Omega }\left[ W_{1,\varepsilon
}^{q-2}+W_{2,\varepsilon }^{q-2}\right] \left( \left\vert \triangle
_{h,s}\left( \nabla u_{\varepsilon }\right) \right\vert ^{2}+\left\vert
\triangle _{h,s}u_{\varepsilon }\right\vert ^{2}\right) \eta ^{2}dx$ \\ 
$+\left( c_{q}\frac{pLm^{2}}{2\tilde{\delta}}+4Lm^{q}+\frac{Lm^{2}}{\tilde{%
\delta}_{2}}+\frac{2Lm^{q}}{\tilde{\delta}_{3}}\right) \int\limits_{\Omega }%
\left[ W_{1,\varepsilon }^{q-2}+W_{2,\varepsilon }^{q-2}\right] \,\left\vert
\triangle _{h,s}u_{\varepsilon }\right\vert ^{2}\left\vert \nabla \eta
\right\vert ^{2}\,dx$ \\
\end{tabular}%
\end{equation}%
If we choose $\delta =\frac{1}{2}$, $\tilde{\delta}=\frac{\nu \vartheta }{%
3pLm^{2}}$, $\tilde{\delta}_{2}=\frac{\nu \vartheta }{6Lm^{2}}$\ and $\tilde{%
\delta}_{3}=\frac{\nu \vartheta }{12Lm^{q}}$\ it follows%
\begin{equation}
\begin{tabular}{l}
$\frac{\left( p-1\right) }{2}\int\limits_{\Omega }\varepsilon
\,\sum\limits_{\alpha =1}^{m}I_{1,\alpha }\left( \triangle _{h,s}\left(
\nabla u_{\varepsilon }^{\alpha }\right) \right) ^{2}\eta ^{2}dx+\frac{\nu
\vartheta }{2}\int\limits_{\Omega }\left[ W_{1,\varepsilon
}^{q-2}+W_{2,\varepsilon }^{q-2}\right] \left( \left\vert \triangle
_{h,s}\left( \nabla u_{\varepsilon }\right) \right\vert ^{2}+\left\vert
\triangle _{h,s}u_{\varepsilon }\right\vert ^{2}\right) \eta ^{2}dx$ \\ 
$\leq 2\left( p-1\right) \int\limits_{\Omega }\varepsilon
\sum\limits_{\alpha =1}^{m}I_{1,\alpha }\,\left\vert \triangle
_{h,s}u_{\varepsilon }^{\alpha }\right\vert ^{2}\,\left\vert \nabla \eta
\right\vert ^{2}\,dx+\left( pLm^{2}+2Lm^{2}\right) \int\limits_{\Omega }%
\left[ W_{1,\varepsilon }^{q-2}+W_{2,\varepsilon }^{q-2}\right] \left\vert
\triangle _{h,s}u_{\varepsilon }\right\vert ^{2}\,\,\eta \,\left\vert \nabla
\eta \right\vert \,dx$ \\ 
$+\left[ c_{q}\frac{\left( pLm^{2}\right) ^{2}}{2\nu \vartheta }++4Lm^{q}+%
\frac{6\left( Lm^{2}\right) ^{2}}{\nu \vartheta }+\frac{24\left(
Lm^{q}\right) ^{2}}{\nu \vartheta }\right] \int\limits_{\Omega }\left[
W_{1,\varepsilon }^{q-2}+W_{2,\varepsilon }^{q-2}\right] \,\left\vert
\triangle _{h,s}u_{\varepsilon }\right\vert ^{2}\left\vert \nabla \eta
\right\vert ^{2}\,dx$ \\ 
$+2Lm^{q}\int\limits_{\Omega }\left\vert \triangle _{h,s}u_{\varepsilon
}\right\vert ^{2}\,\eta ^{2}\,\,\left[ W_{1,\varepsilon
}^{q-2}+W_{2,\varepsilon }^{q-2}\right] \,\,dx$ \\ 
$+Lm^{2}\int\limits_{\Omega }\,\eta \,^{2}\,\left[ W_{1,\varepsilon
}^{q-2}+W_{2,\varepsilon }^{q-2}\right] \,dx$ \\ 
$+Lm^{2}\int\limits_{\Omega }\eta \,^{2}\,\left\vert u_{\varepsilon
}\right\vert ^{2}W_{1,\varepsilon }^{q-2}\,dx$ \\ 
$+Lm^{2}\int\limits_{\Omega }\eta ^{2}\,\,\left\vert \nabla u_{\varepsilon
}\right\vert ^{2}\,\,W_{2,\varepsilon }^{q-2}\,dx$%
\end{tabular}%
\end{equation}%
then, by Lemma 8 and Young Inequality, we have%
\begin{equation}
%
%
\end{equation}%
Recalling the properties of the function $\eta $ we obtain the following two
inequalities%
\begin{equation}
\begin{tabular}{l}
$\frac{\vartheta \nu \varepsilon ^{q}}{2}\int\limits_{B_{s}\left(
x_{0}\right) }\left\vert \triangle _{h,s}\left( \nabla u_{\varepsilon
}\right) \right\vert ^{2}dx$ \\ 
$\leq \frac{2\left( p-1\right) \varepsilon }{\left( t-s\right) ^{2}}%
\int\limits_{B_{t}\left( x_{0}\right) }\left[ \left\vert \nabla
u_{\varepsilon }\left( x\right) \right\vert ^{2}+\left\vert \nabla
u_{\varepsilon }\left( x+h\right) \right\vert ^{2}+\varepsilon \right] ^{%
\frac{p-2}{2}}\,\left\vert \triangle _{h,s}u_{\varepsilon }\right\vert
^{2}\,\,dx$ \\ 
$+\left[ \frac{pLm^{2}+2Lm^{2}}{\left( t-s\right) }+\frac{C_{q,m}}{\left(
t-s\right) ^{2}}\right] \int\limits_{B_{t}\left( x_{0}\right) }\left[
W_{1,\varepsilon }^{q-2}+W_{2,\varepsilon }^{q-2}\right] \left\vert
\triangle _{h,s}u_{\varepsilon }\right\vert ^{2}\,\,\,dx$ \\ 
$+2Lm^{q}\int\limits_{B_{t}\left( x_{0}\right) }\left\vert \triangle
_{h,s}u_{\varepsilon }\right\vert ^{q}\,\,\,\,dx+\left(
4Lm^{2}+2Lm^{q}\right) \int\limits_{B_{t}\left( x_{0}\right) }\,\,\,\left[
W_{1,\varepsilon }^{q}+W_{2,\varepsilon }^{q}\right] \,\,dx$ \\ 
$+Lm^{2}\int\limits_{B_{t}\left( x_{0}\right)
}\,\,\,dx+Lm^{2}\int\limits_{B_{t}\left( x_{0}\right) }\,\left( \left\vert
u_{\varepsilon }\right\vert ^{q}\,+\,\,\left\vert \nabla u_{\varepsilon
}\right\vert ^{q}\right) \,\,\,dx$%
\end{tabular}%
%
\end{equation}%
where $C_{q,m}=\left[ c_{q}\frac{\left( pLm^{2}\right) ^{2}}{2\nu \vartheta }%
++4Lm^{q}+\frac{6\left( Lm^{2}\right) ^{2}}{\nu \vartheta }+\frac{24\left(
Lm^{q}\right) ^{2}}{\nu \vartheta }\right] $ . Using (4.49) and Lemma 7 we
deduce that $u_{\varepsilon }\in W_{loc}^{2,2}\left( \Omega \right) $;
moreover using (4.50) and Lemma 7, sending $h\rightarrow 0$ we get%
\begin{equation}
%
%
\end{equation}%
where $G_{q,m,l}>0$, then, by Young inequality, it follows that 
\begin{equation}
%
%
\end{equation}%
and%
\begin{equation}
\int\limits_{B_{t}\left( x_{0}\right) }\left[ \left[ \left\vert u\left(
x\right) \right\vert ^{2}+\varepsilon \right] ^{\frac{q-2}{2}}+\left[
\left\vert \nabla u_{\varepsilon }\left( x\right) \right\vert
^{2}+\varepsilon \right] ^{\frac{q-2}{2}}\right] \left\vert \nabla
u_{\varepsilon }\right\vert ^{2}\,\,\,dx\leq C_{q}\int\limits_{B_{t}\left(
x_{0}\right) }1+\left\vert u_{\varepsilon }\right\vert ^{q}+\left\vert
\nabla u_{\varepsilon }\right\vert ^{q}\,\,\,dx
\end{equation}%
Moreover by (4.50), (4.51), (4.52) and (4.53) we get%
\begin{equation}
\int\limits_{B_{s}\left( x_{0}\right) }\left[ \left\vert \nabla
u_{\varepsilon }\left( x\right) \right\vert ^{2}+\varepsilon \right] ^{\frac{%
q-2}{2}}\left\vert Hu_{\varepsilon }\right\vert ^{2}dx\leq \frac{C_{p,q}}{%
\left( t-s\right) ^{2}}\int\limits_{B_{t}\left( x_{0}\right) }\varepsilon
\left\vert \nabla u_{\varepsilon }\right\vert ^{p}+\left\vert u_{\varepsilon
}\right\vert ^{q}+\left\vert \nabla u_{\varepsilon }\right\vert ^{q}+1\,\,dx
\end{equation}

\section{Proof of Theorem 5}

In this section we prove the higher local integrability of the\ modulus of
gradient of the local minima of the functional (1.5).

Let us consider $V=\left\vert \nabla u_{\varepsilon }\right\vert ^{\frac{q}{2%
}}$ then%
\begin{equation}
\begin{tabular}{l}
$\left\vert \nabla V\right\vert ^{2}$ \\ 
$\leq \left( \frac{nmq}{2}\right) ^{2}\left\vert \nabla u_{\varepsilon
}\right\vert ^{q-2}\left\vert Hu_{\varepsilon }\right\vert ^{2}$ \\ 
$\leq \left( \frac{nmq}{2}\right) ^{2}\left[ \left\vert \nabla
u_{\varepsilon }\right\vert ^{2}+\varepsilon \right] ^{\frac{q-2}{2}%
}\left\vert Hu_{\varepsilon }\right\vert ^{2}$%
\end{tabular}%
\end{equation}%
and by Theorem 4 it follows 
\begin{equation}
\int\limits_{B_{2r}\left( x_{0}\right) }\left\vert \nabla V\right\vert
^{2}dx\leq \frac{\left( nmq\right) ^{2}C_{p,q}}{r^{2}}\int\limits_{B_{4r}%
\left( x_{0}\right) }\varepsilon \left\vert \nabla u_{\varepsilon
}\right\vert ^{p}+\left\vert u_{\varepsilon }\right\vert ^{q}+\left\vert
\nabla u_{\varepsilon }\right\vert ^{q}+1\,\,dx
\end{equation}%
where $\left( nmq\right) ^{2}C_{p,q}$ is a positive real constan depending
only on $p$,$q$, $n$, $m$.

Let us define $W=\eta V$ where $\eta \in C_{0}^{\infty }\left( B_{2r}\left(
x_{0}\right) \right) $, $0\leq \eta \leq 1$ in $B_{2r}\left( x_{0}\right) $, 
$\eta =1$ in $B_{r}\left( x_{0}\right) $ and $\left\vert \nabla \eta
\right\vert \leq \frac{2}{r}$ in $B_{2r}\left( x_{0}\right) $; then by
Sobolev Inequality we get%
\begin{equation*}
\left[ \int\limits_{B_{2r}\left( x_{0}\right) }\left\vert W\right\vert
^{2^{\ast }}\,dx\right] ^{\frac{1}{2^{\ast }}}\leq c_{n}\left[
\int\limits_{B_{2r}\left( x_{0}\right) }\left\vert \nabla W\right\vert
^{2}\,dx\right] ^{\frac{1}{2}}
\end{equation*}%
moreover, since, $\nabla W=\eta \nabla V+\nabla \eta V$; then%
\begin{equation*}
\left[ \int\limits_{B_{2r}\left( x_{0}\right) }\left\vert W\right\vert
^{2^{\ast }}\,dx\right] ^{\frac{1}{2^{\ast }}}\leq c_{n}\left[ \left( \frac{2%
}{r}\right) ^{2}\int\limits_{B_{2r}\left( x_{0}\right) }\left\vert
V\right\vert ^{2}\,dx\right] ^{\frac{1}{2}}+c_{n}\left[ \int\limits_{B_{2r}%
\left( x_{0}\right) }\left\vert \nabla V\right\vert ^{2}\,dx\right] ^{\frac{1%
}{2}}
\end{equation*}%
Using (5.2), since $W=V$ in $B_{r}\left( x_{0}\right) $, it follows 
\begin{equation}
\begin{tabular}{l}
$\left[ \int\limits_{B_{r}\left( x_{0}\right) }\left\vert \nabla
u_{\varepsilon }\right\vert ^{\frac{2^{\ast }}{2}q}\,dx\right] ^{\frac{1}{%
2^{\ast }}}$ \\ 
$\leq c_{n}\left[ \left( \frac{2}{r}\right) ^{2}\int\limits_{B_{2r}\left(
x_{0}\right) }\left\vert \nabla u_{\varepsilon }\right\vert ^{q}\,dx\right]
^{\frac{1}{2}}$ \\ 
$+c_{n}\left[ \frac{\left( nmq\right) ^{2}C_{p,q}}{r^{2}}\int\limits_{B_{4r}%
\left( x_{0}\right) }\varepsilon \left\vert \nabla u_{\varepsilon
}\right\vert ^{p}+\left\vert u_{\varepsilon }\right\vert ^{q}+\left\vert
\nabla u_{\varepsilon }\right\vert ^{q}+1\,\,dx\right] ^{\frac{1}{2}}$%
\end{tabular}%
\end{equation}%
then $\left\vert \nabla u_{\varepsilon }\right\vert \in L_{loc}^{\frac{%
2^{\ast }}{2}q}\left( \Omega \right) $.

\bigskip

\section{Proof of Theorem 6}

The following inequality is an intresting consequence of the higher local
integrability theorem [Theorem 6];\ let us take $2\leq q<p<\min \left\{ 
\frac{2^{\ast }}{2}q,q^{\ast }\right\} $ then%
\begin{equation}
\int\limits_{B_{r}\left( x_{0}\right) }\left\vert \nabla u_{\varepsilon
}\right\vert ^{p}\,dx\leq \left\vert B_{r}\left( x_{0}\right) \right\vert
^{1-\frac{2p}{2^{\ast }q}}\left[ \int\limits_{B_{r}\left( x_{0}\right)
}\left\vert \nabla u_{\varepsilon }\right\vert ^{\frac{2^{\ast }}{2}q}\,dx%
\right] ^{\frac{2p}{2^{\ast }q}}
\end{equation}%
and by (5.3)%
\begin{equation}
\begin{tabular}{l}
$\left[ \int\limits_{B_{r}\left( x_{0}\right) }\left\vert \nabla
u_{\varepsilon }\right\vert ^{\frac{2^{\ast }}{2}q}\,dx\right] ^{\frac{2}{%
2^{\ast }}}$ \\ 
$\leq 2c_{n}\left[ \left( \frac{2}{r}\right) ^{2}\int\limits_{B_{2r}\left(
x_{0}\right) }\left\vert \nabla u_{\varepsilon }\right\vert ^{q}\,dx\right] $
\\ 
$+2c_{n}\left[ \frac{\left( nmq\right) ^{2}C_{p,q}}{r^{2}}%
\int\limits_{B_{4r}\left( x_{0}\right) }\varepsilon \left\vert \nabla
u_{\varepsilon }\right\vert ^{p}+\left\vert u_{\varepsilon }\right\vert
^{q}+\left\vert \nabla u_{\varepsilon }\right\vert ^{q}+1\,\,dx\right] $%
\end{tabular}%
\end{equation}%
Using hypothesis H.1 it follows%
\begin{equation}
\int\limits_{B_{2r}\left( x_{0}\right) }\left\vert \nabla u_{\varepsilon
}\right\vert ^{q}\,dx\leq \mathcal{F}_{\varepsilon }\left( u_{\varepsilon
},B_{2r}\left( x_{0}\right) \right) ,\int\limits_{B_{4r}\left( x_{0}\right)
}\varepsilon \left\vert \nabla u_{\varepsilon }\left( x\right) \right\vert
^{p}\,\,dx\leq \mathcal{F}_{\varepsilon }\left( u_{\varepsilon
},B_{2r}\left( x_{0}\right) \right)
\end{equation}%
and 
\begin{equation}
\int\limits_{B_{4r}\left( x_{0}\right) }1+\left\vert u_{\varepsilon }\left(
x\right) \right\vert ^{q}+\left\vert \nabla u_{\varepsilon }\left( x\right)
\right\vert ^{q}\,\,dx\leq \mathcal{F}_{\varepsilon }\left( u_{\varepsilon
},B_{2r}\left( x_{0}\right) \right) +\left\vert B_{2r}\left( x_{0}\right)
\right\vert
\end{equation}%
then by (6.2), (6.3) and (6.4) we get%
\begin{equation}
\left[ \int\limits_{B_{r}\left( x_{0}\right) }\left\vert \nabla
u_{\varepsilon }\right\vert ^{\frac{2^{\ast }}{2}q}\,dx\right] ^{\frac{2}{%
2^{\ast }}}\leq 6c_{n}\left[ \left( \frac{2}{r}\right) ^{2}+C_{1}\left( 
\frac{1}{r}\right) ^{2}+C_{2}\left( \frac{1}{r}\right) ^{2}\right] \left( 
\mathcal{F}_{\varepsilon }\left( u_{\varepsilon },B_{2r}\left( x_{0}\right)
\right) +\left\vert B_{2r}\left( x_{0}\right) \right\vert \right)
\end{equation}%
Using (6.1), (6.2) and (6.5) we deduce%
\begin{equation}
\int\limits_{B_{r}\left( x_{0}\right) }\left\vert \nabla u_{\varepsilon
}\right\vert ^{p}\,dx\leq \left\vert B_{r}\left( x_{0}\right) \right\vert
^{1-\frac{2p}{2^{\ast }q}}\left[ \frac{C}{r^{2}}\left( \mathcal{F}%
_{\varepsilon }\left( u_{\varepsilon },B_{2r}\left( x_{0}\right) \right)
+\left\vert B_{2r}\left( x_{0}\right) \right\vert \right) \right] ^{\frac{p}{%
q}}
\end{equation}%
Since 
\begin{equation}
\mathcal{F}_{\varepsilon }\left( u_{\varepsilon },B_{2r}\left( x_{0}\right)
\right) \leq \mathcal{F}_{\varepsilon }\left( u_{1},B_{2r}\left(
x_{0}\right) \right) \leq F_{1}\left( u_{1},B_{2r}\left( x_{0}\right) \right)
\end{equation}%
where $u_{1}$\ is a minimum of the functional 
\begin{equation}
\mathcal{F}_{1}\left( u_{1},B_{2r}\left( x_{0}\right) \right)
=\int\limits_{\Omega }\sum\limits_{\alpha =1}^{m}\left\vert \nabla u^{\alpha
}\right\vert ^{p}+G\left( \sqrt{\left\vert u^{1}\right\vert ^{2}+1},...,%
\sqrt{\left\vert u^{m}\right\vert ^{2}+1},\sqrt{\left\vert \nabla
u^{1}\right\vert ^{2}+1},...,\sqrt{\left\vert \nabla u^{m}\right\vert ^{2}+1}%
\right) \,dx
\end{equation}%
then we get a uniform estimate of $\int\limits_{B_{r}\left( x_{0}\right)
}\left\vert \nabla u_{\varepsilon }\right\vert ^{p}\,dx$:%
\begin{equation}
\int\limits_{B_{r}\left( x_{0}\right) }\left\vert \nabla u_{\varepsilon
}\right\vert ^{p}\,dx\leq \left\vert B_{r}\left( x_{0}\right) \right\vert
^{1-\frac{2p}{2^{\ast }q}}\left[ \frac{C}{r^{2}}\left( \mathcal{F}_{1}\left(
u_{1},B_{2r}\left( x_{0}\right) \right) +\left\vert B_{2r}\left(
x_{0}\right) \right\vert \right) \right] ^{\frac{p}{q}}
\end{equation}%
Using H\"{o}lder Inequality and Sobolev's embedding Theorem we have%
\begin{equation}
\begin{tabular}{l}
$\left[ \int\limits_{B_{r}\left( x_{0}\right) }\left\vert u_{\varepsilon
}\right\vert ^{p}\,dx\right] ^{\frac{1}{p}}$ \\ 
$\leq \left[ \left\vert B_{r}\left( x_{0}\right) \right\vert ^{\frac{q^{+}-p%
}{q^{\ast }}}\left( \int\limits_{B_{r}\left( x_{0}\right) }\left\vert
u_{\varepsilon }\right\vert ^{q^{\ast }}\,dx\right) ^{\frac{p}{q^{\ast }}}%
\right] ^{\frac{1}{p}}$ \\ 
$\leq \left\vert B_{r}\left( x_{0}\right) \right\vert ^{\frac{q^{+}-p}{%
q^{\ast }}}C_{IS}\left\Vert u_{\varepsilon }\right\Vert _{W^{1,q}\left(
B_{2r}\left( x_{0}\right) \right) }$%
\end{tabular}%
\end{equation}%
Now using hypothesis H.1 and inequality (6.10) it follows%
\begin{equation}
\left[ \int\limits_{B_{r}\left( x_{0}\right) }\left\vert u_{\varepsilon
}\right\vert ^{p}\,dx\right] ^{\frac{1}{p}}\leq \left\vert B_{r}\left(
x_{0}\right) \right\vert ^{\frac{q^{+}-p}{q^{\ast }}}C_{IS}\,\mathcal{F}%
_{1}\left( u_{1},B_{2r}\left( x_{0}\right) \right)
\end{equation}%
then from (6.9) and (6.11) a reral positive constant $D_{p,q}$ exists such
that 
\begin{equation}
\left\Vert u_{\varepsilon }\right\Vert _{W^{1,p}\left( B_{r}\left(
x_{0}\right) \right) }<D_{p,q}
\end{equation}%
for every $\varepsilon \in \left( 0,1\right] $.

\section{Proof of Theorem 7}

Let us fix $2\leq q<p<\min \left\{ \frac{2^{\ast }}{2}q,q^{\ast }\right\} $
then by Theorem 6 and Theorem IV.9 of Brezis [8] a sequence $\left\{
\varepsilon _{l}\right\} _{l\in 
\mathbb{N}
}$ exists such that%
\begin{equation}
\lim\limits_{l\rightarrow +\infty }\varepsilon _{l}=0
\end{equation}%
\begin{equation}
\left\Vert u_{\varepsilon _{l}}\right\Vert _{W^{1,p}\left( B_{r}\left(
x_{0}\right) \right) }<C_{p,q}
\end{equation}%
\begin{equation}
u_{\varepsilon _{l}}\rightarrow u_{0}\qquad \text{in }L^{p}\left(
B_{r}\left( x_{0}\right) \right)
\end{equation}%
\begin{equation}
u_{\varepsilon _{l}}\rightarrow u_{0}\qquad \text{a. e. in }B_{r}\left(
x_{0}\right)
\end{equation}%
\begin{equation}
\nabla u_{\varepsilon _{l}}\rightharpoonup \nabla u_{0}\qquad \text{weak in }%
L^{p}\left( B_{r}\left( x_{0}\right) \right)
\end{equation}%
where $u_{0}\in W^{1,p}\left( B_{r}\left( x_{0}\right) \right) $, moreove%
\begin{equation}
\left\vert u_{\varepsilon _{l}}\left( x\right) \right\vert \leq H\left(
x\right) \qquad \forall l\in 
\mathbb{N}
\text{, a. e. in }B_{r}\left( x_{0}\right)
\end{equation}%
where $H\in L^{p}\left( B_{r}\left( x_{0}\right) \right) $. Since the
density $G$ is a strictly convex function then the functional (1.5) is a
lower semicontinuity functional. Lower semicontinuty and the minimality of $%
u_{\varepsilon }$ give%
\begin{equation}
\mathcal{F}_{0}\left( u_{0},B_{r}\left( x_{0}\right) \right) \leq
\liminf\limits_{j\rightarrow +\infty }\mathcal{F}_{\varepsilon _{j}}\left(
u_{\varepsilon _{l}},B_{r}\left( x_{0}\right) \right) \leq
\limsup\limits_{j\rightarrow +\infty }\mathcal{F}_{\varepsilon _{j}}\left(
u_{0}+\varphi ,B_{r}\left( x_{0}\right) \right) =\mathcal{F}_{0}\left(
u_{0}+\varphi ,B_{r}\left( x_{0}\right) \right)
\end{equation}%
for every $\varphi \in W_{0}^{1,p}\left( B_{r}\left( x_{0}\right) ,%
\mathbb{R}
^{m}\right) $. Since $G$\ is a strictly convex function then $u_{0}$ is the
unique minimum of the functional $\mathcal{F}_{0}\left( u,B_{r}\left(
x_{0}\right) \right) $. Moreover if $\tilde{u}_{0}\in W^{1,q}\left( \Omega ,%
\mathbb{R}
^{m}\right) $ is a minimum of $\mathcal{F}_{0}\left( u,\Omega \right) $ then
it is a local minimum in $B_{r}\left( x_{0}\right) $, by the strictly
convexty of $G$ it follows $\tilde{u}_{0}=u_{0}$ a. e. in $B_{r}\left(
x_{0}\right) $. So we can identify $\tilde{u}_{0}$ and $u_{0}$.

Now le us fix $y_{0}\in B_{r}\left( x_{0}\right) $ and let us take $%
0<\varrho <t<s<R<\min \left\{ \frac{r}{8},\frac{dist\left( y_{0},\partial
B_{r}\left( x_{0}\right) \right) }{8}\right\} $ then by Theorem we get

\bigskip 
\begin{equation*}
\begin{tabular}{l}
$\int\limits_{A_{k,s,y_{0}}^{1,\varepsilon _{l}}}\varepsilon _{l}\left\vert
\nabla u_{\varepsilon _{l}}^{1}\right\vert ^{p}\,+\left\vert \nabla
u_{\varepsilon _{l}}^{1}\right\vert ^{q}\,dx$ \\ 
$\leq L\int\limits_{A_{k,s,y_{0}}^{1,\varepsilon _{l}}\backslash
A_{k,t,y_{0}}^{1,\varepsilon _{l}}}\varepsilon _{l}\left\vert \nabla
u_{\varepsilon _{l}}^{1}\right\vert ^{p}+\left\vert \nabla u_{\varepsilon
_{l}}^{1}\right\vert ^{q}\,dx+\frac{2^{p}\left( p^{p}+p^{q}\right) }{\left(
s-t\right) ^{p}}\int\limits_{A_{k,s,y_{0}}^{1,\varepsilon _{l}}\backslash
A_{k,t,y_{0}}^{1,\varepsilon _{l}}}\left( u_{\varepsilon _{l}}^{1}-k\right)
^{p}\,dx+$ \\ 
$D_{B_{r}\left( x_{0}\right) ,\varepsilon _{l}}\left( 1+k^{p}\right) \left[ 
\mathcal{L}^{n}\left( A_{k,s,y_{0}}^{1,\varepsilon _{l}}\right) \right] ^{1-%
\frac{p}{n}+\epsilon }$%
\end{tabular}%
\end{equation*}%
where%
\begin{equation}
A_{k,\varrho ,y_{0}}^{1,\varepsilon _{l}}=\left\{ u_{\varepsilon
_{l}}^{1}>k\right\} \cap B_{r}\left( y_{0}\right)
\end{equation}%
and%
\begin{equation}
D_{B_{r}\left( x_{0}\right) ,\varepsilon _{l}}=2^{p-1}L\left(
p^{q}+L+2m\left( L+1\right) \left\Vert u_{\varepsilon _{l}}\right\Vert
_{W^{1,p}\left( B_{r}\left( x_{0}\right) \right) }^{q}+\left( L+1\right)
\left\Vert a\right\Vert _{L^{\sigma }\left( B_{r}\left( x_{0}\right) \right)
}\right)
\end{equation}%
Using (7.9)\ and (6.12) [Theorem 6] we have%
\begin{equation}
D_{B_{r}\left( x_{0}\right) ,\varepsilon _{l}}\leq D_{B_{r}\left(
x_{0}\right) }
\end{equation}%
where%
\begin{equation}
D_{B_{r}\left( x_{0}\right) }=2^{p-1}L\left( p^{q}+L+2m\left( L+1\right)
\left( D_{p,q}\right) ^{q}+\left( L+1\right) \left\Vert a\right\Vert
_{L^{\sigma }\left( B_{r}\left( x_{0}\right) \right) }\right)
\end{equation}%
then it follows%
\begin{equation*}
\begin{tabular}{l}
$\int\limits_{A_{k,t,y_{0}}^{1,\varepsilon _{l}}}\varepsilon _{l}\left\vert
\nabla u_{\varepsilon _{l}}^{1}\right\vert ^{p}\,+\left\vert \nabla
u_{\varepsilon _{l}}^{1}\right\vert ^{q}\,dx$ \\ 
$\leq \frac{L}{L+1}\int\limits_{A_{k,s,y_{0}}^{1,\varepsilon
_{l}}}\varepsilon _{l}\left\vert \nabla u_{\varepsilon _{l}}^{1}\right\vert
^{p}\,+\left\vert \nabla u_{\varepsilon _{l}}^{1}\right\vert ^{q}\,dx+\frac{%
2^{p}\left( p^{p}+p^{q}\right) }{\left( L+1\right) \left( s-t\right) ^{p}}%
\int\limits_{A_{k,s,y_{0}}^{1,\varepsilon _{l}}}\left( u_{\varepsilon
_{l}}^{1}-k\right) ^{p}\,dx+$ \\ 
$\frac{D_{B_{r}\left( x_{0}\right) }}{L+1}\left( 1+k^{p}\right) \left[ 
\mathcal{L}^{n}\left( A_{k,s,y_{0}}^{1,\varepsilon _{l}}\right) \right] ^{1-%
\frac{p}{n}+\varepsilon }$%
\end{tabular}%
\end{equation*}%
Using Lemma \ref{lemma3} we have the Caccioppoli Inequaity%
\begin{equation*}
\begin{tabular}{l}
$\int\limits_{A_{k,t\varrho y_{0}}^{1,\varepsilon _{l}}}\varepsilon
_{l}\left\vert \nabla u_{\varepsilon _{l}}^{1}\right\vert ^{p}\,+\left\vert
\nabla u_{\varepsilon _{l}}^{1}\right\vert ^{q}\,dx$ \\ 
$\leq \frac{C_{1,B_{r}\left( x_{0}\right) }}{\left( R-\varrho \right) ^{p}}%
\int\limits_{A_{k,R,y_{0}}^{1,\varepsilon _{l}}}\left( u_{\varepsilon
_{l}}^{1}-k\right) ^{p}\,dx+C_{2,B_{r}\left( x_{0}\right) }\left(
1+R^{-\epsilon n}k^{p}\right) \left[ \mathcal{L}^{n}\left(
A_{k,R,y_{0}}^{1,\varepsilon _{l}}\right) \right] ^{1-\frac{p}{n}+\epsilon }$%
\end{tabular}%
\end{equation*}%
Similarly you can proceed for $\alpha =2,...,m$ and we get 
\begin{equation}
\begin{tabular}{l}
$\int\limits_{A_{k,t\varrho y_{0}}^{\alpha ,\varepsilon _{l}}}\varepsilon
_{l}\left\vert \nabla u_{\varepsilon _{l}}^{\alpha }\right\vert
^{p}\,+\left\vert \nabla u_{\varepsilon _{l}}^{\alpha }\right\vert ^{q}\,dx$
\\ 
$\leq \frac{C_{1,B_{r}\left( x_{0}\right) }}{\left( R-\varrho \right) ^{p}}%
\int\limits_{A_{k,R,y_{0}}^{\alpha ,\varepsilon _{l}}}\left( u_{\varepsilon
_{l}}^{\alpha }-k\right) ^{p}\,dx+C_{2,B_{r}\left( x_{0}\right) }\left(
1+R^{-\epsilon n}k^{p}\right) \left[ \mathcal{L}^{n}\left(
A_{k,R,y_{0}}^{\alpha ,\varepsilon _{l}}\right) \right] ^{1-\frac{p}{n}%
+\epsilon }$%
\end{tabular}%
\end{equation}%
for every $\alpha =1,...,m$. Since $-u$ is a minimizer of the int egral
functional 
\begin{equation*}
\mathcal{\tilde{F}}_{0}\left( v,\Omega \right) =\int\limits_{\Omega
}\varepsilon \sum\limits_{\alpha =1}^{m}\left\vert \nabla v^{\alpha
}\right\vert ^{p}+G\left( \sqrt{\left\vert v^{1}\right\vert ^{2}+\varepsilon 
},...,\sqrt{\left\vert v^{m}\right\vert ^{2}+\varepsilon },\sqrt{\left\vert
\nabla v^{1}\right\vert ^{2}+\varepsilon },...,\sqrt{\left\vert \nabla
v^{m}\right\vert ^{2}+\varepsilon }\right) \,dx
\end{equation*}%
then 
\begin{equation}
\begin{tabular}{l}
$\int\limits_{B_{k,t\varrho y_{0}}^{\alpha ,\varepsilon _{l}}}\varepsilon
_{l}\left\vert \nabla u_{\varepsilon _{l}}^{\alpha }\right\vert
^{p}\,+\left\vert \nabla u_{\varepsilon _{l}}^{\alpha }\right\vert ^{q}\,dx$
\\ 
$\leq \frac{C_{1,B_{r}\left( x_{0}\right) }}{\left( R-\varrho \right) ^{p}}%
\int\limits_{B_{k,R,y_{0}}^{\alpha ,\varepsilon _{l}}}\left(
k-u_{\varepsilon _{l}}^{\alpha }\right) ^{p}\,dx+C_{2,B_{r}\left(
x_{0}\right) }\left( 1+R^{-\epsilon n}k^{p}\right) \left[ \mathcal{L}%
^{n}\left( B_{k,R,y_{0}}^{\alpha ,\varepsilon _{l}}\right) \right] ^{1-\frac{%
p}{n}+\epsilon }$%
\end{tabular}%
\end{equation}%
for every $\alpha =1,...,m$. \bigskip Let us observe that 
\begin{equation}
\int\limits_{A_{k,t\varrho y_{0}}^{\alpha ,\varepsilon _{l}}}\varepsilon
_{l}\left\vert \nabla u_{\varepsilon _{l}}^{\alpha }\right\vert
^{p}\,+\left\vert \nabla u_{\varepsilon _{l}}^{\alpha }\right\vert
^{q}\,dx=\int\limits_{B_{r}\left( x_{0}\right) }\left[ \varepsilon
_{l}\left\vert \nabla u_{\varepsilon _{l}}^{\alpha }\right\vert
^{p}\,+\left\vert \nabla u_{\varepsilon _{l}}^{\alpha }\right\vert ^{q}%
\right] 1_{A_{k,t\varrho y_{0}}^{\alpha ,\varepsilon _{l}}}\,dx
\end{equation}%
since, by Holder Inequality and (6.12)\ [Theorem 6],%
\begin{equation}
\begin{tabular}{l}
$\int\limits_{B_{r}\left( x_{0}\right) }\left\vert \nabla u_{\varepsilon
_{l}}^{\alpha }\right\vert ^{q}\left\vert 1_{A_{k,t\varrho y_{0}}^{\alpha
,\varepsilon _{l}}}-1_{A_{k,t\varrho y_{0}}^{\alpha ,0}}\right\vert \,dx$ \\ 
$\leq \left[ \int\limits_{B_{r}\left( x_{0}\right) }\left\vert \nabla
u_{\varepsilon _{l}}^{\alpha }\right\vert ^{p}\,dx\right] ^{\frac{q}{p}}%
\left[ \int\limits_{B_{r}\left( x_{0}\right) }\left\vert 1_{A_{k,t\varrho
y_{0}}^{\alpha ,\varepsilon _{l}}}-1_{A_{k,t\varrho y_{0}}^{\alpha
,0}}\right\vert ^{\frac{p}{p-q}}\,dx\right] ^{\frac{p-q}{p}}$ \\ 
$\leq D_{p,q}^{q}\left[ \int\limits_{B_{r}\left( x_{0}\right) }\left\vert
1_{A_{k,t\varrho y_{0}}^{\alpha ,\varepsilon _{l}}}-1_{A_{k,t\varrho
y_{0}}^{\alpha ,0}}\right\vert ^{\frac{p}{p-q}}\,dx\right] ^{\frac{p-q}{p}}$%
\end{tabular}%
\end{equation}%
where 
\begin{equation}
A_{k,\varrho ,y_{0}}^{1,0}=\left\{ u_{0}^{\alpha }>k\right\} \cap
B_{r}\left( y_{0}\right)
\end{equation}%
and, since, 
\begin{equation}
1_{A_{k,t\varrho y_{0}}^{\alpha ,\varepsilon _{l}}}\rightarrow
1_{A_{k,t\varrho y_{0}}^{\alpha ,0}}\qquad \text{a. e. in }B_{r}\left(
x_{0}\right)
\end{equation}%
and%
\begin{equation}
\left\vert 1_{A_{k,t\varrho y_{0}}^{\alpha ,\varepsilon
_{l}}}-1_{A_{k,t\varrho y_{0}}^{\alpha ,0}}\right\vert ^{\frac{p}{p-q}}\leq
2^{\frac{p}{p-q}}\qquad \text{a. e. in }B_{r}\left( x_{0}\right)
\end{equation}%
then using the Lebesgue Dominated Convergence Theorem [1.34 of [49], see
olso [3, 8, 27]] it follows%
\begin{equation}
\lim\limits_{j\rightarrow +\infty }\int\limits_{B_{r}\left( x_{0}\right)
}\left\vert \nabla u_{\varepsilon _{l}}^{\alpha }\right\vert ^{q}\left\vert
1_{A_{k,t\varrho y_{0}}^{\alpha ,\varepsilon _{l}}}-1_{A_{k,t\varrho
y_{0}}^{\alpha ,0}}\right\vert \,dx=0
\end{equation}%
Moreover by lower semicontinuty we get%
\begin{equation}
\int\limits_{B_{r}\left( x_{0}\right) }\left\vert \nabla u_{0}^{\alpha
}\right\vert ^{q}1_{A_{k,t\varrho y_{0}}^{\alpha ,0}}\,dx\leq
\liminf\limits_{j\rightarrow +\infty }\int\limits_{B_{r}\left( x_{0}\right)
}\left\vert \nabla u_{\varepsilon _{l}}^{\alpha }\right\vert
^{q}1_{A_{k,t\varrho y_{0}}^{\alpha ,0}}\,dx
\end{equation}%
then using (7.14), (7.19) and (7.20)\ we have 
\begin{equation}
\int\limits_{B_{r}\left( x_{0}\right) }\left\vert \nabla u_{0}^{\alpha
}\right\vert ^{q}1_{A_{k,t\varrho y_{0}}^{\alpha ,0}}\,dx\leq
\liminf\limits_{j\rightarrow +\infty }\int\limits_{A_{k,t\varrho
y_{0}}^{\alpha ,\varepsilon _{l}}}\varepsilon _{l}\left\vert \nabla
u_{\varepsilon _{l}}^{\alpha }\right\vert ^{p}\,+\left\vert \nabla
u_{\varepsilon _{l}}^{\alpha }\right\vert ^{q}\,dx
\end{equation}%
Not only, since 
\begin{equation}
\int\limits_{A_{k,R,y_{0}}^{\alpha ,\varepsilon _{l}}}\left( u_{\varepsilon
_{l}}^{\alpha }-k\right) ^{p}\,dx=\int\limits_{B_{r}\left( x_{0}\right)
}\left( u_{\varepsilon _{l}}^{\alpha }-k\right) ^{p}1_{A_{k,t\varrho
y_{0}}^{\alpha ,\varepsilon _{l}}}\,dx
\end{equation}%
then, using (7.22), (7.6), (7.17) and the Lebesgue Dominated Convergence
Theorem [1.34 of [49], see olso [3, 8, 27]] it follows%
\begin{equation}
\lim\limits_{j\rightarrow +\infty }\int\limits_{B_{r}\left( x_{0}\right)
}\left( u_{\varepsilon _{l}}^{\alpha }-k\right) ^{p}1_{A_{k,t\varrho
y_{0}}^{\alpha ,\varepsilon _{l}}}\,dx=\int\limits_{B_{r}\left( x_{0}\right)
}\left( u_{0}^{\alpha }-k\right) ^{p}1_{A_{k,t\varrho y_{0}}^{\alpha ,0}}\,dx
\end{equation}%
and%
\begin{equation}
\lim\limits_{j\rightarrow +\infty }\mathcal{L}^{n}\left(
A_{k,R,y_{0}}^{\alpha ,\varepsilon _{l}}\right) =\mathcal{L}^{n}\left(
A_{k,R,y_{0}}^{\alpha ,0}\right)
\end{equation}%
Using (7.21), (7.23) and (7.24) we obtain%
\begin{equation}
\begin{tabular}{l}
$\int\limits_{A_{k,t\varrho y_{0}}^{\alpha ,0}}\left\vert \nabla
u_{0}^{\alpha }\right\vert ^{q}\,dx$ \\ 
$\leq \frac{C_{1,B_{r}\left( x_{0}\right) }}{\left( R-\varrho \right) ^{p}}%
\int\limits_{A_{k,R,y_{0}}^{\alpha ,0}}\left( u_{0}^{\alpha }-k\right)
^{p}\,dx+C_{2,B_{r}\left( x_{0}\right) }\left( 1+R^{-\epsilon n}k^{p}\right) %
\left[ \mathcal{L}^{n}\left( A_{k,R,y_{0}}^{\alpha ,0}\right) \right] ^{1-%
\frac{p}{n}+\epsilon }$%
\end{tabular}%
\end{equation}%
for every $\alpha =1,...,m$, similary we get 
\begin{equation}
\begin{tabular}{l}
$\int\limits_{B_{k,t\varrho y_{0}}^{\alpha ,0}}\left\vert \nabla
u_{0}^{\alpha }\right\vert ^{q}\,dx$ \\ 
$\leq \frac{C_{1,B_{r}\left( x_{0}\right) }}{\left( R-\varrho \right) ^{p}}%
\int\limits_{B_{k,R,y_{0}}^{\alpha ,0}}\left( k-u_{0}^{\alpha }\right)
^{p}\,dx+C_{2,B_{r}\left( x_{0}\right) }\left( 1+R^{-\epsilon n}k^{p}\right) %
\left[ \mathcal{L}^{n}\left( B_{k,R,y_{0}}^{\alpha ,0}\right) \right] ^{1-%
\frac{p}{n}+\epsilon }$%
\end{tabular}%
\end{equation}%
for every $\alpha =1,...,m$.

\section{Proof of the main Theorem 1}

\bigskip Let us consider $x_{0}\in \Omega $, $0<r<\frac{1}{4}\min \left\{
1,dist\left( \partial \Omega ,x_{0}\right) \right\} $, $y_{0}\in B_{\frac{r}{%
2}}\left( x_{0}\right) \subset B_{r}\left( x_{0}\right) $ and $0<\varrho
<t<s<R<\frac{r}{4}$, since $2\leq q<p<\min \left\{ \frac{2^{\ast }}{2}%
q,q^{\ast }\right\} $, we can choose $0<\epsilon <\frac{\beta q^{\ast }}{2n}$%
, $p=\beta q^{\ast }$ with $\frac{q}{q^{\ast }}<\beta <\min \left\{ \frac{%
2^{\ast }}{2}\frac{q}{q^{\ast }},1,\beta _{0}\right\} $ where $\beta _{0}$\
is a positive real number that we fix later, then the Caccioppoli
Inequalities (7.25) and (7.26) can be rewritten like this%
\begin{equation}
\begin{tabular}{l}
$\int\limits_{A_{k,t\varrho y_{0}}^{\alpha ,0}}\left\vert \nabla
u_{0}^{\alpha }\right\vert ^{q}\,dx$ \\ 
$\leq \frac{C_{1,B_{r}\left( x_{0}\right) }}{\left( R-\varrho \right)
^{\beta q^{\ast }}}\int\limits_{A_{k,R,y_{0}}^{\alpha ,0}}\left(
u_{0}^{\alpha }-k\right) ^{\beta q^{\ast }}\,dx+C_{2,B_{r}\left(
x_{0}\right) }\left( 1+R^{-\epsilon n}k^{\beta q^{\ast }}\right) \left[ 
\mathcal{L}^{n}\left( A_{k,R,y_{0}}^{\alpha ,0}\right) \right] ^{1-\frac{%
\beta q^{\ast }}{n}+\epsilon }$%
\end{tabular}%
\end{equation}%
and%
\begin{equation}
\begin{tabular}{l}
$\int\limits_{B_{k,t\varrho y_{0}}^{\alpha ,0}}\left\vert \nabla
u_{0}^{\alpha }\right\vert ^{q}\,dx$ \\ 
$\leq \frac{C_{1,B_{r}\left( x_{0}\right) }}{\left( R-\varrho \right)
^{\beta q^{\ast }}}\int\limits_{B_{k,R,y_{0}}^{\alpha ,0}}\left(
k-u_{0}^{\alpha }\right) ^{\beta q^{\ast }}\,dx+C_{2,B_{r}\left(
x_{0}\right) }\left( 1+R^{-\epsilon n}k^{\beta q^{\ast }}\right) \left[ 
\mathcal{L}^{n}\left( B_{k,R,y_{0}}^{\alpha ,0}\right) \right] ^{1-\frac{%
\beta q^{\ast }}{n}+\epsilon }$%
\end{tabular}%
\end{equation}%
for every $\alpha =1,...,m$. Let us set $v^{\alpha }\left( x\right)
=u_{0}^{\alpha }\left( x\right) \pm R^{\frac{\epsilon n}{\beta q^{\ast }}}$
and $h=k\pm R^{\frac{\epsilon n}{\beta q^{\ast }}}$\ then the Caccioppoli
Inequalities are written 
\begin{equation}
\begin{tabular}{l}
$\int\limits_{A_{h,t\varrho ,y_{0}}^{\alpha ,0}}\left\vert \nabla v^{\alpha
}\right\vert ^{q}\,dx$ \\ 
$\leq \frac{C_{1,B_{r}\left( x_{0}\right) }}{\left( R-\varrho \right)
^{\beta q^{\ast }}}\int\limits_{A_{h,R,y_{0}}^{\alpha ,0}}\left( v^{\alpha
}-h\right) ^{\beta q^{\ast }}\,dx+2C_{2,B_{r}\left( x_{0}\right) }h^{\beta
q^{\ast }}R^{-\epsilon n}\left[ \mathcal{L}^{n}\left( A_{h,R,y_{0}}^{\alpha
,0}\right) \right] ^{1-\frac{\beta q^{\ast }}{n}+\epsilon }$%
\end{tabular}%
\end{equation}%
and%
\begin{equation}
\begin{tabular}{l}
$\int\limits_{B_{h,t\varrho ,y_{0}}^{\alpha ,0}}\left\vert \nabla v^{\alpha
}\right\vert ^{q}\,dx$ \\ 
$\leq \frac{C_{1,B_{r}\left( x_{0}\right) }}{\left( R-\varrho \right)
^{\beta q^{\ast }}}\int\limits_{B_{h,R,y_{0}}^{\alpha ,0}}\left( h-v^{\alpha
}\right) ^{\beta q^{\ast }}\,dx+2C_{2,B_{r}\left( x_{0}\right) }\left\vert
h\right\vert ^{\beta q^{\ast }}R^{-\epsilon n}\left[ \mathcal{L}^{n}\left(
B_{h,R,y_{0}}^{\alpha ,0}\right) \right] ^{1-\frac{\beta q^{\ast }}{n}%
+\epsilon }$%
\end{tabular}%
\end{equation}%
for every $\alpha =1,...,m$. Moreover let $\frac{1}{2}\leq \varsigma <\tau
\leq 1$\ and $\frac{R}{2}\leq \varsigma R<\tau R\leq R$ then, by (8.3) and
(8.4) it follows 
\begin{equation}
\begin{tabular}{l}
$\int\limits_{A_{h,\varsigma R,y_{0}}^{\alpha ,0}}\left\vert \nabla
v^{\alpha }\right\vert ^{q}\,dx$ \\ 
$\leq \frac{C_{1,B_{r}\left( x_{0}\right) }}{\left( \tau -\varsigma \right)
^{\beta q^{\ast }}}\int\limits_{A_{h,\tau R,y_{0}}^{\alpha ,0}}\left( \frac{%
v^{\alpha }-h}{R}\right) ^{\beta q^{\ast }}\,dx+2C_{2,B_{r}\left(
x_{0}\right) }h^{\beta q^{\ast }}\left( \tau R\right) ^{-\epsilon n}\left[ 
\mathcal{L}^{n}\left( A_{h,\tau R,y_{0}}^{\alpha ,0}\right) \right] ^{1-%
\frac{\beta q^{\ast }}{n}+\epsilon }$%
\end{tabular}%
\end{equation}%
and 
\begin{equation}
\begin{tabular}{l}
$\int\limits_{B_{h,t\varsigma Ry_{0}}^{\alpha ,0}}\left\vert \nabla
v^{\alpha }\right\vert ^{q}\,dx$ \\ 
$\leq \frac{C_{1,B_{r}\left( x_{0}\right) }}{\left( \tau -\varsigma \right)
^{\beta q^{\ast }}}\int\limits_{B_{h,\tau R,y_{0}}^{\alpha ,0}}\left( \frac{%
h-v^{\alpha }}{R}\right) ^{\beta q^{\ast }}\,dx+2C_{2,B_{r}\left(
x_{0}\right) }\left\vert h\right\vert ^{\beta q^{\ast }}\left( \tau R\right)
^{-\epsilon n}\left[ \mathcal{L}^{n}\left( B_{h,\tau R,y_{0}}^{\alpha
,0}\right) \right] ^{1-\frac{\beta q^{\ast }}{n}+\epsilon }$%
\end{tabular}%
\end{equation}%
for every $\alpha =1,...,m$. Now, let us define%
\begin{equation*}
v_{R}\left( x\right) =\frac{v\left( Rx\right) }{R}
\end{equation*}%
since $0<\epsilon <\frac{\beta q^{\ast }}{2n}$, $A_{h,\varsigma
,y_{0}}^{\alpha ,0}=\left\{ x:v_{R}\left( x\right) >\frac{h}{R}\right\} \cap
B_{\varsigma }\left( y_{0}\right) $ and $B_{h,\varsigma ,y_{0}}^{\alpha
,0}=\left\{ x:v_{R}\left( x\right) <\frac{h}{R}\right\} \cap B_{\varsigma
}\left( y_{0}\right) $, then using (8.5) and (8.6), we get

\begin{equation}
\begin{tabular}{l}
$\int\limits_{A_{k,\varsigma ,y_{0}}^{\alpha ,0}}\left\vert \nabla
v_{R}^{\alpha }\right\vert ^{q}\,dx$ \\ 
$\leq \frac{C_{1,B_{r}\left( x_{0}\right) }}{\left( \tau -\varsigma \right)
^{\beta q^{\ast }}}\int\limits_{A_{h,\tau ,y_{0}}^{\alpha ,0}}\left(
v_{R}^{\alpha }-\frac{h}{R}\right) ^{\beta q^{\ast }}\,dx+2^{1+\epsilon
n}C_{2,B_{r}\left( x_{0}\right) }^{\beta q^{\ast }}h^{\beta q^{\ast }}\left[ 
\mathcal{L}^{n}\left( A_{h,\tau ,y_{0}}^{\alpha ,0}\right) \right] ^{1-\frac{%
\beta q^{\ast }}{n}+\epsilon }$%
\end{tabular}%
\end{equation}%
and%
\begin{equation}
\begin{tabular}{l}
$\int\limits_{B_{h,\varsigma ,y_{0}}^{\alpha ,0}}\left\vert \nabla v^{\alpha
}\right\vert ^{q}\,dx$ \\ 
$\leq \frac{C_{1,B_{r}\left( x_{0}\right) }}{\left( \tau -\varsigma \right)
^{\beta q^{\ast }}}\int\limits_{B_{h,\tau ,y_{0}}^{\alpha ,0}}\left( \frac{h%
}{R}-v_{R}^{\alpha }\right) ^{\beta q^{\ast }}\,dx+2^{1+\epsilon
n}C_{2,B_{r}\left( x_{0}\right) }\left\vert h\right\vert ^{\beta q^{\ast }}%
\left[ \mathcal{L}^{n}\left( B_{h,\tau ,y_{0}}^{\alpha ,0}\right) \right]
^{1-\frac{\beta q^{\ast }}{n}+\epsilon }$%
\end{tabular}%
\end{equation}%
for every $\alpha =1,...,m$. \ We define the decreasing sequence%
\begin{equation*}
\varrho _{i}=\frac{1}{2}\left( 1+\frac{1}{2^{i}}\right)
\end{equation*}%
for $i\in 
\mathbb{N}
$, such that $0<\frac{R}{2}\leq \varrho _{i+1}<\varrho _{i}\leq R<2R<R_{0}$.
Moreover let us fix a positive constant $d$, to be chosen lather, and
consider the sequence%
\begin{equation*}
k_{i+1}=k_{i}+\frac{d}{2^{i}}
\end{equation*}%
with $k_{0}=d$ and $i\in 
\mathbb{N}
$. By defining the sequence%
\begin{equation}
W_{\alpha ,i}=\int\limits_{B_{\varrho _{i}}}\left( v_{R}^{\alpha }-\frac{%
k_{i}}{R}\right) _{+}^{\beta q^{\ast }}\,dx=\int\limits_{A_{k_{i},\varrho
_{i}}^{\alpha }}\left( v_{R}^{\alpha }-\frac{k_{i}}{R}\right) ^{\beta
q^{\ast }}\,dx
\end{equation}%
and applying the H\"{o}lder Inequality we get%
\begin{equation}
\begin{tabular}{l}
$W_{\alpha ,i+1}$ \\ 
$=\int\limits_{A_{k_{i+1},\varrho _{i+1}}^{\alpha }}\left( v_{R}^{\alpha }-%
\frac{k_{i+1}}{R}\right) ^{\beta q^{\ast }}\,dx$ \\ 
$\leq \left[ \mathcal{L}^{n}\left( A_{k_{i+1},\varrho _{i+1}}^{\alpha
}\right) \right] ^{1-\beta }\left[ \int\limits_{A_{k_{i+1},\varrho
_{i+1}}^{\alpha }}\left( v_{R}^{\alpha }-\frac{k_{i+1}}{R}\right) ^{q^{\ast
}}\,dx\right] ^{\beta }$%
\end{tabular}%
\end{equation}%
Let $\eta \in C_{0}^{\infty }\left( B_{\tilde{\varrho}_{i}}\right) $ be a
cut-of function satifying%
\begin{equation*}
\begin{tabular}{llllll}
$0\leq \eta \leq 1\quad $in $B_{\tilde{\varrho}_{i}}$ &  & $\eta =1\quad $in 
$B_{\varrho _{i+1}}$ &  & $\left\vert \nabla \eta \right\vert \leq \frac{2}{%
\varrho _{i}-\varrho _{i+1}}\quad $in\ $B_{\tilde{\varrho}_{i}}$ & 
\end{tabular}%
\end{equation*}%
where $\tilde{\varrho}_{i}=\frac{\varrho _{i}+\varrho _{i+1}}{2}$, then%
\begin{equation}
W_{\alpha ,i+1}\leq \left[ \mathcal{L}^{n}\left( A_{k_{i+1},\varrho
_{i+1}}^{\alpha }\right) \right] ^{1-\beta }\left[ \int\limits_{B_{\tilde{%
\varrho}_{i}}}\left[ \eta \left( v_{R}^{\alpha }-\frac{k_{i+1}}{R}\right)
_{+}\right] ^{q^{\ast }}\,dx\right] ^{\beta }
\end{equation}%
Using the Sobolev Inequality it follows%
\begin{equation}
W_{\alpha ,i+1}\leq c_{q,n}\left[ \mathcal{L}^{n}\left( A_{k_{i+1},\varrho
_{i+1}}^{\alpha }\right) \right] ^{1-\beta }\left[ \int\limits_{B_{\tilde{%
\varrho}_{i}}}\left\vert \nabla \left( \eta \left( v_{R}^{\alpha }-\frac{%
k_{i+1}}{R}\right) _{+}\right) \right\vert ^{q}\,dx\right] ^{\frac{q^{\ast
}\beta }{q}}
\end{equation}%
where $c_{q,n}>0$ is a real positive constant depending only on $q$ and $n$.
Since%
\begin{equation}
\begin{tabular}{l}
$\left[ \int\limits_{B_{\tilde{\varrho}_{i}}}\left\vert \nabla \left( \eta
\left( v_{R}^{\alpha }-\frac{k_{i+1}}{R}\right) _{+}\right) \right\vert
^{q}\,dx\right] ^{\frac{q^{\ast }\beta }{q}}$ \\ 
$\leq 2^{\frac{q^{\ast }\beta }{q}-1}\left[ \int\limits_{B_{\tilde{\varrho}%
_{i}}}\left\vert \nabla \eta \right\vert ^{q}\left\vert \left( v_{R}^{\alpha
}-\frac{k_{i+1}}{R}\right) _{+}\right\vert ^{q}\,dx\right] ^{\frac{q^{\ast
}\beta }{q}}$ \\ 
$+2^{\frac{q^{\ast }\beta }{q}-1}\left[ \int\limits_{B_{\tilde{\varrho}%
_{i}}}\eta ^{q}\left\vert \nabla \left( v_{R}^{\alpha }-\frac{k_{i+1}}{R}%
\right) _{+}\right\vert ^{q}\,dx\right] ^{\frac{q^{\ast }\beta }{q}}$%
\end{tabular}%
\end{equation}%
and%
\begin{equation}
\begin{tabular}{l}
$\int\limits_{B_{\tilde{\varrho}_{i}}}\left\vert \nabla \eta \right\vert
^{q}\left\vert \left( v_{R}^{\alpha }-\frac{k_{i+1}}{R}\right)
_{+}\right\vert ^{q}\,dx$ \\ 
$\leq \left( \frac{2}{\varrho _{i}-\varrho _{i+1}}\right)
^{q}\int\limits_{B_{\tilde{\varrho}_{i}}}\left\vert \left( v_{R}^{\alpha }-%
\frac{k_{i+1}}{R}\right) _{+}\right\vert ^{q}\,dx$ \\ 
$\leq \left( \frac{2}{\varrho _{i}-\varrho _{i+1}}\right)
^{q}\int\limits_{A_{k_{i},\varrho _{i}}^{\alpha }}\left\vert v_{R}^{\alpha }-%
\frac{k_{i+1}}{R}\right\vert ^{q}\,dx$ \\ 
$\leq \left( \frac{2}{\varrho _{i}-\varrho _{i+1}}\right)
^{q}\int\limits_{A_{k_{i},\varrho _{i}}^{\alpha }}\left\vert v_{R}^{\alpha }-%
\frac{k_{i}}{R}\right\vert ^{q}\,dx$%
\end{tabular}%
\end{equation}%
usinig (8.12), (8.13) and (8.14) we get%
\begin{equation}
\begin{tabular}{l}
$W_{\alpha ,i+1}$ \\ 
$\leq c_{q,n}2^{\frac{q^{\ast }\beta }{q}-1}\left[ \mathcal{L}^{n}\left(
A_{k_{i+1},\varrho _{i+1}}^{\alpha }\right) \right] ^{1-\beta }2^{q^{\ast
}\beta }\left[ \int\limits_{A_{k_{i},\varrho _{i}}^{\alpha }}\frac{%
\left\vert v_{R}^{\alpha }-\frac{k_{i}}{R}\right\vert ^{q}}{\left( \varrho
_{i}-\varrho _{i+1}\right) ^{q}}\,dx\right] ^{\frac{q^{\ast }\beta }{q}}$ \\ 
$+c_{q,n}2^{\frac{q^{\ast }\beta }{q}-1}\left[ \mathcal{L}^{n}\left(
A_{k_{i+1},\varrho _{i+1}}^{\alpha }\right) \right] ^{1-\beta }\left[
\int\limits_{B_{\tilde{\varrho}_{i}}}\eta ^{q}\left\vert \nabla \left(
v_{R}^{\alpha }-\frac{k_{i+1}}{R}\right) _{+}\right\vert ^{q}\,dx\right] ^{%
\frac{q^{\ast }\beta }{q}}$%
\end{tabular}%
\end{equation}%
Since%
\begin{equation*}
\int\limits_{B_{\tilde{\varrho}_{i}}}\eta ^{q}\left\vert \nabla \left(
v_{R}^{\alpha }-\frac{k_{i+1}}{R}\right) _{+}\right\vert ^{q}\,dx\leq
\int\limits_{A_{k_{i+1},\tilde{\varrho}_{i}}^{\alpha }}\left\vert \nabla
v_{R}^{\alpha }\right\vert ^{q}\,dx
\end{equation*}%
using Caccioppoli Inequality (8.7) it follows%
\begin{equation}
\begin{tabular}{l}
$\int\limits_{B_{\tilde{\varrho}_{i}}}\eta ^{q}\left\vert \nabla \left(
v_{R}^{\alpha }-\frac{k_{i+1}}{R}\right) _{+}\right\vert ^{q}\,dx$ \\ 
$\leq \int\limits_{A_{k_{i+1},\tilde{\varrho}_{i}}^{\alpha }}\left\vert
\nabla v_{R}^{\alpha }\right\vert ^{q}\,dx$ \\ 
$\leq \frac{C_{1}}{\left( \varrho _{i}-\varrho _{i+1}\right) ^{q^{\ast
}\beta }}\int\limits_{A_{k_{i},\varrho _{i}}^{\alpha }}\left\vert
v_{R}^{\alpha }-\frac{k_{i}}{R}\right\vert ^{q^{\ast }\beta
}\,dx+C_{2}\left( 1+\varrho _{i}^{-\epsilon n}k_{i}^{q^{\ast }\beta }\right) %
\left[ \mathcal{L}^{n}\left( A_{k_{i},\varrho _{i}}^{\alpha }\right) \right]
^{1-\frac{q^{\ast }\beta }{n}+\epsilon }$%
\end{tabular}%
\end{equation}%
then by (8.15) and (8.16) we have%
\begin{equation}
\begin{tabular}{l}
$W_{\alpha ,i+1}$ \\ 
$\leq c_{q,n}2^{\frac{q^{\ast }\beta }{q}-1}\left[ \mathcal{L}^{n}\left(
A_{k_{i+1},\varrho _{i+1}}^{\alpha }\right) \right] ^{1-\beta }2^{q^{\ast
}\beta }\left[ \int\limits_{A_{k_{i},\varrho _{i}}^{\alpha }}\frac{%
\left\vert v_{R}^{\alpha }-\frac{k_{i}}{R}\right\vert ^{q}}{\left( \varrho
_{i}-\varrho _{i+1}\right) ^{q}}\,dx\right] ^{\frac{q^{\ast }\beta }{q}}$ \\ 
$+c_{q,n}2^{\frac{q^{\ast }\beta }{q}-1}\left[ \mathcal{L}^{n}\left(
A_{k_{i+1},\varrho _{i+1}}^{\alpha }\right) \right] ^{1-\beta }$ \\ 
$\cdot \left[ \frac{C_{1}}{\left( \varrho _{i}-\varrho _{i+1}\right) ^{\beta
q^{\ast }}}\int\limits_{A_{k_{i},\varrho _{i}}^{\alpha }}\left\vert
v_{R}^{\alpha }-\frac{k_{i}}{R}\right\vert ^{q^{\ast }\beta
}\,dx+C_{2}\left( 1+\varrho _{i}^{-\epsilon n}k_{i}^{q^{\ast }\beta }\right) %
\left[ \mathcal{L}^{n}\left( A_{k_{i},\varrho _{i}}^{\alpha }\right) \right]
^{1-\frac{q^{\ast }\beta }{n}+\epsilon }\right] ^{\frac{q^{\ast }\beta }{q}}$%
\end{tabular}%
\end{equation}%
Moreover by Young Inequality we get%
\begin{equation*}
\int\limits_{A_{k_{i},\varrho _{i}}^{\alpha }}\frac{\left\vert v_{R}^{\alpha
}-\frac{k_{i}}{R}\right\vert ^{q}}{\left( \varrho _{i}-\varrho _{i+1}\right)
^{q}}\,dx\leq \int\limits_{A_{k_{i},\varrho _{i}}^{\alpha }}1+\frac{%
\left\vert v_{R}^{\alpha }-\frac{k_{i}}{R}\right\vert ^{q^{\ast }\beta }}{%
\left( \varrho _{i}-\varrho _{i+1}\right) ^{q^{\ast }\beta }}\,dx
\end{equation*}%
and%
\begin{equation}
\begin{tabular}{l}
$W_{\alpha ,i+1}$ \\ 
$\leq c_{q,n}2^{\frac{q^{\ast }\beta }{q}-1}\left[ \mathcal{L}^{n}\left(
A_{k_{i+1},\varrho _{i+1}}^{\alpha }\right) \right] ^{1-\beta }2^{q^{\ast
}\beta }\left[ \frac{1}{\left( \varrho _{i}-\varrho _{i+1}\right) ^{q^{\ast
}\beta }}\int\limits_{A_{k_{i},\varrho _{i}}^{\alpha }}\left\vert
v_{R}^{\alpha }-\frac{k_{i}}{R}\right\vert ^{q^{\ast }\beta }\,dx+\mathcal{L}%
^{n}\left( A_{k_{i},\varrho _{i}}^{\alpha }\right) \right] ^{\frac{q^{\ast
}\beta }{q}}$ \\ 
$+c_{q,n}2^{\frac{q^{\ast }\beta }{q}-1}\left[ \mathcal{L}^{n}\left(
A_{k_{i+1},\varrho _{i+1}}^{\alpha }\right) \right] ^{1-\beta }$ \\ 
$\cdot \left[ \frac{C_{1}}{\left( \varrho _{i}-\varrho _{i+1}\right)
^{q^{\ast }\beta }}\int\limits_{A_{k_{i},\varrho _{i}}^{\alpha }}\left\vert
v_{R}^{\alpha }-\frac{k_{i}}{R}\right\vert ^{q^{\ast }\beta
}\,dx+C_{2}\left( 1+\varrho _{i}^{-\epsilon n}k_{i}^{q^{\ast }\beta }\right) %
\left[ \mathcal{L}^{n}\left( A_{k_{i},\varrho _{i}}^{\alpha }\right) \right]
^{1-\frac{q^{\ast }\beta }{n}+\epsilon }\right] ^{\frac{q^{\ast }\beta }{q}}$%
\end{tabular}%
\end{equation}%
Moreover it follows%
\begin{equation}
\begin{tabular}{l}
$W_{\alpha ,i+1}$ \\ 
$\leq \frac{c_{q,n}4^{\frac{q^{\ast }\beta }{q}-1}2^{q^{\ast }\beta }}{%
\left( \varrho _{i}-\varrho _{i+1}\right) ^{\frac{\left( q^{\ast }\beta
\right) ^{2}}{q}}}\left[ \mathcal{L}^{n}\left( A_{k_{i+1},\varrho
_{i+1}}^{\alpha }\right) \right] ^{1-\beta }W_{\alpha ,i}^{\frac{q^{\ast
}\beta }{q}}$ \\ 
$+c_{q,n}4^{\frac{q^{\ast }\beta }{q}-1}2^{q^{\ast }\beta }\left[ \mathcal{L}%
^{n}\left( A_{k_{i+1},\varrho _{i+1}}^{\alpha }\right) \right] ^{1-\beta }%
\left[ \mathcal{L}^{n}\left( A_{k_{i},\varrho _{i}}^{\alpha }\right) \right]
^{\frac{q^{\ast }\beta }{q}}$ \\ 
$+\frac{c_{q,n}4^{\frac{q^{\ast }\beta }{q}-1}C_{1}^{\frac{q^{\ast }\beta }{q%
}}}{\left( \varrho _{i}-\varrho _{i+1}\right) ^{\frac{\left( q^{\ast }\beta
\right) ^{2}}{q}}}\left[ \mathcal{L}^{n}\left( A_{k_{i+1},\varrho
_{i+1}}^{\alpha }\right) \right] ^{1-\beta }W_{\alpha ,i}^{\frac{q^{\ast
}\beta }{q}}$ \\ 
$+c_{q,n}4^{\frac{q^{\ast }\beta }{q}-1}C_{2}^{\frac{q^{\ast }\beta }{q}%
}\left( 1+\varrho _{i}^{-\epsilon n}k_{i}^{q^{\ast }\beta }\right) ^{\frac{%
q^{\ast }\beta }{q}}\left[ \mathcal{L}^{n}\left( A_{k_{i+1},\varrho
_{i+1}}^{\alpha }\right) \right] ^{1-\beta }\left[ \mathcal{L}^{n}\left(
A_{k_{i},\varrho _{i}}^{\alpha }\right) \right] ^{\frac{q^{\ast }\beta }{q}-%
\frac{\left( q^{\ast }\beta \right) ^{2}}{nq}+\epsilon \frac{q^{\ast }\beta 
}{q}}$%
\end{tabular}%
\end{equation}%
Since $0<\epsilon <\frac{\beta q^{\ast }}{2n}$, $\frac{1}{2}<\varrho _{i}<1$
and 
\begin{equation*}
\mathcal{L}^{n}\left( A_{k_{i+1},\varrho _{i+1}}^{\alpha }\right) \leq 
\mathcal{L}^{n}\left( A_{k_{i},\varrho _{i}}^{\alpha }\right) \leq \frac{%
W_{\alpha ,i}}{\left( \frac{k_{i+1}-k_{i}}{R}\right) ^{q^{\ast }\beta }}
\end{equation*}%
we get%
\begin{equation}
\begin{tabular}{l}
$W_{\alpha ,i+1}$ \\ 
$\leq \frac{c_{q,n}4^{\frac{q^{\ast }\beta }{q}-1}2^{q^{\ast }\beta }}{%
\left( \varrho _{i}-\varrho _{i+1}\right) ^{\frac{\left( q^{\ast }\beta
\right) ^{2}}{q}}}\left[ \frac{W_{\alpha ,i}}{\left( \frac{k_{i+1}-k_{i}}{R}%
\right) ^{q^{\ast }\beta }}\right] ^{1-\beta -\frac{\left( q^{\ast }\beta
\right) ^{2}}{nq}+\epsilon \frac{q^{\ast }\beta }{q}}W_{\alpha ,i}^{\frac{%
q^{\ast }\beta }{q}}$ \\ 
$+c_{q,n}4^{\frac{q^{\ast }\beta }{q}-1}2^{q^{\ast }\beta }\left[ \frac{%
W_{\alpha ,i}}{\left( \frac{k_{i+1}-k_{i}}{R}\right) ^{q^{\ast }\beta }}%
\right] ^{1-\beta -\frac{\left( q^{\ast }\beta \right) ^{2}}{nq}+\epsilon 
\frac{q^{\ast }\beta }{q}}\left[ \frac{W_{\alpha ,i}}{\left( \frac{%
k_{i+1}-k_{i}}{R}\right) ^{q^{\ast }\beta }}\right] ^{\frac{q^{\ast }\beta }{%
q}}$ \\ 
$+\frac{c_{q,n}4^{\frac{q^{\ast }\beta }{q}-1}C_{1}^{\frac{q^{\ast }\beta }{q%
}}}{\left( \varrho _{i}-\varrho _{i+1}\right) ^{\frac{\left( q^{\ast }\beta
\right) ^{2}}{q}}}\left[ \frac{W_{\alpha ,i}}{\left( \frac{k_{i+1}-k_{i}}{R}%
\right) ^{q^{\ast }\beta }}\right] ^{1-\beta -\frac{\left( q^{\ast }\beta
\right) ^{2}}{nq}+\epsilon \frac{q^{\ast }\beta }{q}}W_{\alpha ,i}^{\frac{%
q^{\ast }\beta }{q}}$ \\ 
$+c_{q,n}4^{\frac{q^{\ast }\beta }{q}-1}C_{2}^{\frac{q^{\ast }\beta }{q}%
}\left( 1+2^{\epsilon ni}k_{i}^{q^{\ast }\beta }\right) ^{\frac{q^{\ast
}\beta }{q}}\left[ \frac{W_{\alpha ,i}}{\left( \frac{k_{i+1}-k_{i}}{R}%
\right) ^{q^{\ast }\beta }}\right] ^{1-\beta }\left[ \frac{W_{\alpha ,i}}{%
\left( \frac{k_{i+1}-k_{i}}{R}\right) ^{q^{\ast }\beta }}\right] ^{\frac{%
q^{\ast }\beta }{q}-\frac{\left( q^{\ast }\beta \right) ^{2}}{nq}+\epsilon 
\frac{q^{\ast }\beta }{q}}$%
\end{tabular}%
\end{equation}%
Since, $\varrho _{i}-\varrho _{i+1}=\frac{1}{2^{1+i}}$ and $k_{i+1}-k_{i}=%
\frac{d}{2^{i}}$, it folows%
\begin{equation}
\begin{tabular}{l}
$W_{\alpha ,i+1}$ \\ 
$\leq \frac{c_{q,n}4^{\frac{q^{\ast }\beta }{q}-1}2^{q^{\ast }\beta }}{%
\left( \frac{1}{2^{1+i}}\right) ^{\frac{\left( q^{\ast }\beta \right) ^{2}}{q%
}}}\frac{W_{\alpha ,i}^{1-\beta +\frac{q^{\ast }\beta }{q}-\frac{\left(
q^{\ast }\beta \right) ^{2}}{nq}+\epsilon \frac{q^{\ast }\beta }{q}}}{\left( 
\frac{d}{2^{i}R}\right) ^{q^{\ast }\beta \left( 1-\beta -\frac{\left(
q^{\ast }\beta \right) ^{2}}{nq}+\epsilon \frac{q^{\ast }\beta }{q}\right) }}
$ \\ 
$+c_{q,n}4^{\frac{q^{\ast }\beta }{q}-1}2^{q^{\ast }\beta }\frac{W_{\alpha
,i}^{1-\beta +\frac{q^{\ast }\beta }{q}-\frac{\left( q^{\ast }\beta \right)
^{2}}{nq}+\epsilon \frac{q^{\ast }\beta }{q}}}{\left( \frac{d}{2^{i}R}%
\right) ^{q^{\ast }\beta \left( 1-\beta +\frac{q^{\ast }\beta }{q}-\frac{%
\left( q^{\ast }\beta \right) ^{2}}{nq}+\epsilon \frac{q^{\ast }\beta }{q}%
\right) }}$ \\ 
$+\frac{c_{q,n}4^{\frac{q^{\ast }\beta }{q}-1}C_{1}^{\frac{q^{\ast }\beta }{q%
}}}{\left( \frac{1}{2^{1+i}}\right) ^{\frac{\left( q^{\ast }\beta \right)
^{2}}{q}}}\frac{W_{\alpha ,i}^{1-\beta +\frac{q^{\ast }\beta }{q}-\frac{%
\left( q^{\ast }\beta \right) ^{2}}{nq}+\epsilon \frac{q^{\ast }\beta }{q}}}{%
\left( \frac{d}{2^{i}R}\right) ^{q^{\ast }\beta \left( 1-\beta -\frac{\left(
q^{\ast }\beta \right) ^{2}}{nq}+\epsilon \frac{q^{\ast }\beta }{q}\right) }}
$ \\ 
$+c_{q,n}4^{\frac{q^{\ast }\beta }{q}-1}C_{2}^{\frac{q^{\ast }\beta }{q}%
}\left( 1+2^{\epsilon ni}k_{i}^{q^{\ast }\beta }\right) ^{\frac{q^{\ast
}\beta }{q}}\frac{W_{\alpha ,i}^{1-\beta +\frac{q^{\ast }\beta }{q}-\frac{%
\left( q^{\ast }\beta \right) ^{2}}{nq}+\epsilon \frac{q^{\ast }\beta }{q}}}{%
\left( \frac{d}{2^{i}R}\right) ^{q^{\ast }\beta \left( 1-\beta +\frac{%
q^{\ast }\beta }{q}-\frac{\left( q^{\ast }\beta \right) ^{2}}{nq}+\epsilon 
\frac{q^{\ast }\beta }{q}\right) }}$%
\end{tabular}%
\end{equation}

Moreover $0<\epsilon <\frac{\beta q^{\ast }}{2n}$ and $k_{i}=d\sum%
\limits_{j=0}^{i}2^{-j}\leq 2d$, then%
\begin{equation}
\begin{tabular}{l}
$W_{\alpha ,i+1}$ \\ 
$\leq \frac{c_{q,n}4^{\frac{q^{\ast }\beta }{q}-1}\left( 2^{q^{\ast }\beta
}+C_{1}^{\frac{q^{\ast }\beta }{q}}\right) }{\left( \frac{1}{2^{1+i}}\right)
^{\frac{\left( q^{\ast }\beta \right) ^{2}}{q}}}\frac{W_{\alpha ,i}^{1-\beta
+\frac{q^{\ast }\beta }{q}-\frac{\left( q^{\ast }\beta \right) ^{2}}{nq}%
+\epsilon \frac{q^{\ast }\beta }{q}}}{\left( \frac{d}{2^{i}R}\right)
^{q^{\ast }\beta \left( 1-\beta -\frac{\left( q^{\ast }\beta \right) ^{2}}{nq%
}+\epsilon \frac{q^{\ast }\beta }{q}\right) }}$ \\ 
$+c_{q,n}4^{\frac{q^{\ast }\beta }{q}-1}\left[ 2^{q^{\ast }\beta }+C_{2}^{%
\frac{q^{\ast }\beta }{q}}\left( 1+2^{q^{\ast }\beta i}d^{q^{\ast }\beta
}\right) ^{\frac{q^{\ast }\beta }{q}}\right] \frac{W_{\alpha ,i}^{1-\beta +%
\frac{q^{\ast }\beta }{q}-\frac{\left( q^{\ast }\beta \right) ^{2}}{nq}%
+\epsilon \frac{q^{\ast }\beta }{q}}}{\left( \frac{d}{2^{i}R}\right)
^{q^{\ast }\beta \left( 1-\beta +\frac{q^{\ast }\beta }{q}-\frac{\left(
q^{\ast }\beta \right) ^{2}}{nq}+\epsilon \frac{q^{\ast }\beta }{q}\right) }}
$%
\end{tabular}%
\end{equation}%
if we choose $d>1$ we get%
\begin{equation}
\begin{tabular}{l}
$W_{\alpha ,i+1}$ \\ 
$\leq \frac{c_{q,n}4^{\frac{q^{\ast }\beta }{q}-1}\left( 2^{q^{\ast }\beta
}+C_{1}^{\frac{q^{\ast }\beta }{q}}\right) }{\left( \frac{1}{2^{1+i}}\right)
^{\frac{\left( q^{\ast }\beta \right) ^{2}}{q}}}\frac{W_{\alpha ,i}^{1-\beta
+\frac{q^{\ast }\beta }{q}-\frac{\left( q^{\ast }\beta \right) ^{2}}{nq}%
+\epsilon \frac{q^{\ast }\beta }{q}}}{\left( \frac{d}{2^{i}R}\right)
^{q^{\ast }\beta \left( 1-\beta -\frac{\left( q^{\ast }\beta \right) ^{2}}{nq%
}+\epsilon \frac{q^{\ast }\beta }{q}\right) }}$ \\ 
$+c_{q,n}4^{\frac{q^{\ast }\beta }{q}-1}\left[ 2^{q^{\ast }\beta }+d^{\frac{%
\left( q^{\ast }\beta \right) ^{2}}{q}}C_{2}^{\frac{q^{\ast }\beta }{q}%
}\left( 1+2^{q^{\ast }\beta i}\right) ^{\frac{q^{\ast }\beta }{q}}\right] 
\frac{W_{\alpha ,i}^{1-\beta +\frac{q^{\ast }\beta }{q}-\frac{\left( q^{\ast
}\beta \right) ^{2}}{nq}+\epsilon \frac{q^{\ast }\beta }{q}}}{\left( \frac{d%
}{2^{i}R}\right) ^{q^{\ast }\beta \left( 1-\beta +\frac{q^{\ast }\beta }{q}-%
\frac{\left( q^{\ast }\beta \right) ^{2}}{nq}+\epsilon \frac{q^{\ast }\beta 
}{q}\right) }}$%
\end{tabular}%
\end{equation}%
and%
\begin{equation}
\begin{tabular}{l}
$W_{\alpha ,i+1}$ \\ 
$\leq C_{1,q,n,\beta }\,2^{\left( 1+i\right) \frac{\left( q^{\ast }\beta
\right) ^{2}}{q}}\frac{W_{\alpha ,i}^{1-\beta +\frac{q^{\ast }\beta }{q}-%
\frac{\left( q^{\ast }\beta \right) ^{2}}{nq}+\epsilon \frac{q^{\ast }\beta 
}{q}}}{\left( \frac{d}{2^{i}R}\right) ^{q^{\ast }\beta \left( 1-\beta -\frac{%
\left( q^{\ast }\beta \right) ^{2}}{nq}+\epsilon \frac{q^{\ast }\beta }{q}%
\right) }}$ \\ 
$+C_{2,q,n,\beta }\,d^{\frac{\left( q^{\ast }\beta \right) ^{2}}{q}%
}2^{\left( 1+i\right) \frac{\left( q^{\ast }\beta \right) ^{2}}{q}}\frac{%
W_{\alpha ,i}^{1-\beta +\frac{q^{\ast }\beta }{q}-\frac{\left( q^{\ast
}\beta \right) ^{2}}{nq}+\epsilon \frac{q^{\ast }\beta }{q}}}{\left( \frac{d%
}{2^{i}R}\right) ^{q^{\ast }\beta \left( 1-\beta +\frac{q^{\ast }\beta }{q}-%
\frac{\left( q^{\ast }\beta \right) ^{2}}{nq}+\epsilon \frac{q^{\ast }\beta 
}{q}\right) }}$%
\end{tabular}%
\end{equation}%
\bigskip where $C_{1,q,n,\beta }=c_{q,n}4^{\frac{q^{\ast }\beta }{q}%
-1}\left( 2^{q^{\ast }\beta }+C_{1}^{\frac{q^{\ast }\beta }{q}}\right) $ and 
$C_{2,q,n,\beta }=c_{q,n}4^{\frac{q^{\ast }\beta }{q}-1}\left( 1+C_{2}^{%
\frac{q^{\ast }\beta }{q}}\right) 2^{\frac{q^{\ast }\beta }{q}}$. Since $%
0<R<r<1$ it follows%
\begin{equation}
W_{\alpha ,i+1}\leq \frac{C_{3,q,n,\beta }\,\left[ 4^{q^{\ast }\beta \left( 
\frac{q^{\ast }\beta }{q}+1-\beta -\frac{\left( q^{\ast }\beta \right) ^{2}}{%
nq}+\epsilon \frac{q^{\ast }\beta }{q}\right) }\right] ^{i}}{\left( \frac{d}{%
R}\right) ^{q^{\ast }\beta \left( 1-\beta -\frac{\left( q^{\ast }\beta
\right) ^{2}}{nq}+\epsilon \frac{q^{\ast }\beta }{q}\right) }}W_{\alpha
,i}^{1-\beta +\frac{q^{\ast }\beta }{q}-\frac{\left( q^{\ast }\beta \right)
^{2}}{nq}+\epsilon \frac{q^{\ast }\beta }{q}}
\end{equation}%
where $C_{3,q,n,\beta }=\left( C_{1,q,n,\beta }+C_{2,q,n,\beta }\right)
4^{q^{\ast }\beta \left( \frac{q^{\ast }\beta }{q}+1-\beta -\frac{\left(
q^{\ast }\beta \right) ^{2}}{nq}+\epsilon \frac{q^{\ast }\beta }{q}\right) }$%
. Let us consider the inequality $1-\beta -\frac{\left( q^{\ast }\beta
\right) ^{2}}{nq}+\epsilon \frac{q^{\ast }\beta }{q}>0$ that it is
equivalent to the inequality 
\begin{equation*}
\frac{\left( q^{\ast }\right) ^{2}}{nq}\beta ^{2}+\left( 1-\epsilon \frac{%
q^{\ast }}{q}\right) \beta -1<0
\end{equation*}%
which has as a set of solution $\beta _{1}<\beta <\beta _{2}$ where $\beta
_{1}=\frac{-\left( 1-\epsilon \frac{q^{\ast }}{q}\right) -\sqrt{\left(
1-\epsilon \frac{q^{\ast }}{q}\right) ^{2}+4\frac{\left( q^{\ast }\right)
^{2}}{nq}}}{2\frac{\left( q^{\ast }\right) ^{2}}{nq}}$ and\ $\beta _{2}=%
\frac{-\left( 1-\epsilon \frac{q^{\ast }}{q}\right) +\sqrt{\left( 1-\epsilon 
\frac{q^{\ast }}{q}\right) ^{2}+4\frac{\left( q^{\ast }\right) ^{2}}{nq}}}{2%
\frac{\left( q^{\ast }\right) ^{2}}{nq}}$; since $\frac{q}{q^{\ast }}<\beta
<\min \left\{ \frac{2^{\ast }}{2}\frac{q}{q^{\ast }},1,\beta _{0}\right\} $
then after simple algebraic calculation we have 
\begin{equation*}
\frac{q}{q^{\ast }}<\frac{-\left( 1-\epsilon \frac{q^{\ast }}{q}\right) +%
\sqrt{\left( 1-\epsilon \frac{q^{\ast }}{q}\right) ^{2}+4\frac{\left(
q^{\ast }\right) ^{2}}{nq}}}{2\frac{\left( q^{\ast }\right) ^{2}}{nq}}
\end{equation*}%
We put $\beta _{0}=$ $\frac{-\left( 1-\epsilon \frac{q^{\ast }}{q}\right) +%
\sqrt{\left( 1-\epsilon \frac{q^{\ast }}{q}\right) ^{2}+4\frac{\left(
q^{\ast }\right) ^{2}}{nq}}}{2\frac{\left( q^{\ast }\right) ^{2}}{nq}}$ and
we choose%
\begin{equation*}
d=C_{4,q,n,\beta }\left[ \frac{1}{R^{n+\frac{\frac{\left( q^{\ast }\beta
\right) ^{2}}{q}-q^{\ast }\beta }{\frac{q^{\ast }\beta }{q}-\beta -\frac{%
\left( q^{\ast }\beta \right) ^{2}}{nq}+\epsilon \frac{q^{\ast }\beta }{q}}}}%
\int\limits_{B_{R}\left( y_{0}\right) }u_{0,+}^{q^{\ast }\beta }\,dx\right]
^{\frac{\frac{q^{\ast }\beta }{q}-\beta -\frac{\left( q^{\ast }\beta \right)
^{2}}{nq}+\epsilon \frac{q^{\ast }\beta }{q}}{q^{\ast }\beta \left[ 1-\beta -%
\frac{\left( q^{\ast }\beta \right) ^{2}}{nq}+\epsilon \frac{q^{\ast }\beta 
}{q}\right] }}+1
\end{equation*}%
where $C_{4,q,n,\beta }=C_{3,q,n,\beta }^{\frac{1}{q^{\ast }\beta \left[
1-\beta -\frac{\left( q^{\ast }\beta \right) ^{2}}{nq}+\epsilon \frac{%
q^{\ast }\beta }{q}\right] }}4^{\frac{q^{\ast }\beta \left[ \frac{q^{\ast
}\beta }{q}+1-\beta -\frac{\left( q^{\ast }\beta \right) ^{2}}{nq}+\epsilon 
\frac{q^{\ast }\beta }{q}\right] }{\left[ \frac{q^{\ast }\beta }{q}-\beta -%
\frac{\left( q^{\ast }\beta \right) ^{2}}{nq}+\epsilon \frac{q^{\ast }\beta 
}{q}\right] \left[ 1-\beta -\frac{\left( q^{\ast }\beta \right) ^{2}}{nq}%
+\epsilon \frac{q^{\ast }\beta }{q}\right] }}$and $u_{0,+}=\max \left\{
u_{0},0\right\} $, then using Lemma 4 we get%
\begin{equation}
ess-\sup_{B_{\frac{R}{2}}\left( y_{0}\right) }\left( u_{0}^{\alpha }\right)
<C_{4,q,n,\beta }\left[ \frac{1}{R^{n+\frac{\frac{\left( q^{\ast }\beta
\right) ^{2}}{q}-q^{\ast }\beta }{\frac{q^{\ast }\beta }{q}-\beta -\frac{%
\left( q^{\ast }\beta \right) ^{2}}{nq}+\epsilon \frac{q^{\ast }\beta }{q}}}}%
\int\limits_{B_{R}\left( y_{0}\right) }u_{0,+}^{q^{\ast }\beta }\,dx\right]
^{\frac{\frac{q^{\ast }\beta }{q}-\beta -\frac{\left( q^{\ast }\beta \right)
^{2}}{nq}+\epsilon \frac{q^{\ast }\beta }{q}}{q^{\ast }\beta \left[ 1-\beta -%
\frac{\left( q^{\ast }\beta \right) ^{2}}{nq}+\epsilon \frac{q^{\ast }\beta 
}{q}\right] }}+1
\end{equation}%
for every $\alpha =1,...,m$. \bigskip Similary we get%
\begin{equation}
ess-\inf_{B_{\frac{R}{2}}\left( y_{0}\right) }\left( u_{0}^{\alpha }\right)
>-C_{4,q,n,\beta }\left[ \frac{1}{R^{n+\frac{\frac{\left( q^{\ast }\beta
\right) ^{2}}{q}-q^{\ast }\beta }{\frac{q^{\ast }\beta }{q}-\beta -\frac{%
\left( q^{\ast }\beta \right) ^{2}}{nq}+\epsilon \frac{q^{\ast }\beta }{q}}}}%
\int\limits_{B_{R}\left( y_{0}\right) }u_{0,-}^{q^{\ast }\beta }\,dx\right]
^{\frac{\frac{q^{\ast }\beta }{q}-\beta -\frac{\left( q^{\ast }\beta \right)
^{2}}{nq}+\epsilon \frac{q^{\ast }\beta }{q}}{q^{\ast }\beta \left[ 1-\beta -%
\frac{\left( q^{\ast }\beta \right) ^{2}}{nq}+\epsilon \frac{q^{\ast }\beta 
}{q}\right] }}-1
\end{equation}%
Prooceding as in Thorem 7.4 of [27] we obtain $u_{0}^{\alpha }\in L^{\infty
}\left( B_{\frac{r}{2}}\left( x_{0}\right) \right) $ for every $\alpha
=1,...,m$, then it follows $u_{0}\in L_{loc}^{\infty }\left( \Omega ,%
\mathbb{R}
^{m}\right) $ or $\left\vert u_{0}\right\vert \in L_{loc}^{\infty }\left(
\Omega \right) $.

\section{\protect\bigskip Declarations:}

Data sharing not applicable to this article as no datasets were generated or
analysed during the current study. The author has no conicts of interest to
declare that are relevant to the content of this article.

I would like to thank all my family and friends for the support given to me
over the years: Elisa Cirri, Caterina Granucci, Delia Granucci, Irene
Granucci, Laura and Fiorenza Granucci, Massimo Masi and Monia Randolfi.

\end{document}